\documentclass[a4paper,10pt]{article}
\usepackage{stmaryrd}
\usepackage{amsfonts}
\usepackage{bbm}
\usepackage{amscd}
\usepackage{mathrsfs}
\usepackage{latexsym,amssymb,amsmath,amscd,amscd,amsthm,amsxtra}
\usepackage[dvips]{graphicx}
\usepackage[utf8]{inputenc}
\usepackage[T1]{fontenc}
\usepackage{lmodern}
\usepackage{amssymb}
\usepackage[all]{xy}
\usepackage{nicefrac,mathtools,enumitem}
\usepackage{microtype}

\usepackage{ifpdf}
\ifpdf
 \usepackage[colorlinks,final,backref=page,hyperindex]{hyperref}
\else
 \usepackage[colorlinks,final,backref=page,hyperindex,hypertex]{hyperref}
\fi
\newtheorem{thm}{Theorem}[section]
\newtheorem{lem}[thm]{Lemma}
\newtheorem{cor}[thm]{Corollary}
\newtheorem{pro}[thm]{Proposition}
\newtheorem{ex}[thm]{Example}
\newtheorem{rmk}[thm]{Remark}
\newtheorem{defi}[thm]{Definition}

\newcommand {\emptycomment}[1]{} 
\newcommand{\yh}[1]{\textcolor{blue}{ #1}}
\newcommand{\liu}[1]{\textcolor{red}{ #1}}
\usepackage{color}

\newcommand{\lon }{\,\rightarrow\,}
\newcommand{\be }{\begin{equation}}
\newcommand{\ee }{\end{equation}}

\newcommand{\pf}{\noindent{\bf Proof.}\ }

\newcommand{\bl}{\bar{l} }
\newcommand{\tl}{\tilde{\frkl} }
\newcommand{\lhg}{l^\huaG }

\newcommand{\Real}{\mathbb R}

\newcommand{\huaB}{\mathcal{B}}

\newcommand{\huaA}{\mathcal{A}}
\newcommand{\huaL}{\mathcal{L}}

\newcommand{\huaE}{\mathcal{E}}
\newcommand{\huaF}{\mathcal{F}}
\newcommand{\huaG}{\mathcal{G}}
\newcommand{\huaM}{\mathcal{M}}

\newcommand{\huaV}{\mathcal{V}}

\newcommand{\huaD}{\mathcal{D}}

\newcommand{\CWM}{C^{\infty}(M)}

\newcommand{\frka}{\mathfrak a}

\newcommand{\frkg}{\mathfrak g}

\newcommand{\frkl}{\mathfrak l}

\newcommand{\frkD}{\mathfrak D}

\newcommand{\frkL}{\mathfrak L}

\newcommand{\frkX}{\mathfrak X}

\def\qed{\hfill ~\vrule height6pt width6pt depth0pt}

\newcommand{\half}{\frac{1}{2}}

\newcommand{\Sym}{\mathsf{Sym}}

\newcommand{\Courant}[1]{\left\llbracket  #1\right\rrbracket }
\newcommand{\Poisson}[1]{\{ #1\}}
\newcommand{\Schouten}[1]{   [    #1   ]_{S}   }
\newcommand{\br}[1]{   [ \cdot,    \cdot  ]   }
\newcommand{\degree}[1]{   \mid    #1  \mid  }

\newcommand{\id}{\rm{id}}

\newcommand{\g}{\mathfrak g}

\newcommand{\ii}{\mathbbm{i}}

\newcommand{\dM}{\mathrm{d}}

\newcommand{\LWX}{\mathrm{CLWX}}
\newcommand{\Hom}{\mathrm{Hom}}

\newcommand{\ad}{\mathrm{ad}}

\newcommand{\pr}{\mathrm{pr}}
\newcommand{\Jac}{\mathrm{Jac}}
\newcommand{\bas}{\mathrm{bas}}

\newcommand{\img}{\mathrm{im}}

\newcommand{\R}{\mathbb{R}}

\begin{document}
\title{Homotopy Poisson algebras, Maurer-Cartan elements and Dirac structures of $\LWX$ 2-algebroids\thanks
 {
Research supported by NSFC (11922110,11901501), NSF of Jilin Province (20170101050JC), Nanhu Scholars Program for Young Scholars of XYNU.
 }
 }
\author{Jiefeng Liu and Yunhe Sheng
}

\date{}
\footnotetext{{Keyword}: $\LWX$ $2$-algebroid,  split Lie $2$-algebroid, split Lie $2$-bialgebroid, Manin triple, weak Dirac structure, Maurer-Cartan element, homotopy Poisson algebra }
\footnotetext{{{MSC}}:16T10, 17B63, 53D17.}
\maketitle
\begin{abstract}
In this paper, we construct a homotopy Poisson algebra of degree 3 associated to a split Lie 2-algebroid, by which we give a new approach to characterize a split Lie 2-bialgebroid. We develop the differential calculus associated to a split Lie 2-algebroid and establish the Manin triple theory for split Lie 2-algebroids. More precisely, we give the notion of a strict Dirac structure and define a Manin triple for split Lie 2-algebroids to be a $\LWX$ 2-algebroid with two transversal strict Dirac structures. We show that there is a one-to-one correspondence between Manin triples for split Lie 2-algebroids and split Lie 2-bialgebroids. We further introduce the notion of a weak Dirac structure of a $\LWX$ 2-algebroid and show that the graph of a Maurer-Cartan element of the homotopy Poisson algebra of degree 3  associated to a split Lie 2-bialgebroid is a weak Dirac structure. Various examples including the string Lie 2-algebra, split Lie 2-algebroids constructed from integrable distributions and  left-symmetric algebroids are given.
\end{abstract}
\tableofcontents

\section{Introduction}

The notion of a Lie algebroid was introduced by Pradines in 1967, which
is a generalization of Lie algebras and tangent bundles. Just as Lie algebras are the
infinitesimal objects of Lie groups, Lie algebroids are the infinitesimal objects of Lie groupoids. See \cite{General theory of Lie groupoid and Lie algebroid} for the general theory about Lie algebroids. The notion of a Lie bialgebroid was introduced by Mackenzie and Xu in \cite{MackenzieX:1994} as the infinitesimal of a Poisson groupoid. To study the double of a Lie bialgebroid, Liu, Weinstein and Xu introduced the notion of a Courant algebroid in \cite{lwx} and established the Manin triple theory for Lie algebroids. There are many applications of Courant algebroids. See \cite{BCG,DLB,gualtieri,hitchin,Schwarzbach4,Roytenbergphdthesis,rw,royt,RoyCF}  for more details. In particular, the relation between Dirac structures and Maurer-Cartan elements were studied in detail in \cite{lwx,LiuDirac}. The notion of a Dirac structure was originally introduced by Courant in \cite{CourantDirac} to unify symplectic structures and Poisson structures. Then it was widely studied and has many applications, e.g. in generalized complex geometry \cite{gualtieri,hitchin}, in the theory of D-branes for the Wess–Zumino–Witten model \cite{CGM}, in moment map theories \cite{BC}, in Poisson geometry \cite{LiM,Mei18,Mei17}, in Dixmier-Douady bundles \cite{AM}, in reduction theory \cite{Jotz13}, in $L_\infty$-algebras \cite{zambon:l-infty} and in integrable systems \cite{CGdS}. See \cite{Bur} for more details.

Recently, people have paid more attention to higher categorical
structures by reasons in both mathematics and physics.  A Lie $2$-algebra
is the categorification of a Lie algebra \cite{baez:2algebras}. The 2-category of Lie 2-algebras is equivalent to the 2-category of 2-term $L_\infty$-algebras. Due to this reason, an $n$-term $L_\infty$-algebra will be called a Lie $n$-algebra. See \cite{LadaMartin,stasheff:introductionSHLA,Stasheff1} for more details of $L_\infty$-algebras. Usually an NQ-manifold of degree $n$ is considered as a Lie $n$-algebroid \cite{Voronov:2010halgd}. In \cite{sz}, a split Lie $n$-algebroid is defined  using the language of graded vector bundles.   The equivalence between  the category of  split Lie $n$-algebroids and the category of NQ-manifolds of degree $n$ was given in \cite{BP}.

 The notion of a $\LWX$ 2-algebroid was introduced in \cite{LiuSheng} as the categorification of a Courant algebroid. There is a one-to-one correspondence between $\LWX$ 2-algebroids and symplectic NQ-manifolds of degree 3, and the later can be used to construct 4D topological field theory \cite{aksz,Ikeda}. There is also a close connection between $\LWX$ 2-algebroids and the first Pontryagin classes of   quadratic Lie 2-algebroids introduced in \cite{sheng}, which are represented by closed 5-forms.
 The notion of a split Lie 2-bialgebroid was also introduced in \cite{LiuSheng} and it is shown that the double of a split Lie 2-bialgebroid is a $\LWX$ 2-algebroid.     The split Lie 2-bialgebroid used here is a direct geometric generalization of the Lie 2-bialgebra introduced in \cite{BSZ,CSX1}. See \cite{olga} and \cite{BV} for the more general notions of an $L_\infty$-bialgebra and an $L_\infty$-bialgebroid.

The first purpose of this paper is to establish the Manin triple theory for split Lie 2-algebroids, i.e. to show that two transversal Dirac structures of a $\LWX$ 2-algebroid  constitute a split Lie 2-bialgebroid.  The second purpose of this paper is to study the   homotopy Poisson algebra associated to a split Lie 2-algebroid, and establish the relation between its Maurer-Cartan elements  and Dirac structures of the $\LWX$ 2-algebroid.  Upon careful study, we  found  that we need to take two kinds of Dirac structures into account, one is called a strict Dirac structure, served for the first purpose, and the other is called a weak Dirac structure, served for the second purpose.  It is the weak Dirac structures that reflect properties of higher structures and make the paper containing more meaningful contents. Due to aforementioned important applications of Dirac geometry, it is natural to explore similar applications of Dirac structures introduced in this paper, which will be studied in the future.

To study the Manin triple theory for split Lie 2-algebroids, we develop the differential calculus for split Lie 2-algebroids in  Section \ref{sec:D}. In particular, we define the coboundary operator, Lie derivatives, the contraction operator and give their properties.

In Section \ref{sec:H}, first we define   bracket operations  $[\cdot]_S,~[\cdot,\cdot]_S,~[\cdot,\cdot,\cdot]_S$ on $\Sym(\huaA[-3])$ associated to a split Lie 2-algebroid $\huaA$ using the derived bracket \cite{getzler:higher-derived,Kosmann-Schwarzbach,VoronovHigherP}, and show that $(\Sym(\huaA[-3]),[\cdot]_S,~[\cdot,\cdot]_S,~[\cdot,\cdot,\cdot]_S)$  is a homotopy Poisson algebra of degree 3. The notion of a homotopy Poisson manifold of degree $n$ was given in \cite{LSX} in the study of the dual structure of a Lie 2-algebra. See also \cite{derivedPoisson,Bruce,CF,VoronovHigherP,Mehta,Voronov1} for more applications of similar structures. Then we use the usual differential geometry language to give a new characterization of a split Lie 2-bialgebroid.

In Section \ref{sec:M}, we introduce the notion of a strict Dirac structure of a $\LWX$ 2-algebroid and establish the Manin triple theory for split Lie 2-algebroids. More precisely, we show that there is a one-to-one correspondence between Manin triples  of split Lie $2$-algebroids and split Lie $2$-bialgebroids.  Note  that a strict Dirac structure of a $\LWX$ 2-algebroid is defined to be a maximal isotropic graded subbundle whose section space is closed under the multiplication.

In Section \ref{sec:MC}, first we introduce the notion of a weak Dirac structure of a $\LWX$ 2-algebroid, which is a Lie 2-algebroid such that there is a Leibniz 2-algebra morphism from the underlying Lie 2-algebra  to the underlying Leibniz 2-algebra of the original $\LWX$ 2-algebroid satisfying some compatibility conditions.  Note that the image is not closed under the multiplication in the $\LWX$ 2-algebroid anymore, and this is the main difference between strict Dirac structures and weak Dirac structures. Such ideas had already been used in \cite{sz:intsemi} to integrate semidirect product Lie 2-algebras. Then we show that  a Maurer-Cartan element $H+K$, where $H\in A_{-1}\odot A_{-2}$ and $K\in\wedge^3 A_{-2}$,  of the homotopy Poisson algebra $(\Sym(\huaA[-3]),[\cdot]_S,[\cdot,\cdot]_S,[\cdot,\cdot,\cdot]_S)$ associated to a split Lie 2-algebroid $\huaA=(A_{-2},A_{-1},l_1,l_2,l_3,a)$ gives rise to a split Lie 2-algebroid structure on $\huaA^*[3]$ such that $(H^\sharp,-H^\natural,-K^\flat)$ is a morphism from the split Lie $2$-algebroid $\huaA^*[3]$ to the split Lie $2$-algebroid $\huaA$. Consequently, the graph of $(H^\sharp,-H^\natural)$ is a weak Dirac structure of the $\LWX$ 2-algebroid $\huaA\oplus \huaA^*[3]$. Finally, we generalize the above result to the case of split Lie 2-bialgebroids. We also give various examples including the string Lie 2-algebra, integrable distributions and the split Lie 2-algebroids constructed from  left-symmetric algebroids.

\vspace{2mm}
\noindent
{\bf Acknowledgements. }  We give our warmest thanks to the referee for valuable comments and suggestions that improved the paper a lot.
\section{Leibniz $2$-algebras}

The notion of a strongly homotopy
Leibniz algebra, or a $Lod_\infty$-algebra was given in \cite{livernet} by Livernet,
 which was further studied by Ammar and Poncin
in \cite{ammardefiLeibnizalgebra}.
  In \cite{Leibniz2al}, the authors introduced the notion of a Leibniz 2-algebra, which is the categorification of a Leibniz algebra, and proved that the category of Leibniz 2-algebras and the category of 2-term $Lod_\infty$-algebras are equivalent. Here we use a shift 1 version of Leibniz 2-algebras.

\begin{defi}\label{defi:2leibniz}
  A   {\bf Leibniz $2$-algebra} $\huaV$ consists of the following data:
\begin{itemize}
\item[$\bullet$] a complex of vector spaces $\huaV:V_{-2}\stackrel{\dM}{\longrightarrow}V_{-1},$

\item[$\bullet$] bilinear maps $l_2:V_{i}\times V_{j}\longrightarrow
V_{i+j+1}$, where  $-3\leq i+j\leq-2$,

\item[$\bullet$] a  trilinear map $l_3:V_{-1}\times V_{-1}\times V_{-1}\longrightarrow
V_{-2}$,
   \end{itemize}
   such that for all $w,x,y,z\in V_{-1}$ and $m,n\in V_{-2}$, the following equalities are satisfied:
\begin{itemize}
\item[$\rm(a)$] $\dM l_2(x,m)=l_2(x,\dM m),$
\item[$\rm(b)$]$\dM l_2(m,x)=-l_2(\dM m,x),$
\item[$\rm(c)$]$l_2(\dM m,n)=-l_2(m,\dM n),$
\item[$\rm(d)$]$\dM l_3(x,y,z)=l_2(x,l_2(y,z))-l_2(l_2(x,y),z)-l_2(y,l_2(x,z)),$
\item[$\rm(e_1)$]$ l_3(x,y,\dM m)=l_2(x,l_2(y,m))-l_2(l_2(x,y),m)-l_2(y,l_2(x,m)),$
\item[$\rm(e_2)$]$ -l_3(x,\dM m,y)=l_2(x,l_2(m,y))-l_2(l_2(x,m),y)-l_2(m,l_2(x,y)),$
\item[$\rm(e_3)$]$-l_3(\dM m,x,y)=l_2(m,l_2(x,y))+l_2(l_2(m,x),y)-l_2(x,l_2(m,y)),$
\item[$\rm(f)$] the Jacobiator identity:\begin{eqnarray*}
&&l_2(x,l_3(y,z,w))-l_2(y,l_3(x,z,w))+l_2(z,l_3(x,y,w))-l_2(l_3(x,y,z),w)\\
&&-l_3(l_2(x,y),z,w)-l_3(y,l_2(x,z),w)-l_3(y,z,l_2(x,w))\\
&&+l_3(x,l_2(y,z),w)+l_3(x,z,l_2(y,w))-l_3(x,y,l_2(z,w))=0.\end{eqnarray*}
   \end{itemize}
\end{defi}
We usually denote a  Leibniz 2-algebra by
$(V_{-2},V_{-1},\dM,l_2,l_3)$, or simply by
$\huaV$. In particular, if $l_2$ is graded symmetric and $l_3$ is totally skew-symmetric, we call it a {\bf Lie 2-algebra}. A Lie 2-algebra used in this paper is equivalent to a 2-term $L_\infty[1]$-algebra.

\begin{defi}\label{defi:Leibniz morphism}
 Let $\huaV$ and $\huaV^\prime$ be  Leibniz $2$-algebras. A {\bf morphism} $F$
 from $\huaV$ to $\huaV^\prime$ consists of
\begin{itemize}
  \item[$\bullet$] linear maps $F_1:V_{-1}\longrightarrow V_{-1}^\prime$
  and $F_2:V_{-2}\longrightarrow V_{-2}^\prime$ commuting with the
  differential, i.e.
  \begin{equation}\label{eq:Leibniz morphism1}
F_1\circ \dM=\dM^\prime\circ F_2;
\end{equation}
  \item[$\bullet$] a bilinear map $F_3:V_{-1}\times V_{-1}\longrightarrow
  V_{-2}^\prime$,
\end{itemize}
  such that  for all $x,y,z\in V_{-1},~m\in
V_{-2}$, we have
\begin{equation}\label{eq:Leibniz morphism2}\left\{\begin{array}{rll}
F_1l_2(x,y)-l_2^\prime(F_1(x),F_1(y))&=&\dM^\prime
F_3(x,y),\\
F_2l_2(x,m)-l_2^\prime(F_1(x),F_2(m))&=&F_3(x,\dM
m),\\
F_2l_2(m,x)-l_2^\prime(F_2(m),F_1(x))&=&-F_3(\dM
m,x),\end{array}\right.
\end{equation}
and
\begin{eqnarray}
\nonumber&&-F_2(l_3(x,y,z))+l_2^\prime(F_1(x),F_3(y,z))-l_2^\prime(F_1(y),F_3(x,z))+l_2^\prime(F_3(x,y),F_1(z))\\
\label{eq:Leibniz morphism3}
&&-F_3(l_2(x,y),z)+F_3(x,l_2(y,z))-F_3(y,l_2(x,z))+l_3^\prime(F_1(x),F_1(y),F_1(z))=0.
\end{eqnarray}
\end{defi}
In particular, if $\huaV$ and $\huaV'$ are Lie $2$-algebras and $F_3$ is skew-symmetric, we obtain the definition of a morphism between Lie $2$-algebras.

\emptycomment{
If $(F_1,F_2)$ is an injective map of underlying complexes, we say that $(F_1,F_2,F_3)$ is injective. In particular, if $(F_1,F_2)$ are inclusion maps, we call $\huaV$ a {\bf sub-Leibniz $2$-algebra of $\huaV'$ with $F_3$}.
\begin{pro}\label{pro:equivalent sub-Leibniz 2}
Let $\huaV^\prime=(V'_{-2},V'_{-1},\dM',l'_2,l'_3)$ be Leibniz $2$-algebra, $(V_{-2}, V_{-1},\dM')$ a subcomplex of $(V'_{-2}, V'_{-1},\dM')$, and $F:V_{-1}\times V_{-1}\longrightarrow
  V_{-2}^\prime$.
Then $\huaV=(V_{-2},V_{-1},\dM',l_2,l_3)$ is a sub-Leibniz $2$-algebra of $\huaV'$ with $F$ if and only if the operations $l_2$ and $l_3$ determined by

\begin{equation}\label{eq:Leibniz morphism4}\left\{\begin{array}{rll}
l_2(x,y)&=&\dM^\prime F(x,y)+l_2^\prime(x,y),\\
l_2(x,m)&=&F(x,\dM m)+l_2^\prime(x,m),\\
l_2(m,x)&=&-F(\dM m,x)+l_2^\prime(m,x),\end{array}\right.
\end{equation}
and
\begin{eqnarray}
l_3(x,y,z)&=&l_2^\prime(x,F(y,z))-l_2^\prime(y,F(x,z))+l_2^\prime(F(x,y),z)\label{eq:Leibniz morphism5}\\
\nonumber&&-F(l_2(x,y),z)+F(x,l_2(y,z))-F(y,l_2(x,z))+l_3^\prime(x,y,z)
\end{eqnarray}
are closed for $x,y,z\in V_{-1}$ and $m\in V_{-2}$.
\end{pro}
}
\section{Differential calculus on split Lie $2$-algebroids}\label{sec:D}

\subsection{Characterization of split Lie $2$-algebroids via the big bracket}

The notion of a split Lie $n$-algebroid was introduced in \cite{sz} using graded vector bundles. The equivalence between  the category of  split Lie $n$-algebroids and the category of NQ-manifolds of degree $n$ was given in \cite{BP}.

\begin{defi}
A split Lie $2$-algebroid is a graded vector bundle $\huaA=A_{-1}\oplus A_{-2}$ over a manifold $M$ equipped with a bundle map $a:A_{-1}\longrightarrow TM$, and  brackets $l_i:\Gamma(\wedge^i\huaA)\longrightarrow \Gamma(\huaA)$ of degree $1$ for $i=1,2,3$, such that
\begin{itemize}
\item[$\rm(1)$]$(\Gamma(A_{-2}),\Gamma(A_{-1}),l_1,l_2,l_3)$ is a Lie $2$-algebra;
\item[$\rm(2)$]$l_2$ satisfies the Leibniz rule with respect to $a$:
$$l_2(X^1,fY)=fl_2(X^1,Y)+a(X^1)(f)Y,\quad \forall~X^1\in\Gamma(A_{-1}),~f\in\CWM,~Y\in\Gamma(\huaA);$$
\item[$\rm(3)$]for $i\neq2$, $l_i$ are $\CWM$-linear.
 \end{itemize}
 Denote a split Lie $2$-algebroid by $(A_{-2},A_{-1},l_1,l_2,l_3,a),$ or simply by $\huaA$.
\end{defi}

Split Lie 2-algebroids become active research objects recently. See \cite{Jotz2,Jotz3,MZ,Leibniz2al,sz} for more examples and applications of split Lie 2-algebroids.

\begin{lem}
Let $(A_{-2},A_{-1},l_1,l_2,l_3,a)$ be a split Lie $2$-algebroid. Then we have
\begin{eqnarray}
\label{eq:al1}a\circ l_1&=&0,\\
\label{eq:amorphism}a(l_2(X^1,Y^1))&=&[a(X^1),a(Y^1)],\quad\forall~X^1,Y^1\in\Gamma(A_{-1}).
\end{eqnarray}
\end{lem}

\begin{defi}
Let $\huaA=(A_{-2},A_{-1},l_1,l_{2},l_3,a)$  and $\huaA'=(A_{-2}',A_{-1}',l_1',l_{2}',l_3',a')$ be split Lie $2$-algebroids over the same base manifold $M$. A morphism $F$ from $\huaA$ to $\huaA'$ consists of bundle maps  $F^{1}:A_{-1}\longrightarrow A_{-1}',$ $F^{2}:A_{-2}\longrightarrow A_{2}'$ and $F^{3}: \wedge^2 A_{-1}\longrightarrow A_{-2}'$ such that  $a'\circ F^1=a$ and $(F^1,F^2,F^3)$ is a morphism between the underlying Lie $2$-algebras.
 \end{defi}

In the sequel, we describe a split Lie 2-algebroid structure on $\huaA=A_{-1}\oplus A_{-2}$ using the graded Poisson bracket $\Poisson{\cdot,\cdot}$ on the symplectic manifold $\huaM:=T^*[3](A_{-1}\oplus A_{-2})$. Denote by $(x^i,\xi^j,\theta^k,p_i,\xi_j,\theta_k)$ a canonical Darboux coordinate on $\huaM$, where $x^i$ is a coordinate on $M$, $(\xi^j,\theta^k)$ is the fiber coordinate on $A_{-1}\oplus A_{-2}$, $(p_i,\xi_j,\theta_k)$ is the momentum coordinate on $\huaM$ for $(x^i,\xi^j,\theta^k)$. The degrees of variables $(x^i,\xi^j,\theta^k,p_i,\xi_j,\theta_k)$ are respectively $(0,1,2,3,2,1)$.
  \emptycomment{\textcolor{blue}{We introduce the tridegree of the coordinates $(x^i,\xi^j,\theta^k,p_i,\xi_j,\theta_k)$ as follows:
 $$((0,0),(0,1),(2,0),(2,1),(2,0),(0,1)).$$
 It is straightforward to check the tridegree is globally well-defined. Denote $C^n(\huaM)$ by the space of functions on $\huaM$ of degree $n$. Then the space $C^n(\huaM)$ is uniquely decomposed into the homogenous subspaces with respect to the tridegree,
 $$C^n(\huaM)=\sum_{2i+j=n}C^{2i,j}(\huaM).$$
 In particular, the space $C^4(\huaM)$ is decomposed into
  $$C^4(\huaM)=C^{4,0}(\huaM)+C^{2,2}(\huaM)+C^{0,4}(\huaM).$$}}

We introduce the {\bf tridegree} of the coordinates $(x^i,\xi^j,\theta^k,p_i,\xi_j,\theta_k)$ as follows\footnote{{We thank the referee for pointing that the third entry is the fiberwise polynomial degree, the second entry is the fiberwise
polynomial degree after applying the Legendre transformation and the first entry is chosen  so that the sum is the total degree. }}:
 $$((0,0,0),(0,1,0),(1,1,0),(1,1,1),(1,0,1),(0,0,1)).$$
 It is straightforward to check that the tridegree  is globally well-defined.  Denote by $C^n(\huaM)$ the space of functions on $\huaM$ of degree $n$. Then the space $C^n(\huaM)$ is uniquely decomposed into the homogenous subspaces with respect to the tridegree,
 $$C^n(\huaM)=\sum_{i+j+k=n}C^{(i,j,k)}(\huaM).$$

 The degree and tridegree of the symplectic structure $\omega=dx^idp_i+d\xi^j d\xi_j+d\theta_kd\theta^k$ is $3$ and $(1,1,1)$ respectively and the degree and tridegree  of the corresponding graded Poisson structure $\Poisson{\cdot,\cdot}$  is $-3$ and $(-1,-1,-1)$ respectively.

Now we consider the following fiberwise  linear function $\mu$ of degree $4$ on $\huaM$:
\begin{equation} \label{eq:L2Amu}
\mu={\mu_1}^i_j(x)p_i\xi^j+{\mu_2}^{i}_j(x)\xi_i\theta^j+\half{\mu_3}^k_{ij}(x)\xi_k\xi^i\xi^j+ {\mu_4}_{ij}^k\theta_k\xi^i\theta^j+\frac{1}{6}{\mu_5}_{ijk}^l(x)\theta_l\xi^i\xi^j\xi^k,
\end{equation}
where ${\mu_1}^i_j,~{\mu_2}^{i}_j,~{\mu_3}^k_{ij},~{\mu_4}_{ij}^k,~{\mu_5}_{ijk}^l$ are functions on $ M$. The function $\mu$ can be uniquely decomposed into
 $$\mu=\mu^{(2,1,1)}+\mu^{(1,2,1)}+\mu^{(0,3,1)},$$
where
\begin{eqnarray*}
\mu^{(2,1,1)}={\mu_2}^{i}_j(x)\xi_i\theta^j,\quad
\mu^{(1,2,1)}={\mu_1}^i_j(x)p_i\xi^j+\half{\mu_3}^k_{ij}(x)\xi_k\xi^i\xi^j+{\mu_4}_{ij}^k\theta_k\xi^i\theta^j,\quad
\mu^{(0,3,1)}=\frac{1}{6}{\mu_5}_{ijk}^l(x)\theta_l\xi^i\xi^j\xi^k.
\end{eqnarray*}


Define $l_1:A_{-2}\longrightarrow A_{-1}$, $l_2:\Gamma(A_{i})\times\Gamma( A_{j})\longrightarrow \Gamma(A_{i+j+1}),~-3\leq i+j\leq -2$,  $l_3:\wedge^3A_{-1}\longrightarrow A_{-2}$  and a bundle map $a:A_{-1}\longrightarrow TM$ by
\begin{eqnarray}\label{eq:Lie 2-algebroid brackets}
\left\{\begin{array}{rcl}
l_1(X^2)&=&-\Poisson{\mu^{(2,1,1)} ,X^2},\\
l_2(X^1,Y^1)&=&-\Poisson{\Poisson{\mu^{(1,2,1)},X^1},Y^1},\\
l_2(X^1,Y^2)&=&-\Poisson{\Poisson{\mu^{(1,2,1)},X^1},Y^2},\\
l_3(X^1,Y^1,Z^1)&=&-\Poisson{\Poisson{\Poisson{\mu^{(0,3,1)},X^1},Y^1},Z^1},\\
 a(X^1)(f)&=&-\Poisson{\Poisson{\mu^{(1,2,1)},X^1},f},
\end{array}\right.
\end{eqnarray}
for all $X^1,Y^1,Z^1\in\Gamma(A_{-1})$, $X^2,Y^2\in\Gamma(A_{-2})$ and $f\in\CWM$.

\begin{thm}\label{thm:L2A-MST}
Let  $\huaA=A_{-2} \oplus A_{-1}$ be a graded vector bundle and $\mu$ a degree $4$ function given by \eqref{eq:L2Amu}. If $\Poisson{\mu,\mu}=0$, then  $(A_{-2},A_{-1},l_1,l_2,l_3,a)$ is a split Lie $2$-algebroid, where $l_1$, $l_2$, $l_3$ and $a$ are given by \eqref{eq:Lie 2-algebroid brackets}.

Conversely, if $(A_{-2},A_{-1},l_1,l_2,l_3,a)$ is a split Lie $2$-algebroid, then we have $\Poisson{\mu,\mu}=0$, where $\mu$ is given by \eqref{eq:L2Amu}, in which ${\mu_1}^i_j,~{\mu_2}^{i}_j,~{\mu_3}^k_{ij},~{\mu_4}_{ij}^k,~{\mu_5}_{ijk}^l$ are given by:
\begin{eqnarray*}
a(\xi_j)={\mu_1}^i_j\frac{\partial}{\partial x^i},\quad l_1(\theta_j)={\mu_2}^{i}_j\xi_i,\quad
l_2(\xi_i,\xi_j)={\mu_3}^k_{ij}\xi_k,\quad l_2(\xi_i,\theta_j)= {\mu_4}^k_{ij}\theta_k,\quad l_3(\xi_i,\xi_j,\xi_k)={\mu_5}_{ijk}^l\theta_l.
\end{eqnarray*}
\end{thm}

\subsection{Differential calculus on Lie $2$-algebroids}
Let $\huaA=(A_{-2},A_{-1},l_1,l_2,l_3,a)$ be a split Lie $2$-algebroid with the structure function $\mu$. Then we have the generalized Chevalley-Eilenberg cochain complex $(\Sym(\huaA^*)=\bigoplus_k \Sym^k(\huaA^*),\delta)$, where the set of $k$-cochains  $\Sym^k(\huaA^*)$ is given by
\begin{eqnarray*}
\Sym^k(\huaA^*)=\sum_{p+2q=k}\Sym^{p,q}(\huaA^*)\quad\mbox{with}\quad \Sym^{p,q}(\huaA^*)=\Gamma(\wedge^p A^*_{-1})\odot\Gamma(\Sym^q(A^*_{-2}))
\end{eqnarray*}
and the differential $\delta:\Sym^k(\huaA^*)\longrightarrow \Sym^{k+1}(\huaA^*)$ is defined by\footnote{ $\Sym(\huaA^*)$ can be embedded into the Poisson algebra  $C^\infty(\huaM)$ through the pullback of the canonical map $T^*[3](A_{-1}\oplus A_{-2})\lon A_{-1}\oplus A_{-2}$. On the other hand, by the Legendre transformation,    $T^*[3](A_{-1}\oplus A_{-2})$ is isomorphic to $T^*[3](A^*_{-1}\oplus A^*_{-2})[3]$ as graded symplectic manifolds. Thus, $\Sym(\huaA)$ can also be embedded into the Poisson algebra  $C^\infty(\huaM)$. See \cite{Roytenbergphdthesis} for more details on the Legendre transformation.}
\begin{equation}
\delta(\cdot)=\{\mu,\cdot\}.
\end{equation}

If there is a Lie $2$-algebroid structure on $\huaA^*[3]$, we use $\delta_*$ to denote the corresponding differential. The differential   $\delta$ can be written as
$\delta=\bar{\delta}+\dM+\hat{\delta},$
where $\bar{\delta}:\Sym^{p,q}(\huaA^*)\longrightarrow \Sym^{p-1,q+1}(\huaA^*)$, $\dM:\Sym^{p,q}(\huaA^*)\longrightarrow \Sym^{p+1,q}(\huaA^*)$ and $\hat{\delta}:\Sym^{p,q}(\huaA^*)\longrightarrow \Sym^{p+3,q-1}(\huaA^*)$ are given by
\begin{eqnarray*}
\bar{\delta}(\cdot)=\Poisson{\mu^{(2,1,1)},\cdot},\quad \dM(\cdot)=\Poisson{\mu^{(1,2,1)},\cdot},\quad\hat{\delta}(\cdot)=\Poisson{\mu^{(0,3,1)},\cdot}.
\end{eqnarray*}

In particular, for $\alpha^1\in\Gamma(A^*_{-1})$, we have
\begin{eqnarray*}
  \bar{\delta}\alpha^1(X^2)=-\langle \alpha^1,l_1(X^2)\rangle,\quad\forall~X^2\in\Gamma(A_{-2}).
\end{eqnarray*}

For all $f\in\CWM,~\alpha^1\in\Gamma(A^*_{-1}),~\alpha^2\in\Gamma(A_{-2}^*)$, $X^1,~Y^1\in\Gamma(A_{-1}),~Y^2\in\Gamma(A_{-2})$, we have
\begin{eqnarray*}\label{eq:delta}
\left\{\begin{array}{rcl}
\dM (f)(X^1)&=&a(X^1)(f),\\
\dM(\alpha^1)(X^1,Y^1)&=&a(X^1)\langle\alpha^1,Y^1\rangle-a(Y^1)\langle\alpha^1,X^1\rangle-\langle\alpha^1,l_2(X^1,Y^1)\rangle,\\
\dM(\alpha^2)(X^1,Y^2)&=&a(X^1)\langle\alpha^2,Y^2\rangle-\langle\alpha^2,l_2(X^1,Y^2)\rangle.
\end{array}\right.
\end{eqnarray*}

For all $\alpha^2\in\Gamma(A_{-2}^*)$, we have
\begin{eqnarray*}
  \hat{\delta}\alpha^2(X^1,Y^1,Z^1)=-\langle l_3(X^1,Y^1,Z^1),\alpha^2\rangle,\quad\forall~X^1,Y^1,Z^1\in\Gamma(A_{-1}).
\end{eqnarray*}

By the properties of graded Poisson bracket, for all $\phi_1\in \Sym^k(\huaA^*)$ and $\phi_2\in \Sym^l(\huaA^*)$, we have
\begin{equation}
  \delta(\phi_1\odot\phi_2)=\delta(\phi_1)\odot\phi_2+(-1)^{k}\phi_1\odot\delta(\phi_2).
\end{equation}

Define the Lie derivative $L^0:\Sym^{p,q}(\huaA^*)\longrightarrow \Sym^{p-1,q+1}(\huaA^*)$ by
 \begin{equation}
   L^0(\phi)=-\Poisson{\mu^{(2,1,1)},\phi},\quad\forall~\phi\in \Sym^{p,q}(\huaA^*).
 \end{equation}
In particular, for all $\alpha^1\in\Gamma(A^*_{-1})$, we have
\begin{equation*}
\langle  L^0(\alpha^1),X^2\rangle=\langle\alpha^1,l_1(X^2)\rangle,\quad\forall~X^2\in\Gamma(A_{-2}).
\end{equation*}
It is obvious that $ L^0(\alpha^1)=l^*_1(\alpha^1)$.

For all $X^1\in\Gamma(A_{-1})$, define the Lie derivative $L^1_{X^1}:\Sym^{p,q}(\huaA^*)\longrightarrow \Sym^{p,q}(\huaA^*)$ by
\begin{eqnarray}\label{eq:L1}
 L^1_{X^1}\phi=-\Poisson{\Poisson{\mu^{(1,2,1)},X^1},\phi},\quad \forall~\phi\in \Sym^{p,q}(\huaA^*).
\end{eqnarray}

In particular, for all $\alpha^i\in\Gamma(A^*_{-i}),~i=1,2$, we have
\begin{eqnarray*}
\langle L^1_{X^1}\alpha^i,Y^i\rangle&=&a(X^1)\langle Y^i,\alpha^i\rangle-\langle \alpha^i,l_2(X^1,Y^i)\rangle,\quad \forall~ Y^i\in\Gamma(A_{-i}).
\end{eqnarray*}
It is straightforward to deduce that
\begin{equation}
  L^1_{X^1}(\phi\odot\psi)= L^1_{X^1}{\phi}\odot\psi+\phi\odot L^1_{X^1}\psi,\quad \forall~\phi\in \Sym^k(\huaA^*),~\psi\in \Sym^l(\huaA^*).
\end{equation}

\emptycomment{For all $X^2\in\Gamma(A_{-2})$, define $L^2_{X^2}:\Sym^{p,q}\longrightarrow \Sym^{p,{q-1}}$ by
\begin{equation}
 L^2_{X^2}\alpha^2=-\Poisson{\Poisson{\mu^{(1,2,1)},X^2},\alpha^2},\quad \forall~\alpha^2\in\Gamma(A^*_{-2}),X^2\in\Gamma(A_{-2}),
\end{equation}
or equivalently,}

For all $X^2\in\Gamma(A_{-2})$, define the Lie derivative $L^2_{X^2}:\Sym^{p,q}(\huaA^*)\longrightarrow \Sym^{p+1,{q-1}}(\huaA^*)$ by
\begin{equation}
 L^2_{X^2}\phi=-\Poisson{\Poisson{\mu^{(1,2,1)},X^2},\phi},\quad \forall~\phi\in \Sym^{p,q}(\huaA^*).
\end{equation}
In particular, for all $\alpha^2\in\Gamma(A_{-2}^*)$, we have
\begin{equation*}
\langle L^2_{X^2}\alpha^2,Y^1\rangle=-\langle \alpha^2,l_2(X^2,Y^1)\rangle,\quad \forall~Y^1\in\Gamma(A_{-1}).
\end{equation*}
It is easy to see that
\begin{equation}
  L^2_{X^2}(\phi\odot\psi)= L^2_{X^2}{\phi}\odot\psi+(-1)^k\phi\odot L^2_{X^2}\psi,\quad \forall~\phi\in \Sym^k(\huaA^*),~\psi\in \Sym^l(\huaA^*).
\end{equation}

For all $X^1,Y^1\in\Gamma(A_{-1})$, define the Lie derivative $L^3_{X^1,Y^1}:\Sym^{p,q}(\huaA^*)\longrightarrow \Sym^{p+1,q-1}(\huaA^*)$ by
\begin{equation}
 L^3_{X^1,Y^1}\phi=-\Poisson{\Poisson{\Poisson{\mu^{(0,3,1)} ,X^1},Y^1},\phi},\quad \forall~\phi\in \Sym^{p,q}(\huaA^*).
\end{equation}
In particular, for all $\alpha^2\in\Gamma(A_{-2}^*)$, we have
\begin{equation*}
\langle L^3_{X^1,Y^1}\alpha^2,Z^1\rangle=-\langle \alpha^2,l_3(X^1,Y^1,Z^1)\rangle,\quad \forall~Z^1\in\Gamma(A_{-1}).
\end{equation*}
It is easy to see that
\begin{equation}
 L^3_{X^1,Y^1}(\phi\odot\psi)=  L^3_{X^1,Y^1}{\phi}\odot\psi+(-1)^k\phi\odot L^3_{X^1,Y^1}\psi,\quad \forall \phi\in \Sym^k(\huaA^*),~\psi\in \Sym^l(\huaA^*).
\end{equation}

For all $ X\in\Gamma(\huaA),~ \alpha^1_i\in\Gamma(A^*_{-1}),~\alpha^2_j\in\Gamma(A^*_{-2})$, define the contraction operator $\iota_{X}:\Sym^{p,q}(\huaA^*)\lon  \Sym^{p-1,q}(\huaA^*)\oplus \Sym^{p,q-1}(\huaA^*)$ by
\begin{eqnarray*}
 \iota_{X} (\alpha^1_1\odot \cdots\odot\alpha^1_p\odot\alpha^2_{1}\odot\cdots\odot \alpha^2_{q})
  &=&\sum^{p}_{i=1}(-1)^{i+1}\langle X,\alpha^1_i\rangle \alpha^1_1\odot  \cdots\odot  \hat{\alpha^1_i}\odot\cdots\odot\alpha^1_p \odot\alpha^2_{1}\odot\cdots\odot \alpha^2_{q}\\
  &&+\sum^{q}_{j=1}\langle X,\alpha^2_{j}\rangle \alpha^1_1\odot\cdots\odot\alpha^1_p \odot\alpha^2_{1}\odot\cdots\odot  \hat{\alpha^2_{j}}\odot\cdots \odot\alpha^2_{q}.
\end{eqnarray*}

For any $\phi\in \Sym^{k}(\huaA^*)$, let us denote by
$$\phi(X_1,X_2,\cdots,X_{k})=\iota_{X_{k}}\iota_{X_{k-1}}\cdots \iota_{X_1}\phi,\quad \forall~X_i\in\Gamma(\huaA).$$

Thus, for all $X^1\in\Gamma(A_{-1}),~\phi\in \Sym^{k}(\huaA^*)$ and $Y_i\in\Gamma(\huaA)$, we have
\begin{eqnarray*}
   (L^1_{X^1}\phi)(Y_1, \cdots,Y_{k})&=&a(X^1)\phi(Y_1, \cdots,Y_{k})- \sum_{i=1}^{p+q}\phi(Y_1, \cdots,l_2(X^1,Y_i),\cdots,Y_{k}).
\end{eqnarray*}

The following lemmas list some properties of the above operators.
\begin{lem}
 For all $X^1\in\Gamma(A_{-1}),~X^2\in\Gamma(A_{-2}),~f\in\CWM$ and $\phi\in \Sym^k(\huaA^*)$, we have
 \begin{eqnarray*}\begin{array}{rclrcl}
  L^1_{X^1}f\phi&=&f(L^1_{X^1}\phi)+a(X^1)(f)\phi, &
 L^1_{fX^1}\phi&=&f(L^1_{X^1}\phi)+\dM f\odot \iota_{X^1}\phi,\\
  L^2_{X^2}f\phi&=&f(L^2_{X^2}\phi), &
   L^2_{fX^2}\phi&=&f( L^2_{X^2}\phi)-\dM f\odot \iota_{X^2}\phi,\\
  L^1_{X^1}\phi&=&\iota_{X^1}\dM \phi+\dM\iota_{X^1}\phi, &
  L^2_{X^2}\phi&=&\iota_{X^2}\dM \phi-\dM\iota_{X^2}\phi.\end{array}
 \end{eqnarray*}
 \end{lem}

 \begin{lem}\label{lem:LieJacobi}
 For all $X^1,Y^1\in\Gamma(A_{-1}),~X^2\in\Gamma(A_{-2}),~\phi\in \Sym^k(\huaA^*)$, we have
 \begin{eqnarray}
\label{eq:jacobi-eq1}L^1_{l_2(X^1, Y^1)}\phi- L^1_{X^1}L^1_{Y^1}\phi+L^1_{Y^1}L^1_{X^1}\phi&=&-L^3_{X^1,Y^1}L^0(\phi)-L^0(L^3_{X^1,Y^1}\phi),\\
 L^2_{l_2(X^1, Y^2)}\phi-L^1_{X^1}L^2_{Y^2}\phi+L^2_{Y^2}L^1_{X^1}\phi&=&-L^3_{l_1(Y^2),X^1}\phi.
 \end{eqnarray}
 \end{lem}

 \begin{lem}
For all $X^1,Y^1\in\Gamma(A_{-1}),~X^2\in\Gamma(A_{-2}),~\alpha^1\in\Gamma(A^*_{-1}),~\alpha^2\in\Gamma(A_{-2}^*)$, we have
 \begin{eqnarray}
\label{eq:interi-eq1}\iota_{l_2(X^1, Y^1)}\dM\alpha^1- L^1_{X^1}\iota_{Y^1}\dM\alpha^1+\iota_{Y^1}L^1_{X^1}\dM\alpha^1&=&-L^3_{X^1,Y^1}l^*_1\alpha^1,\\
 \label{eq:interi-eq2}\iota_{l_2(X^1, Y^1)}\dM\alpha^2- L^1_{X^1}\iota_{Y^1}\dM\alpha^2+\iota_{Y^1}L^1_{X^1}\dM\alpha^2&=&-l^*_1L^3_{X^1,Y^1}\alpha^2,\\
\label{eq:interi-eq3} \iota_{l_2(X^1, Y^2)}\dM\alpha^2-L^1_{X^1}\iota_{Y^2}\dM\alpha^2+\iota_{Y^2}L^2_{X^1}\dM\alpha^2&=&-L^3_{l_1(Y^2),X^1}\alpha^2.
 \end{eqnarray}
 \end{lem}

If $\huaA^*[3]$ is a split Lie 2-algebroid, we use $\huaL^0,~\huaL^1,~\huaL^2$ and $\huaL^3$ to denote the corresponding Lie derivatives.

\section{Homotopy Poisson algebras of degree 3 associated to split Lie $2$-algebroids and split Lie 2-bialgebroids}\label{sec:H}

\subsection{Homotopy Poisson algebras of degree 3 associated to split Lie $2$-algebroids}

Associated to a Lie algebroid $A$, we have the Gerstenhaber algebra $(\Gamma(\wedge^\bullet A),\wedge,[\cdot,\cdot])$, which is also a Poisson algebra of degree $-1$, where $\wedge^\bullet A=\sum_k\wedge^kA$ and $[\cdot,\cdot]$ is the Schouten bracket on $\Gamma(\wedge^\bullet A)$. This algebra plays very important role in the theory of Lie algebroids. Note that $\Gamma(\wedge^\bullet A)$ can be understood as the symmetric algebra of $A[-1]$ and elements in $\Gamma(\wedge^kA)$ are of degree $k$. This algebra can be obtained by the derived bracket as follows. Consider the shifted cotangent bundle $T^*[2]A[1]$, which is a symplectic manifold of degree 2. The corresponding Poisson structure is of degree $-2$. The Lie algebroid structure is equivalent to a degree 3 function $\kappa$ on $T^*[2]A[1]$ which is fiberwise linear. Then the Schouten bracket on $\Gamma(\wedge^\bullet A)$ can be obtained by
$$
[P,Q]=-\{\{\kappa,P\},Q\},\quad \forall~P\in\Gamma(\wedge^kA), ~Q\in\Gamma(\wedge^lA).
$$
By the fact that the degree of $\kappa$ is 3 and the degree of the Poisson bracket is $-2$, we deduce that the degree of the Schouten bracket is $3-2-2=-1$.

Now for a split Lie 2-algebroid, using the above idea, we define higher bracket operations on its symmetric algebra using the canonical Poisson bracket on the shifted cotangent bundle.

 Given a split Lie $2$-algebroid $\huaA=(A_{-2},A_{-1},l_1,l_2,l_3,a)$ with the structure function $\mu$ given by \eqref{eq:L2Amu}, denote by $(\Sym(\huaA[-3])=\bigoplus_k \Sym^k(\huaA[-3]),\odot)$ the symmetric algebra of $\huaA[-3]$, in which $\Sym^k(\huaA[-3])$ is given by\footnote{For the shifted vector bundle $\huaA[-3]$, elements in $\Gamma(A_{-2})$ are of degree $1$ and  elements in $\Gamma(A_{-1})$ are of degree $2$.}
\begin{eqnarray*}
\Sym^k(\huaA[-3])&=&\sum_{p+2q=k}\Sym^{p,q}(\huaA[-3]),
\end{eqnarray*}
where $\Sym^{p,q}(\huaA[-3])=\Gamma(\wedge^p A_{-2})\odot\Gamma(\Sym^q(A_{-1})).$ We will use $|\cdot|$ to denote the degree of a homogeneous element in $\Sym(\huaA[-3])$.

For all $P\in \Sym^{k}(\huaA[-3]),~Q\in \Sym^{l}(\huaA[-3]),~R\in \Sym^{m}(\huaA[-3])$, define $\Schouten{\cdot}:\Sym^{k}(\huaA[-3])\longrightarrow \Sym^{k+1}(\huaA[-3])$ by
\begin{equation}
\Schouten{P}=-\Poisson{\mu^{(2,1,1)},P}.
\end{equation}
Define $\Schouten{\cdot,\cdot}: \Sym^{k}(\huaA[-3])\times \Sym^{l}(\huaA[-3])\longrightarrow \Sym^{k+l-2}(\huaA[-3])$ by
\begin{equation}\label{eq:defiSB}
\Schouten{P,Q}=-\Poisson{\Poisson{\mu^{(1,2,1)},P},Q}.
\end{equation}
Define $\Schouten{\cdot,\cdot,\cdot}:\Sym^{k}(\huaA[-3])\times \Sym^{l}(\huaA[-3])\times \Sym^{m}(\huaA[-3])\longrightarrow \Sym^{k+l+m-5}(\huaA[-3])$ by
\begin{eqnarray}
\Schouten{P,Q,R}=-\Poisson{\Poisson{\Poisson{\mu^{(0,3,1)} ,P},Q},R}.
\end{eqnarray}

Comparing with \eqref{eq:Lie 2-algebroid brackets}, it is straightforward to obtain that
\begin{lem}
 With the above notations, for all $X^1,Y^1,Z^1\in\Gamma(A_{-1}),~X^2\in\Gamma(A_{-2})$, we have
\begin{eqnarray*}
\Schouten{X^2}=l_1(X^2),\quad
{[X^1,f]}_S=a(X^1)(f),\quad
{[X^1,Y^1]}_S=l_2(X^1,Y^1),\\
{[X^1,X^2]}_S=l_2(X^1,X^2),\quad
\Schouten{X^1,Y^1,Z^1}=l_3(X^1,Y^1,Z^1).
\end{eqnarray*}
\end{lem}
For $\alpha\in\Gamma(\huaA^*) $, define $\iota_{\alpha}:\Sym^{p,q}(\huaA[-3])\longrightarrow \Sym^{p-1,q}(\huaA[-3])\oplus \Sym^{p,q-1}(\huaA[-3])$ by
\begin{eqnarray*}
\iota_{\alpha} (X^1_1\odot \cdots\odot X^1_p\odot X^2_{1}\odot\cdots\odot X^2_{q})
 =\sum^{p}_{i=1}\langle \alpha,X^1_i\rangle X^1_1\odot\cdots\odot  \hat{X^1_i}\odot\cdots\odot X^1_p \odot X^2_{1}\odot\cdots\odot X^2_{q}\\
  +\sum^{q}_{j={1}}(-1)^{j+1}\langle \alpha,X^2_{j}\rangle X^1_1\odot\cdots\odot X^1_p \odot X^2_{1}\odot\cdots\odot  \hat{X^2_{j}}\odot\cdots \odot X^2_{q},
\end{eqnarray*}
where $ X^1_i\in\Gamma(A_{-1}[-3]),~X^2_j\in\Gamma(A_{-2}[-3])$.

For all $P\in \Sym^{p}(\huaA[-3])$, let us denote by
$$P(\alpha_1,\cdots,\alpha_{p})=\iota_{\alpha_{p}}\cdots \iota_{\alpha_1} P,\quad \forall~ \alpha_i\in\Gamma(\huaA^*).$$

Using the properties of the graded Poisson bracket $\Poisson{\cdot,\cdot}$, we get the following formulas.
\begin{pro}\label{lem:Higher GA}
For all $P\in  \Sym^{\degree{P}}(\huaA[-3]), ~ Q\in \Sym^{\degree{Q}}(\huaA[-3]),~R\in \Sym^{\degree{R}}(\huaA[-3])$ and $W\in \Sym^{\degree{W}}(\huaA[-3])$,   we have
\begin{eqnarray*}
\Schouten{P,fQ}&=&f\Schouten{P,Q}+(-1)^{\degree{P}}\iota_{\delta(f)}P\odot Q,\\
\Schouten{P\odot Q}&=&\Schouten{P}\odot Q+(-1)^{\degree{P}}P\odot\Schouten{Q},\\
{[P,Q]}_S&=&(-1)^{(\degree{P}-1)(\degree{Q}-1)}{[Q,P]_S},\\
{[P,Q\odot R]}_S&=&{[P,Q]}_S\odot R+(-1)^{\degree{P}\degree{Q}}Q\odot {[P,R]}_S,\\
\Schouten{P,Q,R}&=&(-1)^{(\degree{P}-3)(\degree{Q}-3)}\Schouten{Q,P,R}=(-1)^{(\degree{R}-3)(\degree{Q}-3)}\Schouten{P,R,Q},\\
\Schouten{P,Q,R\odot W}&=&\Schouten{P,Q,R}\odot W +(-1)^{(\degree{P}+\degree{Q}-5)\degree{R}}R\odot\Schouten{P,Q,W}.
\end{eqnarray*}
\end{pro}

It is easy to see that for all $X\in \Gamma(\huaA)$, we have
\begin{eqnarray*}
[X,Y_1\odot \cdots\odot Y_n]_S&=&\sum_{i=1}^nY_1\odot\cdots\odot l_2(X,Y_i) \odot\cdots\odot Y_n,\quad\forall~Y_i\in\Gamma(\huaA).
\end{eqnarray*}

For all $X^1\in \Gamma(A_{-1}),~D\in \Sym^{p,q}(\huaA[-3])$ and $\alpha_i\in\Gamma(\huaA^*)$, we have
\begin{eqnarray*}
   ([X^1,D]_S)(\alpha_1, \cdots,\alpha_{p+q})=a(X^1) D(\alpha_1, \cdots,\alpha_{p+q}) - \sum_{i=1}^{p+q}D(\alpha_1, \cdots,L^1_{X^1}\alpha_i,\cdots,\alpha_{p+q}).
\end{eqnarray*}

\begin{pro}\label{pro:hP}
For all $P\in  \Sym^{\degree{P}}(\huaA[-3]),~Q\in \Sym^{\degree{Q}}(\huaA[-3]),~R\in \Sym^{\degree{R}}(\huaA[-3])$ and $ W\in \Sym^{\degree{W}}(\huaA[-3])$, we have
\begin{eqnarray*}
[{[P,Q]}_S]_S=-[[P]_S,Q]_S+(-1)^{\degree{P}}[P,[Q]_S]_S,\\
{[P,{[Q,R]}_S]}_S-(-1)^{\degree{P}}{[{[P,Q]}_S,R]}_S-(-1)^{\degree{P}\degree{Q}}{[Q,{[P,R]}_S]}_S
=(-1)^{\degree{P}}[[P,Q,R]_S]_S\\
+(-1)^{\degree{Q}}[P,Q,[R]_S]_S-[P,[Q]_S,R]_S+(-1)^{\degree{P}}[[P]_S,Q,R]_S,
\end{eqnarray*}
and
\begin{eqnarray*}
&&-(-1)^{\degree{P}}{[P,[Q,R,W]_S]}_S+(-1)^{\degree{P}(\degree{Q}-1)}{[Q,[P,R,W]_S]}_S+(-1)^{(\degree{P}+\degree{Q}-1)(\degree{R}-1)}{[R,[P,Q,W]_S]}_S\nonumber\\
&&+{[[P,Q,R]_S ,W]}_S-(-1)^{(\degree{R}-1)(\degree{W}-1)}[[P,Q]_S,R,W]_S+(-1)^{\degree{P}(\degree{Q}-1)}[Q,{[P,R]}_S ,W]_S\\
&&-(-1)^{\degree{P}(\degree{R}+\degree{Q})}[Q,R,{[P ,W]}_S]_S-(-1)^{\degree{P}}[P,{[Q,R]}_S, W]_S-(-1)^{\degree{P}+\degree{Q}}[P,Q,{[R ,W]}_S]_S\\
&&+(-1)^{\degree{R}\degree{Q}-\degree{P}-\degree{Q}}[P,R,{[Q ,W]}_S]_S=0.
\end{eqnarray*}
\end{pro}
\pf It follows from the graded Jacobi identity for the Poisson bracket $\{\cdot,\cdot\}$. We omit details. \qed\vspace{3mm}

To understand the meaning of the above brackets, we need the notion of a homotopy Poisson algebra (\cite{LSX,Mehta}).

\begin{defi}
\begin{itemize}
\item[\rm(i)]
A {\bf  homotopy Poisson algebra} of degree $n$   is a graded commutative algebra $\frka$ over a field of characteristic zero with an $L_\infty[1]$-algebra structure $\{l_m\}_{m\geq1}$ on $\frka[n]$, such that the map
\[x\longrightarrow{l_m(x_1,\cdots,x_{m-1},x)},\ \ \ \ x_1,\cdots,x_{m-1},x\in \frka\]
is a  derivation of $\frka$ of degree $\kappa:=\sum_{i=1}^{m-1}|x_i|+1-n(m-1)$, i.e. for all $x,y\in \frka$, we have
\begin{eqnarray*}
l_m(x_1,\cdots,x_{m-1},xy)=l_m(x_1,\cdots,x_{m-1},x)y+(-1)^{\kappa|x|}
x l_m(x_1,\cdots,x_{m-1},y).
\end{eqnarray*}
Here, $|x_i|$ denotes the degree of a homogeneous element $x_i\in\frka.$

\item[\rm(ii)]A  homotopy Poisson algebra of degree $n$ is of {\bf finite type} if there exists a $q$ such that $l_m=0$ for all $m>q$.

\item[\rm(iii)]A {\bf homotopy Poisson manifold} of degree $n$ is a graded manifold $\huaM$ whose algebra of functions $C^\infty(\huaM)$ is equipped with a degree $n$ homotopy Poisson algebra structure of finite type.
    \item[\rm(iv)]
   A {\bf Maurer-Cartan element} of a  homotopy Poisson algebra  of degree $n$ is a degree $n$ element $m$ satisfying
    $$
    l_1(m)+\frac{1}{2}l_2(m,m)+\frac{1}{6}l_3(m,m,m)+\cdots=0.
    $$
      \end{itemize}
\end{defi}

The only difference between the above definition and the one provided in \cite{LSX} is that we use $L_\infty[1]$-algebra here, while in \cite{LSX} the authors used $L_\infty$-algebra. Since $L_\infty[1]$-algebras and $L_\infty$-algebras are equivalent,  there is no intrinsic difference.

The following theorem can be proved quickly using the derived bracket construction by Voronov in \cite{Voronov1}. Here we give another proof using the properties of $[\cdot]_S,~[\cdot,\cdot]_S$ and $[\cdot,\cdot,\cdot]_S$.

\begin{thm}\label{thm:hpa}
  Let  $\huaA$ be a split Lie $2$-algebroid. Then $(\Sym(\huaA[-3]),[\cdot]_S,[\cdot,\cdot]_S,[\cdot,\cdot,\cdot]_S)$   is a homotopy Poisson algebra  of degree $3$.
\end{thm}
\pf By Proposition \ref{pro:hP}, it is obvious that $[\cdot]_S,~[\cdot,\cdot]_S,~[\cdot,\cdot,\cdot]_S$ define  an $L_\infty[1]$-algebra structure on $\Sym(\huaA[-3])[3]$. By Proposition \ref{lem:Higher GA}, the graded derivation conditions for $[\cdot]_S,~[\cdot,\cdot]_S,~[\cdot,\cdot,\cdot]_S$ are satisfied. Thus, $(\Sym(\huaA[-3]),[\cdot]_S,[\cdot,\cdot]_S,[\cdot,\cdot,\cdot]_S)$   is a homotopy Poisson algebra  of degree $3$. \qed

\begin{cor}
    Let  $\huaA=(A_{-2},A_{-1},l_1,l_2,l_3,a)$ be a split Lie $2$-algebroid. Then $\huaA^*[3]=A_{-2}^*[3]\oplus A_{-1}^*[3]$ is a homotopy Poisson manifold of degree $3$.
\end{cor}

The method of defining higher bracket operations using the Poisson  bracket of the shifted cotangent bundle can be easily generalized to a split Lie $n$-algebroid. Consequently, one can obtain a  homotopy Poisson algebra  of degree $n+1$ on the symmetric algebra $\Sym(\huaA[-n-1])$ of a split Lie $n$-algebroid $\huaA$ and $\huaA^*[n+1]$ is a homotopy Poisson manifold of degree $n+1$. In \cite[Example 3.6]{LSX}, it is  {stated (but not proved)} that the dual of a Lie $n$-algebroid is a homotopy Poisson manifold of degree $n$. The difference is generated by the difference between an $L_\infty[1]$-algebra and an $L_\infty$-algebra, which is not intrinsic.

\subsection{A new approach to split Lie 2-bialgebroids}

In this subsection, we extend some results given in \cite{MackenzieX:1994,Roytenbergphdthesis} to Lie 2-bialgebroids. The notion of a split Lie 2-bialgebroid was introduced in \cite{LiuSheng} using the canonical Poisson bracket of the shifted cotangent bundle.

Now assume that there is a  split Lie $2$-algebroid structure on the dual bundle $\huaA^*[3]=A^*_{-1}[3]\oplus A^*_{-2}[3]$. The two cotangent bundles $\huaM=T^*[3](A_{-1}\oplus A_{-2})$ and $T^*[3](A^*_{-1}\oplus A^*_{-2})[3]$ are naturally isomorphic as graded symplectic manifold by the Legendre transformation. By Theorem \ref{thm:L2A-MST}, the split Lie $2$-algebroid $(\huaA^*[3],\frkl_1,\frkl_2,\frkl_3,\frka)$ gives rise to a degree $4$ function $\gamma$ on $\huaM$ satisfying $\Poisson{\gamma,\gamma}=0.$ It is given in local coordinates $(x^i,\xi^j,\theta^k, p_i,\xi_j,\theta_k)$ by
\begin{equation}\label{eq:gamma}
\gamma={\gamma_1}^{ij}(x)p_j\theta_i+{\gamma_2}^{j}_i(x)\xi_j\theta^i+\half{\gamma_3}_k^{ij}(x)\theta^k\theta_i\theta_j+ {\gamma_4}_{k}^{ij}\xi_i\theta_j\xi^k+\frac{1}{6}{\gamma_5}^{ijk}_l(x)\xi^l\theta_i\theta_j\theta_k.
\end{equation}
According to the tridegree, $\gamma$ can be decomposed into
$$\gamma=\gamma^{(2,1,1)}+\gamma^{(1,1,2)}+\gamma^{(0,1,3)}.$$
\begin{defi}{\rm(\cite{LiuSheng})}
Let $\huaA$ and $\huaA^*[3]$ be  split Lie $2$-algebroids with the structure functions $\mu$ and $\gamma$ respectively. The pair $(\huaA,\huaA^*[3])$ is called a split Lie $2$-bialgebroid if $\mu^{(2,1,1)}=\gamma^{(2,1,1)}$ and
\begin{equation}
\Poisson{\mu+\gamma-\mu^{(2,1,1)},\mu+\gamma-\mu^{(2,1,1)}}=0,
\end{equation}
where $\Poisson{\cdot,\cdot}$ is the graded Poisson bracket on the shifted cotangent bundle $T^*[3](A_{-1}\oplus A_{-2}).$
\end{defi}
Denote a split Lie $2$-bialgebroid by $(\huaA,\huaA^*[3])$.
\begin{rmk}
\begin{itemize}
\item[$\rm(1)$]
  The condition $\mu^{(2,1,1)}=\gamma^{(2,1,1)}$ is due to the invariant condition {\rm(iii)} in the definition of a $\LWX$ $2$-algebroid (see Definition \ref{defi:Courant-2 algebroid}).
\item[$\rm(2)$]Note that the function $\mu+\gamma$ contains two copies of the term $\mu^{(2,1,1)}$, which is  of the tridegree $(2,1,1)$. Thus, we use the degree $4$ function $\mu+\gamma-\mu^{(2,1,1)}$ in the definition of a split Lie $2$-bialgebroid.
\end{itemize}
\end{rmk}

Now using the homotopy Poisson algebra associated to a split Lie 2-algebroid, we can describe a split Lie 2-bialgebroid using the usual language of differential geometry similar as the case of a Lie bialgebroid (\cite{MackenzieX:1994}).

 \begin{thm}\label{thm:Lie2bi}
Let $\huaA$ and $\huaA^*[3]$ be  split Lie $2$-algebroids with the structure functions $\mu$ and   $\gamma$ respectively, such that $\mu^{(2,1,1)}=\gamma^{(2,1,1)}$. Then $(\huaA,\huaA^*[3])$ is a split Lie $2$-bialgebroid if and only if the following two conditions hold:
\begin{eqnarray}
\delta_*{[X,Y]}_S&=&-{[\delta_*(X),Y]}_S+(-1)^{|X|}{[X,\delta_*(Y)]}_S,\quad \forall ~X,Y\in\Gamma(\huaA[-3]),\label{L2bialgebroid1}\\
\delta{[\alpha,\beta]}_S&=&-{[\delta(\alpha),\beta]}_S+(-1)^{|\alpha|}{[\alpha,\delta(\beta)]}_S,\quad\forall~ \alpha,\beta\in\Gamma(\huaA^*),\label{L2bialgebroid2}
\end{eqnarray}
where $\delta_*$ and $\delta$ are coboundary operators associated to split Lie $2$-algebroids $\huaA^*[3]$ and $\huaA$ respectively.
 \end{thm}
\pf Let $\huaA$ and $\huaA^*[3]$ be  split Lie $2$-algebroids. Then by the {tridegree} reason, the following equalities are automatically satisfied:
\begin{eqnarray}
\Poisson{\mu^{(1,2,1)},\mu^{(2,1,1)}}&=&0,\quad \Poisson{\mu^{(1,2,1)},\mu^{(1,2,1)}}+2\Poisson{\mu^{(2,1,1)},\mu^{(0,3,1)}}=0,\quad\Poisson{\mu^{(1,2,1)},\mu^{(0,3,1)}}=0;\label{eq:mu-L2A}\\
\Poisson{\gamma^{(1,1,2)},\gamma^{(2,1,1)}}&=&0,\quad \Poisson{\gamma^{(1,1,2)},\gamma^{(1,1,2)}}+2\Poisson{\gamma^{(2,1,1)},\gamma^{(0,1,3)}}=0,\quad\Poisson{\gamma^{(1,1,2)},\gamma^{(0,1,3)}}=0.\label{eq:gama-L2A}
\end{eqnarray}
 If $(\huaA,\huaA^*[3])$ is a split Lie $2$-bialgebroid, by the {tridegree} reason, $\Poisson{\mu+\gamma-\mu^{(2,1,1)},\mu+\gamma-\mu^{(2,1,1)}}=0$ is equivalent to
\begin{eqnarray}
\Poisson{\mu^{(1,2,1)},\gamma^{(1,1,2)}}&=&0,\quad\Poisson{\mu^{(1,2,1)},\gamma^{(0,1,3)}}=0,\quad\Poisson{\gamma^{(1,1,2)},\mu^{(0,3,1)}}=0\label{eq:compatibility-L2A}.
\end{eqnarray}
By \eqref{eq:mu-L2A} and \eqref{eq:gama-L2A}, \eqref{eq:compatibility-L2A} is equivalent to
\begin{eqnarray*}
\label{eq:compatibility-L2A21}\Poisson{\gamma,\mu^{(1,2,1)}}=0,\quad
\label{eq:compatibility-L2A2}\Poisson{\mu,\gamma^{(1,1,2)}}=0.
\end{eqnarray*}

\emptycomment{
For all $ X^1, Y^1\in\Gamma(A_{-1}),$ we have
\begin{eqnarray*}
&&\Poisson{Y^1,\Poisson{X^1,\Poisson{\mu^{(1,2,1)},\gamma}}}\\
&=&\Poisson{\Poisson{Y^1,\Poisson{X^1,\mu^{(1,2,1)}}},\gamma}+\Poisson{\Poisson{X^1,\mu^{(1,2,1)}},\Poisson{Y^1,\gamma}}\\
&&-\Poisson{\Poisson{Y^1,\mu^{(1,2,1)}},\Poisson{X^1,\gamma}}+\Poisson{\mu^{(1,2,1)},\Poisson{Y^1,\Poisson{X^1,\gamma}}}\\
&=&-\Poisson{\gamma,\Poisson{Y^1,\Poisson{X^1,\mu^{(1,2,1)}}}}+\Poisson{\Poisson{\gamma,Y^1},\Poisson{X^1,\mu^{(1,2,1)}}}-\Poisson{\Poisson{\gamma,X^1},\Poisson{Y^1,\mu^{(1,2,1)}}}\\
&=&\delta_*{[X^1,Y^1]}_S-{[X^1,\delta_*(Y^1)]}_S-{[\delta_*(X^1),Y^1]}_S,
\end{eqnarray*}
which implies that \eqref{L2bialgebroid1} holds for all $ X^1, Y^1\in\Gamma(A_{-1})$. Similarly, we can show that \eqref{L2bialgebroid1} holds for all $ X, Y\in\Gamma(\huaA[-3])$.
}

For all $ X, Y\in\Gamma(\huaA),$ by the graded Jacobi identity and \eqref{eq:defiSB}, we have
\begin{eqnarray*}
&&\Poisson{\Poisson{\Poisson{\mu^{(1,2,1)},\gamma},X},Y}\\
&=&\Poisson{\Poisson{\mu^{(1,2,1)},\Poisson{\gamma,X}},Y}+(-1)^{(|X|-3)}\Poisson{\Poisson{\Poisson{\mu^{(1,2,1)},X},\gamma},Y}\\
&=&\Poisson{\mu^{(1,2,1)},\Poisson{\Poisson{\gamma,X},Y}}+(-1)^{(|Y|-3)(|X|)}\Poisson{\Poisson{\mu^{(1,2,1)},Y},\Poisson{\gamma,X}}\\
&&+(-1)^{(|X|-3)}\Poisson{\Poisson{\mu^{(1,2,1)},X},\Poisson{\gamma,Y}}+(-1)^{(|X|+|Y|)}\Poisson{\Poisson{\Poisson{\mu^{(1,2,1)},X},Y},\gamma}\\
&=&\Poisson{\Poisson{\mu^{(1,2,1)},\Poisson{\gamma,X}},Y}+(-1)^{(|X|-3)}\Poisson{\Poisson{\mu^{(1,2,1)},X},\Poisson{\gamma,Y}}+\Poisson{\gamma,\Poisson{\Poisson{\mu^{(1,2,1)},X},Y}}\\
&=&-\delta_*{[X,Y]}_S-{[\delta_*(X),Y]}_S+(-1)^{|X|}{[X,\delta_*(Y)]}_S.
\end{eqnarray*}
Thus, $\Poisson{\gamma,\mu^{(1,2,1)}}=0$ if and only if \eqref{L2bialgebroid1} holds. Similarly, we can show that $\Poisson{\mu,\gamma^{(1,1,2)}}=0$   if and only if \eqref{L2bialgebroid2} holds.
 We finish the proof.\qed\vspace{3mm}

\emptycomment{
For $X=X^1\in\Gamma(A_{-1}),Y^2\in\Gamma(A_{-2})$, we have
\begin{eqnarray*}
&&\Poisson{Y^2,\Poisson{X^1,\Poisson{\mu^{(1,2,1)},\gamma}}}\\
&=&\Poisson{\Poisson{Y^2,\Poisson{X^1,\mu^{(1,2,1)}}},\gamma}+\Poisson{\Poisson{X^1,\mu^{(1,2,1)}},\Poisson{Y^2,\gamma}}\\
&&-\Poisson{\Poisson{Y^2,\mu^{(1,2,1)}},\Poisson{X^1,\gamma}}-\Poisson{\mu^{(1,2,1)},\Poisson{Y^2,\Poisson{X^1,\gamma}}}\\
&=&-\Poisson{\gamma,\Poisson{Y^2,\Poisson{X^1,\mu^{(1,2,1)}}}}+\Poisson{\Poisson{\gamma,Y^2},\Poisson{X^1,\mu^{(1,2,1)}}}-\Poisson{\Poisson{\gamma,X^1},\Poisson{Y^2,\mu^{(1,2,1)}}}\\
&=&\delta_*{[X^1,Y^2]}_S-{[X^1,\delta_*(Y^2)]}_S-{[\delta_*(X^1),Y^2]}_S,
\end{eqnarray*}
which implies that $\delta_*{[X^1,Y^2]}_S={[X^1,\delta_*(Y^2)]}_S+{[\delta_*(X^1),Y^2]}_S$ holds. Thus \eqref{L2bialgebroid1} follows immediately.
}

At the end of this section, we give some useful formulas that will be used in the next section.
\begin{pro}
Let $(\huaA,\huaA^*[3])$ be a split Lie $2$-bialgebroid. Then we have
\begin{eqnarray*}
\label{eq:properL2A1}\huaL^2_{\dM f}X^1=-\Schouten{X^1,\dM_{\ast}f},\quad
\label{eq:properL2A2}L^2_{\dM_* f}\alpha^2=-\Schouten{\alpha^2,\dM f},\quad \forall X^1\in\Gamma(A_{-1}),~\alpha^2\in\Gamma(A^*_{-2}),~f\in\CWM.
\end{eqnarray*}
\end{pro}
\pf For all $X,Y\in\Gamma(\huaA[-3])$ and $f\in\CWM$, by \eqref{L2bialgebroid1}, we have
\begin{eqnarray*}
\dM_*\Schouten{X,fY}&=&\dM_*(f\Schouten{X,Y})+\dM_*(a(X)(f)Y)\\
&=&f\dM_*\Schouten{X,Y}+\dM_*(f)\odot\Schouten{X,Y}+\big(\dM_*(a(X)(f))\odot Y+a(X)(f)\dM_*(Y)\big).
\end{eqnarray*}
On the other hand, we have
\begin{eqnarray*}
\dM_*\Schouten{X,fY}&=&-{[\dM_*(X),fY]}_S+(-1)^{|X|}{[X,\dM_*(fY)]}_S\\
&=&-f{[\dM_*(X),Y]}_S-(-1)^{\degree{X}+1}\iota_{\dM f}\dM_*{X}\odot Y+(-1)^{\degree{X}}{\Schouten{X,\dM_{*}f}}\odot Y+\dM_*(f)\odot\Schouten{X,Y}\\
&&+(-1)^{\degree{X}}f\Schouten{X,\dM_{*}Y}+a(X)(f)\dM_*(Y).
\end{eqnarray*}
Therefore, we have
$$\dM_*(a(X)(f))\odot Y=(-1)^{\degree{X}}\iota_{\dM f}\dM_*{X}\odot Y+(-1)^{\degree{X}}{\Schouten{X,\dM_{*}f}}\odot Y.$$
For $X=X^1$, we have
$\dM_*(a(X^1)(f))-\iota_{\dM f}\dM_*{X^1}=-{\Schouten{X^1,\dM_{*}f}},$
which implies that
$\huaL^2_{\dM f}X^1=-\Schouten{X^1,\dM_{\ast}f}.$

The other one can be proved similarly. We omit details.\qed

\section{Manin triples of split Lie 2-algebroids}\label{sec:M}
The notion of a $\LWX$ 2-algebroid (named
after Courant-Liu-Weinstein-Xu) was introduced in \cite{LiuSheng} as the categorification of a Courant algebroid \cite{lwx,Roytenbergphdthesis}.

\begin{defi}\label{defi:Courant-2 algebroid}
A {\bf $\LWX$ $2$-algebroid} is a graded vector bundle $\huaE=E_{-2}\oplus E_{-1}$ over $M$ equipped with a non-degenerate graded symmetric bilinear form 
 $S$ on $\huaE$, a bilinear operation $\diamond:\Gamma(E_{i})\times \Gamma(E_{j})\longrightarrow \Gamma(E_{i+j+1})$, $-3\leq i+j\leq -2$, which is skewsymmetric  on $\Gamma(E_{-1})\times \Gamma(E_{-1})$, an $E_{-2}$-valued $3$-form $\Omega$ on $E_{-1}$, two bundle maps $\partial:E_{-2}\longrightarrow E_{-1}$ and $\rho:E_{-1}\longrightarrow TM$, such that $E_{-2}$ and $E_{-1}$ are isotropic and the following axioms are satisfied:
\begin{itemize}
\item[$\rm(i)$]$(\Gamma(E_{-2}),\Gamma(E_{-1}),\partial,\diamond,\Omega)$ is a Leibniz $2$-algebra;
\item[$\rm(ii)$]for all $e^1\in\Gamma(E_{-1}),~e^2\in\Gamma(E_{-2}) $, we have $e^1\diamond e^2-e^2\diamond e^1=\huaD S(e^1,e^2)$, where $\huaD:\CWM\longrightarrow \Gamma(E_{-2})$ is defined by
$$
 S(\huaD f,e^1)=\rho(e^1)(f),\quad \forall~e^1\in\Gamma(E_{-1});
$$
\item[$\rm(iii)$]for all $e^2_1,e^2_2\in\Gamma(E_{-2})$, we have $S( \partial(e^2_1),e^2_2)=S(e^2_1,\partial(e^2_2))$;
\item[$\rm(iv)$]for all $e_1,e_2,e_3\in\Gamma(\huaE)$, we have $\rho(e_1)S( e_2,e_3)=S( e_1\diamond e_2,e_3)+S(e_2,e_1\diamond e_3)$;
\item[$\rm(v)$]for all $e^1_1,e^1_2,e^1_3,e^1_4\in\Gamma(E_{-1})$, we have $S(\Omega(e^1_1,e^1_2,e^1_3),e^1_4)=-S(e^1_3,\Omega(e^1_1,e^1_2,e^1_4))$.
 \end{itemize}
\end{defi}
 Denote a $\LWX$ 2-algebroid  by $(E_{-2},E_{-1},\partial,\rho,S,\diamond,\Omega)$, or simply by $\huaE$.  The following lemma lists some properties of a $\LWX$ 2-algebroid.
 \begin{lem}
Let $(E_{-2},E_{-1},\partial,\rho,S,\diamond,\Omega)$ be a $\LWX$ $2$-algebroid. Then for all $e_1,e_2\in\Gamma(\huaE)$, $e^1, e^1_1, e^1_2 \in\Gamma(E_{-1})$ and $ f\in C^\infty(M)$, we have
$$\begin{array}{rclrcl}
e_1\diamond fe_2&=&f(e_1\diamond e_2)+\rho(e_1)(f)e_2,&
(fe_1)\diamond e_2&=&f(e_1\diamond e_2)+\rho(e_2)(f)e_1+ S(e_1,e_2)\huaD f,\\
 \rho\diamond \partial&=&0,&
\partial\diamond\huaD&=&0,\\
e^1\diamond \huaD f&=&\huaD S(e^1,\huaD f),&
 \huaD f\diamond e^1&=&0.
 \end{array}
$$
\end{lem}
 \emptycomment{\begin{defi}
  A {\bf Maurer-Cartan element} in an $L_\infty$-algebra is a degree $1$ element $\alpha$   satisfying the Maurer-Cartan equation:
\begin{equation}
  \sum_i\frac{(-1)^i}{i!}l_i(\alpha,\cdots,\alpha)=0.
\end{equation}
 \end{defi}
 }

\begin{defi}\label{defi:sDirac}
Let $\huaE=(E_{-2},E_{-1},\partial,\rho,S,\diamond,\Omega)$ be a $\LWX$ $2$-algebroid.

 \begin{itemize}
 \item[\rm(a)]A graded subbundle $L=L_{-2}\oplus L_{-1}$ of $\huaE$ is called {\bf isotropic} if
 $S(X,Y)=0,$ for all $X,Y\in\Gamma(L).$

  \item[\rm(b)] A graded subbundle $L=L_{-2}\oplus L_{-1}$ of $\huaE$  is called {\bf integral} if
\begin{itemize}
  \item[\rm(i)] $\partial\big(\Gamma(L_{-2})\big)\subseteq \Gamma(L_{-1})$;

   \item[\rm(ii)]$\Gamma(L)$ is closed under the operation $\diamond$;

   \item[\rm(iii)]$\Omega\big(\Gamma(L_{-1}),\Gamma(L_{-1}),\Gamma(L_{-1})\big)\subseteq \Gamma(L_{-2}).$
\end{itemize}
 \item[\rm(c)] A maximal isotropic and integral graded  subbundle $L$ of $\huaE$ is called a  {\bf strict Dirac structure} of a $\LWX$ $2$-algebroid.
 \end{itemize}
\end{defi}

The following proposition follows immediately from the definition.
\begin{pro}\label{pro:Dirac structures}
Let $L$ be a strict Dirac structure of a $\LWX$ $2$-algebroid $\huaE=(E_{-2},E_{-1},\partial,\rho,S,\diamond,\Omega)$. Then $(L_{-2},L_{-1},\partial\mid_{L},\diamond\mid_{L},\Omega\mid_{L},\rho\mid_{L})$ is a split Lie $2$-algebroid.
\end{pro}

\begin{defi}
  A {\bf Manin triple  of split Lie $2$-algebroids} $(\huaE;\huaA,\huaB)$ consists of a  $\LWX$ $2$-algebroid $\huaE=(E_{-2},E_{-1},\partial,\rho,S,\diamond,\Omega)$ and two transversal strict Dirac structures $\huaA$ and $\huaB$.
\end{defi}

\begin{thm}
  There is a one-to-one correspondence between Manin triples  of split Lie $2$-algebroids and split Lie $2$-bialgebroids.
\end{thm}
\pf It follows from the following Proposition \ref{thm:Lie2biC2} and Proposition \ref{pro:C2L2}.\qed\vspace{3mm}

Assume that $\huaA=(A_{-2},A_{-1},l_1,l_2,l_3,a)$ is a split Lie $2$-algebroid with structure function $\mu$ and $\huaA^*[3]=(A_{-1}^*[3],A_{-2}^*[3],\frkl_1,\frkl_2,\frkl_3,\frka)$ a split Lie $2$-algebroid with the structure function $\gamma$.
Let $E_{-1}=A_{-1}\oplus A^*_{-2}$ and $E_{-2}=A_{-2}\oplus A^*_{-1}$ and $\huaE=E_{-2}\oplus E_{-1}$.

Let $\partial:E_{-2}\longrightarrow E_{-1}$ and $\rho:E_{-1}\longrightarrow TM$ be bundle maps defined by
\begin{eqnarray}
\label{eq:parbi}\partial(X^2+\alpha^1)&=&l_1(X^2)+\frkl_1(\alpha^1),\label{eq:L2Bpartial}\\
\rho(X^1+\alpha^2)(f)&=&a(X^1)(f)+\frka(\alpha^2)(f)\label{eq:L2Banchor}.
\end{eqnarray}
On $\Gamma(\huaE)$, there is a natural symmetric bilinear form $(\cdot,\cdot)_+$ given by
\begin{equation} \label{eq:naturalsymform}
(X^1+\alpha^2+X^2+\alpha^1,Y^1+\beta^2+Y^2+\beta^1)_+=\langle X^1,\beta^1 \rangle+\langle Y^1,\alpha^1 \rangle+\langle X^2,\beta^2 \rangle+\langle Y^2,\alpha^2 \rangle.
\end{equation}

On $\Gamma(\huaE)$, we introduce operations $\diamond:E_i\times E_j\longrightarrow E_{i+j+1}$, $-3\leq i+j\leq -2$, by
\begin{eqnarray}\label{eq:L2Bbracket0}
\left\{\begin{array}{rcl}
(X^1+\alpha^2)\diamond(Y^1+\beta^2)&=&l_2(X^1,Y^1)+L^1_{X^1}\beta^2-L^1_{Y^1}\alpha^2 +\frkl_2(\alpha^2,\beta^2)+\huaL^1_{\alpha^2}Y^1-\huaL^1_{\beta^2}X^1,\\
(X^1+\alpha^2)\diamond(X^2+\alpha^1)&=&l_2(X^1,X^2)+L^1_{X^1}\alpha^1+\iota_{X^2}\dM(\alpha^2) +\frkl_2(\alpha^2,\alpha^1)+\huaL^1_{\alpha^2}X^2+\iota_{\alpha^1}\dM_*(X^1), \\
(X^2+\alpha^1)\diamond(X^1+\alpha^2)&=&l_2(X^2,X^1)+L^2_{X^2}\alpha^2+\iota_{X^1}\dM(\alpha^1) +\frkl_2(\alpha^1,\alpha^2)+\huaL^2_{\alpha^1}X^1+\iota_{\alpha^2}\dM_*(X^2).
\end{array}\right.
\end{eqnarray}

\emptycomment{\begin{eqnarray}
(X^1+\alpha^2)\diamond(Y^1+\beta^2)&=&-\Poisson{\Poisson{\mu^{(1,2,1)}+\gamma^{(1,1,2)},X^1+\alpha^2},Y^1+\beta^2}\\
\nonumber&=&l_2(X^1,Y^1)+L^1_{X^1}\beta^2-L^1_{Y^1}\alpha^2 +\frkl_2(\alpha^2,\beta^2)\\
\nonumber&&+\huaL^1_{\alpha^2}Y^1-\huaL^1_{\beta^2}X^1,\\
(X^1+\alpha^2)\diamond(X^2+\alpha^1)&=&-\Poisson{\Poisson{\mu^{(1,2,1)}+\gamma^{(1,1,2)},X^1+\alpha^2},X^2+\alpha^1}\\
\nonumber&=&l_2(X^1,X^2)+L^1_{X^1}\alpha^1+\iota_{X^2}\dM(\alpha^2) +\frkl_2(\alpha^2,\alpha^1)\\
\nonumber&&+\huaL^1_{\alpha^2}X^2+\iota_{\alpha^1}\dM_*(X^1), \\
(X^2+\alpha^1)\diamond(X^1+\alpha^2)&=&-\Poisson{\Poisson{\mu^{(1,2,1)}+\gamma^{(1,1,2)},X^2+\alpha^1},X^1+\alpha^2}\\
\nonumber&=&l_2(X^2,X^1)+L^2_{X^2}\alpha^2+\iota_{X^1}\dM(\alpha^1) +\frkl_2(\alpha^1,\alpha^2)\\
\nonumber&&+\huaL^2_{\alpha^1}X^1+\iota_{\alpha^2}\dM_*(X^2).
\end{eqnarray}}

\emptycomment{The skew-symmetric of operation $\circ$ is given by
\begin{equation}\label{eq:L2Bbracket0}
\left\{\begin{array}{rcl}
\Courant{X^1+\alpha^2,Y^1+\beta^2}&=&l_2(X^1,Y^1)+L^1_{X^1}\beta^2-L^1_{Y^1}\alpha^2 +\frkl_2(\alpha^2,\beta^2)+\huaL^1_{\alpha^2}Y^1-\huaL^1_{\beta^2}X^1,\\
\Courant{X^1+\alpha^2,X^2+\alpha^1}&=&l_2(X^1,X^2)+L^1_{X^1}\alpha^1-L^2_{X^2}\alpha^2 +\half\delta\langle X^2,\alpha^2\rangle-\half\delta\langle X^1,\alpha^1\rangle\\
&&\frkl_2(\alpha^2,\alpha^1)+\huaL^1_{\alpha^2}X^2-\huaL^2_{\alpha^1}(X^1)-\half\delta_*\langle X^2,\alpha^2\rangle+\half\delta_*\langle X^1,\alpha^1\rangle.\\
\end{array}\right.
\end{equation}}

An $E_{-2}$-valued $3$-form $\Omega$ on $E_{-1}$ is defined by
\begin{eqnarray}\label{eq:L2B3-form}
\nonumber\Omega(X^1+\alpha^2,Y^1+\beta^2,Z^1+\gamma^2)&=&
l_3(X^1,Y^1,Z^1)+L^3_{X^1,Y^1}\gamma^2+L^3_{Y^1,Z^1}\alpha^2+L^3_{Z^1,X^1}\beta^2\\
&&+\frkl_3(\alpha^2,\beta^2,\gamma^2)+\huaL^3_{\alpha^2,\beta^2}Z^1+\huaL^3_{\beta^2,\gamma^2}X^1+\huaL^3_{\gamma^2,\alpha^2}Y^1,
\end{eqnarray}
for all $X^1,Y^1,Z^1\in\Gamma(A_{-1}),~X^2,Y^2\in\Gamma(A_{-2}),~\alpha^1,\beta^1\in\Gamma(A^*_{-1}),~\alpha^2,\beta^2,\gamma^2\in\Gamma(A_{-2}^*).$

It is proved in \cite{LiuSheng} that
\begin{pro}\label{thm:Lie2biC2}
Let $(\huaA,\huaA^*[3])$ be a split Lie $2$-bialgebroid. Then $(E_{-2},E_{-1},\partial,\rho,(\cdot,\cdot)_+,\diamond,\Omega)$ is a $\LWX$ $2$-algebroid, where  $E_{-1}=A_{-1}\oplus A^*_{-2}$, $E_{-2}=A_{-2}\oplus A^*_{-1}$, $\partial$ is given by \eqref{eq:L2Bpartial}, $\rho$ is given by \eqref{eq:L2Banchor}, $(\cdot,\cdot)_+$ is given by \eqref{eq:naturalsymform}, $\diamond$ is given by \eqref{eq:L2Bbracket0} and $\Omega$ is given by \eqref{eq:L2B3-form}.
\end{pro}

Conversely, we have
\begin{pro}\label{pro:C2L2}
 Let $\huaA$ and $\huaB$ be two transversal strict Dirac structures of   a $\LWX$ $2$-algebroid $\huaE=(E_{-2},E_{-1},\partial,\rho,S,\diamond,\Omega)$, i.e. $\huaE=\huaA\oplus \huaB$ as graded vector bundles. Then $(\huaA,\huaB)$ is a split Lie $2$-bialgebroid, where $\huaB$ is considered as the shifted dual bundle of $\huaA$ under the bilinear form $S$.
\end{pro}

The proof is a long and tedious calculation and we put it in Appendix.

\section{Weak Dirac structures and Maurer-Cartan elements}\label{sec:MC}
Recall that given a Lie algebroid $(A,[\cdot,\cdot],\rho)$, there is naturally a Courant algebroid $A\oplus A^*$. Given an element $\pi\in\Gamma(\wedge^2A)$ such that $[\pi,\pi]=0$, one can define a Lie bracket $[\cdot,\cdot]_\pi$ on $\Gamma(A^*)$ such that $(A^*,[\cdot,\cdot]_\pi,\rho\circ \pi^\sharp)$ is a Lie algebroid, which we denote by $A^*_\pi$. Moreover, $\pi^\sharp$ is a Lie algebroid morphism from $A^*_\pi$ to $A$ and $(A,A^*_\pi)$ is a triangular Lie bialgebroid (\cite{MackenzieX:1994}). The graph of $\pi^\sharp$, which we denote by $\huaG_\pi\subset A\oplus A^*$ is a Dirac structure of the Courant algebroid $A\oplus A^*$. In this section, we generalize the above story to Lie 2-algebroids. First we introduce the notion of a weak Dirac structure of a $\LWX$ $2$-algebroid. Then we study Maurer-Cartan elements of the homotopy Poisson algebra associated to a Lie 2-algebroid $\huaA$ given in Section \ref{sec:H}. We show that a Maurer-Cartan element  gives rise to a split Lie 2-algebroid structure on the shifted dual bundle $\huaA^*[3]$ as well as a morphism from $\huaA^*[3]$ to $\huaA$. We also study Maurer-Cartan elements associated to a Lie 2-bialgebroid and show that the graph of such a Maurer-Cartan element is a weak Dirac structure of the corresponding $\LWX$ $2$-algebroid. Finally we give various examples including the string Lie 2-algebra, integrable distributions and left-symmetric algebroids.

\begin{defi}\label{defi:wDirac}
A split Lie $2$-algebroid $( L_{-2},L_{-1},l_1,l_2,l_3,a)$ is called a {\bf weak Dirac structure} of a $\LWX$ $2$-algebroid $(E_{-2},E_{-1},\partial,\rho,S,\diamond,\Omega)$ if there exists bundle maps $F_1:L_{-1}\longrightarrow E_{-1}$, $F_2:L_{-2}\longrightarrow E_{-2}$ and $F_3:\wedge^2L_{-1}\longrightarrow E_{-2}$ such that
\begin{itemize}
  \item[\rm(i)] $F_1$ and $F_2$ are injective such that the image $\img(F_2)\oplus\img(F_1)$ is a maximal isotropic graded subbundle of $E_{-2}\oplus E_{-1}$;
  \item[\rm(ii)] $(F_1,F_2,F_3)$ is a morphism from the Lie $2$-algebra $(\Gamma(L_{-2}),\Gamma(L_{-1}),l_1,l_2,l_3)$ to the Leibniz $2$-algebra $(\Gamma(E_{-2}),\Gamma(E_{-1}),\partial,\diamond,\Omega)$;

  \item[\rm(iii)] $\rho\circ F_1=a.$
\end{itemize}
\end{defi}

It is obvious that a strict Dirac structure $L$ given in Definition \ref{defi:sDirac} is a weak Dirac structure, in which $F_1$ and $F_2$ are inclusion maps and $F_3=0$.

\subsection{Maurer-Cartan elements associated to a split Lie 2-algebroid}
\begin{defi}
 Let $\huaA$ be a split Lie $2$-algebroid. A {\bf Maurer-Cartan element} of the associated  homotopy Poisson algebra $(\Sym(\huaA[-3]),[\cdot]_S,[\cdot,\cdot]_S,[\cdot,\cdot,\cdot]_S)$ given in Theorem \ref{thm:hpa} is an element $$m\in \Sym^3(\huaA[-3])=A_{-1}[-3]\odot A_{-2}[-3]\oplus \wedge^3 A_{-2}[-3]$$ such that
  \begin{equation}
    [m]_S+\half[m,m]_S+\frac{1}{6}[m,m,m]_S=0.
  \end{equation}
\end{defi}

An element $m\in\Sym^3(\huaA[-3])$ consists of an $H\in \Gamma(A_{-1}\odot A_{-2})$ and  a  $K\in\Gamma(\wedge^3A_{-2})$. For $H\in \Gamma(A_{-1}\odot A_{-2})$, define $H^\natural:\Gamma(A^*_{-1})\rightarrow \Gamma(A_{-2})$ and $H^\sharp:\Gamma(A^*_{-2})\rightarrow \Gamma(A_{-1})$ by
\begin{eqnarray*}
 \langle H^\natural(\alpha^1),\alpha^2\rangle=H(\alpha^1,\alpha^2),\quad
 \langle H^\sharp(\alpha^2),\alpha^1\rangle=H(\alpha^2,\alpha^1),\quad\forall~\alpha^1\in\Gamma(A^*_{-1}),~\alpha^2\in\Gamma(A^*_{-2}).
\end{eqnarray*}
We have
$
  \Poisson{\alpha^1,H}=H^\natural(\alpha^1),$ $\Poisson{\alpha^2,H}=-H^\sharp(\alpha^2).
$
For $K\in\Gamma(\wedge^3A_{-2}) $, define $K^\flat:\wedge^2A_{-2}^*\longrightarrow A_{-2}$ by
\begin{equation}\label{eq:Kb}
 \langle K^\flat(\alpha^2,\beta^2),\gamma^2\rangle=K(\alpha^2,\beta^2,\gamma^2),\quad\forall~\alpha^2,\beta^2,\gamma^2\in\Gamma(A^*_{-2}).
\end{equation}
It is not hard to see that
\begin{eqnarray*}
\Poisson{\Poisson{K,\alpha^2},\beta^2}=-K^\flat(\alpha^2,\beta^2),\quad
\Poisson{\Poisson{\Poisson{K,\alpha^2},\beta^2},\gamma^2}=-K(\alpha^2,\beta^2,\gamma^2).
\end{eqnarray*}

Let $\mu$ be the degree 4 function on $T^*[3](A_{-1}\oplus A_{-2})$ corresponding to the split Lie 2-algebroid $\huaA$. For $H\in \Gamma(A_{-1}\odot A_{-2})$ and $K\in\Gamma(\wedge^3A_{-2}) $, define
\begin{eqnarray}\label{eq:r-matrix}
\left\{\begin{array}{rcl}
\frkl^H_1(\alpha^1)&=&-\Poisson{\mu^{(2,1,1)},\alpha^1},\\
{\frkl^H_2(\alpha^2,\beta^2)}&=&-\Poisson{\Poisson{\Poisson{\mu^{(1,2,1)},H},\alpha^2},\beta^2},\\
{\frkl^H_2(\alpha^2,\beta^1)}&=&-\Poisson{\Poisson{\Poisson{\mu^{(1,2,1)},H},\alpha^2},\beta^1},\\
\frkl^{H,K}_3(\alpha^2,\beta^2,\gamma^2)&=&-\Poisson{\Poisson{\Poisson{\Poisson{\mu^{(1,2,1)},K},\alpha^2},\beta^2},\gamma^2}
-\half\Poisson{\Poisson{\Poisson{\Poisson{\Poisson{\mu^{(0,3,1)},H},H},\alpha^2},\beta^2},\gamma^2},\\
\frka_H(\alpha^2)(f)&=&-\Poisson{\Poisson{\Poisson{\mu^{(1,2,1)},H},\alpha^2},f},
\end{array}\right.
\end{eqnarray}
for all $\alpha^2,\beta^2,\gamma^2\in\Gamma(A^*_{-2}),~\alpha^1,\beta^1\in\Gamma(A^*_{-1}),~f\in C^\infty(M)$.
We use the Lie derivatives introduced in   Section \ref{sec:D} to give a precise description of the above operations.
\begin{lem}\label{lem:formulalH}
For all $\alpha^2,\beta^2,\gamma^2\in\Gamma(A^*_{-2}),~\beta^1\in\Gamma(A^*_{-1})$, we have
\begin{eqnarray}
\frkl^H_1&=&l^*_1,\label{eq:r-matrix1}\\
{\frkl^H_2(\alpha^2,\beta^2)}&=&L^1_{H^\sharp(\alpha^2)}\beta^2-L^1_{H^\sharp(\beta^2)}\alpha^2,\label{eq:r-matrix2}\\
{\frkl^H_2(\alpha^2,\beta^1)}&=&L^1_{H^\sharp(\alpha^2)}\beta^1-L^2_{H^\natural(\beta^1)}\alpha^2-\dM H(\alpha^2,\beta^1), \label{eq:r-matrix3}\\
\frkl^{H,K}_3(\alpha^2,\beta^2,\gamma^2)&=&-L^2_{K^\flat(\alpha^2,\beta^2)}\gamma^2-L^2_{K^\flat(\gamma^2,\alpha^2)}\beta^2-L^2_{K^\flat(\beta^2,\gamma^2)}\alpha^2-2\dM K(\alpha^2,\beta^2,\gamma^2)\nonumber\\
&&+L^3_{H^\sharp(\alpha^2),H^\sharp(\beta^2)}\gamma^2+L^3_{H^\sharp(\beta^2),H^\sharp(\gamma^2)}\alpha^2+L^3_{H^\sharp(\gamma^2),H^\sharp(\alpha^2)}\beta^2\label{eq:r-matrix4},\\
\frka_H&=&a\circ H^\sharp\label{eq:r-matrix0}.
\end{eqnarray}
\end{lem}
We need the following preparation before we give the main result in this subsection.

\begin{lem}\label{lem:MCformulars}
 Let $\huaA=(A_{-2},A_{-1},l_1,l_2,l_3,a)$ be a split Lie $2$-algebroid with the structure function $\mu$ and $H\in\Gamma(A_{-1}\odot A_{-2}),~K\in\Gamma(\wedge^3A_{-2})$. Then for all $\alpha^2,\beta^2,\gamma^2\in\Gamma(A^*_{-2}),~\alpha^1,\beta^1\in\Gamma(A^*_{-1})$, we have
\begin{eqnarray}
\label{eq:HS}{[H]_S}(\alpha^1)&=&l_1(H^\natural(\alpha^1))+H^\sharp(\frkl^H_1(\alpha^1)),\\
{[K]}_S(\alpha^2,\beta^2)&=&l_1K^\flat(\alpha^2,\beta^2),\\
{[K]}_S(\alpha^2,\beta^1)&=&-K^\flat(\alpha^2,\frkl^H_1(\beta^1)),\\
 \half{[H,H]}_S(\alpha^2,\beta^2,\cdot)&=&H^\sharp(\frkl^H_2(\alpha^2,\beta^2))-l_2(H^\sharp(\alpha^2),H^\sharp(\beta^2)),\label{eq:MC-formular1}\\
 \half{[H,H]}_S(\alpha^2,\beta^1,\cdot)&=&H^\natural(\frkl^H_2(\alpha^2,\beta^1))- l_2(H^\sharp(\alpha^2),H^\natural(\beta^1)),\label{eq:MC-formular2}\\
 {[H,K]_S}(\alpha^2,\beta^2,\gamma^2)
 &=&K^\flat(\frkl^H_2(\alpha^2,\beta^2),\gamma^2)+K^\flat(\frkl^H_2(\gamma^2,\alpha^2),\beta^2)+K^\flat(\frkl^H_2(\beta^2,\gamma^2),\alpha^2)\nonumber\\
 &&-l_2(H^\sharp(\alpha^2),K^\flat(\beta^2,\gamma^2))-l_2(H^\sharp(\gamma^2),K^\flat(\alpha^2,\beta^2))-l_2(H^\sharp(\beta^2),K^\flat(\gamma^2,\alpha^2))\nonumber\\
\label{eq:MC-formular3} &&-H^\natural(L^2_{K^\flat(\alpha^2,\beta^2)}\gamma^2+L^2_{K^\flat(\gamma^2,\alpha^2)}\beta^2+L^2_{K^\flat(\beta^2,\gamma^2)}\alpha^2+2\dM K(\alpha^2,\beta^2,\gamma^2)),\\
 \frac{1}{6}[H,H,H]_S(\alpha^2,\beta^2,\gamma^2)&=&H^\natural(L^3_{H^\sharp(\alpha^2),H^\sharp(\beta^2)}\gamma^2+L^3_{H^\sharp(\beta^2),H^\sharp(\gamma^2)}\alpha^2+L^3_{H^\sharp(\gamma^2),H^\sharp(\alpha^2)}\beta^2)\nonumber\\
\label{eq:MC-formular4} &&+l_3(H^\sharp(\alpha^2),H^\sharp(\beta^2),H^\sharp(\gamma^2)).
\end{eqnarray}
 \end{lem}
 \pf
 By the graded Jacobi identity for the canonical Poisson bracket on $T^*[3](A_{-1}\oplus A_{-2})$, we have
 \begin{eqnarray*}
   {[H]_S}(\alpha^1)&=&\Poisson{\Poisson{\mu^{(2,1,1)},H},\alpha^1}=-\Poisson{\mu^{(2,1,1)},H^\natural(\alpha^1)}+\Poisson{\Poisson{\mu^{(2,1,1)},\alpha^1},H}\\
   &=&l_1(H^\sharp(\alpha^1))+H^\sharp(\frkl^H_1(\alpha_1)),\\
   {[K]}_S(\alpha^2,\beta^2)&=&\Poisson{\Poisson{\Poisson{\mu^{(2,1,1)},K},\alpha^2},\beta^2}=\Poisson{\Poisson{\mu^{(2,1,1)},\Poisson{K,\alpha^2}},\beta^2}\\
   &=&\Poisson{\mu^{(2,1,1)},\Poisson{\Poisson{K,\alpha^2},\beta^2}}=-\Poisson{\mu^{(2,1,1)},K^\flat(\alpha^2,\beta^2)}=l_1K^\flat(\alpha^2,\beta^2),\\
   {[K]}_S(\alpha^2,\beta^1)&=&-\Poisson{\Poisson{\Poisson{\mu^{(2,1,1)},K},\alpha^2},\beta^1}=-\Poisson{\Poisson{\mu^{(2,1,1)},\Poisson{K,\alpha^2}},\beta^1}\\
   &=&-\Poisson{\Poisson{\mu^{(2,1,1)},\beta^1},\Poisson{K,\alpha^2}}=\Poisson{\frkl^H_1(\beta^1),\Poisson{K,\alpha^2}}=-K^\flat(\alpha^2,\frkl^H_1(\beta^1)),\\
 {[H,H]}_S(\alpha^2,\beta^2,\cdot)&=&\Poisson{\Poisson{\Poisson{\Poisson{\mu^{(1,2,1)},H},H},\alpha^2},\beta^2}\\
 &=&\Poisson{\Poisson{\Poisson{\mu^{(1,2,1)},H},\Poisson{H,\alpha^2},\beta^2}+\Poisson{\Poisson{\Poisson{\Poisson{\mu^{(1,2,1)},H},H},\alpha^2}},\beta^2}\\
 &=&-\Poisson{\Poisson{\Poisson{\mu^{(1,2,1)},H},\beta^2},\Poisson{H,\alpha^2}}+\Poisson{\Poisson{\Poisson{\mu^{(1,2,1)},H},\Poisson{H,\alpha^2}},\beta^2}\\
 &&+\Poisson{\Poisson{\Poisson{\Poisson{\mu^{(1,2,1)},H},\alpha^2},H},\beta^2}\\
 &=&-\Poisson{\Poisson{\mu^{(1,2,1)},\Poisson{H,\beta^2}},\Poisson{H,\alpha^2}}-\Poisson{\Poisson{\Poisson{\mu^{(1,2,1)},\beta^2},H},\Poisson{H,\alpha^2}}\\
 &&+\Poisson{\Poisson{\mu^{(1,2,1)},\Poisson{H,\alpha^2}},\Poisson{H,\beta^2}}+\Poisson{\Poisson{\Poisson{\mu^{(1,2,1)},\alpha^2},H},\Poisson{H,\beta^2}}\\
 &&+\Poisson{\Poisson{\Poisson{\Poisson{\mu^{(1,2,1)},H},\alpha^2},\beta^2},H}\\
 &=&2\Poisson{\Poisson{\mu^{(1,2,1)},\Poisson{H,\alpha^2}},\Poisson{H,\beta^2}}+\Poisson{\Poisson{\Poisson{\mu^{(1,2,1)},\Poisson{H,\alpha^2}},\beta^2},H}\\
 &&-\Poisson{\Poisson{\Poisson{\mu^{(1,2,1)},\Poisson{H,\beta^2}},\alpha^2},H}+\Poisson{\Poisson{\Poisson{\Poisson{\mu^{(1,2,1)},H},\alpha^2},\beta^2},H}\\
 &=&2H^\sharp\frkl^H_2(\alpha^2,\beta^2)-2l_2(H^\sharp(\alpha^2),H^\sharp(\beta^2)),
 \end{eqnarray*}
 which imply that \eqref{eq:HS}-\eqref{eq:MC-formular1} hold.

  By direct calculation, we have
 \begin{eqnarray*}
\langle H^\natural(\frkl^H_2(\alpha^2,\beta^1))-l_2(H^\sharp(\alpha^2),H^\natural(\beta^1)),\beta^2\rangle=\langle H^\sharp(\frkl^H_2(\alpha^2,\beta^2))-l_2(H^\sharp(\alpha^2),H^\sharp(\beta^2)),\beta^1\rangle.
 \end{eqnarray*}
 By \eqref{eq:MC-formular1}, \eqref{eq:MC-formular2} follows immediately.

\eqref{eq:MC-formular3} and \eqref{eq:MC-formular4} can be proved similarly. We omit the details.  \qed\vspace{3mm}

Now we are ready to give the main result in this subsection. Consider the following function\footnote{The function $\gamma_{H,K}$ is obtained in the following intrinsic way. The map $\{H+K,\cdot\}:C^\infty(\huaM)\lon C^\infty(\huaM)$ is an inner derivation of $C^\infty(\huaM)$. It follows that $e^{\{H+K,\cdot\}}$ is an automorphism of $C^\infty(\huaM)$. Thus, $\mu':=e^{\{H+K,\cdot\}}\mu$ is also a degree 4 function satisfying $\{\mu',\mu'\}=0$. $\gamma_{H,K}$ is exactly the projection of $\mu'$  to the subspace $C^{(2,1,1)}(\huaM)\oplus C^{(1,1,2)}(\huaM)\oplus C^{(0,1,3)}(\huaM)$ of $C^\infty(\huaM)$. See \cite{roytwist} for a similar discussion in Lie algebroids.} $\gamma_{H,K}$ of degree $4$ on $\huaM=T^*[3](A_{-1}\oplus A_{-2})$:
\begin{equation}\label{eq:gammaHK}
\gamma_{H,K}=\mu^{(2,1,1)}+\Poisson{\mu^{(1,2,1)},H}+\Poisson{\mu^{(1,2,1)},K}+\half\Poisson{\Poisson{\mu^{(0,3,1)},H},H}.
\end{equation}

\begin{thm}\label{thm:MCE1}
Let $\huaA=(A_{-2},A_{-1},l_1,l_2,l_3,a)$ be a split Lie $2$-algebroid and  $H+K$ is a Maurer-Cartan element of the associated homotopy Poisson algebra $(\Sym(\huaA[-3]),[\cdot]_S,[\cdot,\cdot]_S,[\cdot,\cdot,\cdot]_S)$, i.e.
\begin{eqnarray}\label{eq:MC}
[H+K]_S+\half[H+K,H+K]_S+\frac{1}{6}[H+K,H+K,H+K]_S=0.
\end{eqnarray}
Then we have
 \begin{itemize}
   \item[\rm(i)] $\huaA^*[3]=(A_{-1}^*[3],A_{-2}^*[3],\frkl^H_1,\frkl^H_2,\frkl^{H,K}_3,\frka_H)$ is a split Lie $2$-algebroid, where $\frkl^H_1,~ \frkl^H_2,~\frkl^{H,K}_3$ and $\frka_H$ are given by \eqref{eq:r-matrix1}-\eqref{eq:r-matrix0} respectively;
     \item[\rm(ii)]   $(H^\sharp,-H^\natural,-K^\flat)$ is a morphism from the split Lie $2$-algebroid $\huaA^*[3]=(A_{-1}^*[3],A_{-2}^*[3],\frkl^H_1,\frkl^H_2,\frkl^{H,K}_3,\frka_H)$ to the split Lie $2$-algebroid $\huaA=(A_{-2},A_{-1},l_1,l_2,l_3,a)$;
          \item[\rm(iii)] $(\huaA,\huaA^*[3])$ is a split Lie $2$-bialgebroid if and only if $$\Poisson{\mu^{(2,1,1)},\Poisson{\mu^{(0,3,1)},H}}=0,\quad\Poisson{\mu^{(1,2,1)},\Poisson{\mu^{(0,3,1)},H}}=0,\quad\Poisson{\mu^{(0,3,1)},\Poisson{\mu^{(2,1,1)},K}}=0.$$
 \end{itemize}
\end{thm}
\pf It is straightforward to deduce that \eqref{eq:MC} is equivalent to the following equations
\begin{eqnarray}
  [H]_S&=&0,\label{eq:MC-Cond 1}\\
  {[K]}_S+\half [H,H]_S&=&0,\label{eq:MC-Cond 2}\\
  {[H,K]_S}+\frac{1}{6}[H,H,H]_S&=&0.\label{eq:MC-Cond 3}
\end{eqnarray}

(i) To show that
$\huaA^*[3]=(A_{-1}^*[3],A_{-2}^*[3],\frkl^H_1,\frkl^H_2,\frkl^{H,K}_3,\frka_H)$ is a split Lie $2$-algebroid, we only need to prove $$\Poisson{\gamma_{H,K},\gamma_{H,K}}=0,$$
which is equivalent to the following equations:
\begin{eqnarray}
  \Poisson{\Poisson{\mu^{(1,2,1)},H},\mu^{(2,1,1)}}&=&0,\label{eq:MC-Lie 2-algebroid 1}\\
\Poisson{\Poisson{\mu^{(1,2,1)},H},\Poisson{\mu^{(1,2,1)},H}}+2\Poisson{\mu^{(2,1,1)},\Poisson{\mu^{(1,2,1)},K}}+\Poisson{\mu^{(2,1,1)},\Poisson{\Poisson{\mu^{(0,3,1)},H},H}}&=&0,\label{eq:MC-Lie 2-algebroid 2}\\
  \Poisson{\Poisson{\mu^{(1,2,1)},H},\Poisson{\mu^{(1,2,1)},K}}+\half \Poisson{\Poisson{\mu^{(1,2,1)},H},\Poisson{\Poisson{\mu^{(0,3,1)},H},H}}&=&0\label{eq:MC-Lie 2-algebroid 3}.
\end{eqnarray}

By \eqref{eq:MC-Cond 1} and the fact that $\Poisson{\mu^{(1,2,1)},\mu^{(2,1,1)}}=0$, we have
\begin{eqnarray*}
 \Poisson{\Poisson{\mu^{(1,2,1)},H},\mu^{(2,1,1)}}=\Poisson{\mu^{(1,2,1)},\Poisson{H,\mu^{(2,1,1)}}}+ \Poisson{\Poisson{\mu^{(1,2,1)},\mu^{(2,1,1)}},H}=0,
\end{eqnarray*}
which implies that \eqref{eq:MC-Lie 2-algebroid 1} holds.

By \eqref{eq:MC-Cond 1}, \eqref{eq:MC-Cond 2} and the fact that $\half\Poisson{\mu^{(1,2,1)},\mu^{(1,2,1)}}+\Poisson{\mu^{(2,1,1)},\mu^{(0,3,1)}}=0$, we have
\begin{eqnarray*}
&&\Poisson{\Poisson{\mu^{(1,2,1)},H},\Poisson{\mu^{(1,2,1)},H}}+2\Poisson{\mu^{(2,1,1)},\Poisson{\mu^{(1,2,1)},K}}+\Poisson{\mu^{(2,1,1)},\Poisson{\Poisson{\mu^{(0,3,1)},H},H}}\\
&=&\Poisson{\Poisson{\mu^{(1,2,1)},H},\Poisson{\mu^{(1,2,1)},H}}+\Poisson{\mu^{(1,2,1)},\Poisson{\Poisson{\mu^{(1,2,1)},H},H}}+\Poisson{\Poisson{\Poisson{\mu^{(2,1,1)},\mu^{(0,3,1)}},H},H}\\
&=&\Poisson{\Poisson{\half\Poisson{\mu^{(1,2,1)},\mu^{(1,2,1)}}+\Poisson{\mu^{(2,1,1)},\mu^{(0,3,1)}},H},H}=0,
\end{eqnarray*}
which implies that \eqref{eq:MC-Lie 2-algebroid 2} holds.

By \eqref{eq:MC-Cond 2} and \eqref{eq:MC-Cond 3}, we have
\begin{eqnarray*}
&&\Poisson{\Poisson{\mu^{(1,2,1)},H},\Poisson{\mu^{(1,2,1)},K}}\\
&=&\Poisson{\Poisson{\Poisson{\mu^{(1,2,1)},H},\mu^{(1,2,1)}},K}-\Poisson{\mu^{(1,2,1)},\Poisson{\Poisson{\mu^{(1,2,1)},H},K}}\\
&=&\Poisson{\Poisson{\half\Poisson{\mu^{(1,2,1)},\mu^{(1,2,1)}},H},K}-\Poisson{\mu^{(1,2,1)},\Poisson{\Poisson{\mu^{(1,2,1)},H},K}}\\
&=&-\Poisson{\Poisson{\Poisson{\mu^{(2,1,1)},\mu^{(0,3,1)}},H},K}-\Poisson{\mu^{(1,2,1)},\Poisson{\Poisson{\mu^{(1,2,1)},H},K}}\\
&=&-\Poisson{\Poisson{\mu^{(2,1,1)},K},\Poisson{\mu^{(0,3,1)},H}}-\Poisson{\mu^{(1,2,1)},\Poisson{\Poisson{\mu^{(1,2,1)},H},K}}\\
&=&\half\Poisson{\Poisson{\Poisson{\mu^{(1,2,1)},H},H},\Poisson{\mu^{(0,3,1)},H}}+\frac{1}{6}\Poisson{\mu^{(1,2,1)},\Poisson{\Poisson{\Poisson{\mu^{(0,3,1)},H},H},H}}.
\end{eqnarray*}
On the other hand, by the fact that $\Poisson{\mu^{(1,2,1)},\mu^{(0,3,1)}}=0$, we have
\begin{eqnarray*}
&& \frac{1}{6}\Poisson{\mu^{(1,2,1)},\Poisson{\Poisson{\Poisson{\mu^{(0,3,1)},H},H},H}}\\
 &=&-\half\Poisson{\Poisson{\Poisson{\mu^{(1,2,1)},H},H},\Poisson{\mu^{(0,3,1)},H}}-\half \Poisson{\Poisson{\mu^{(1,2,1)},H},\Poisson{\Poisson{\mu^{(0,3,1)},H},H}}.
\end{eqnarray*}
Therefore we have
$$\Poisson{\Poisson{\mu^{(1,2,1)},H},\Poisson{\mu^{(1,2,1)},K}}+\half \Poisson{\Poisson{\mu^{(1,2,1)},H},\Poisson{\Poisson{\mu^{(0,3,1)},H},H}}=0,$$
which implies that \eqref{eq:MC-Lie 2-algebroid 3} holds. Thus, $(A_{-1}^*[3],A_{-2}^*[3],\frkl^H_1,\frkl^H_2,\frkl^{H,K}_3,\frka_H=a\circ H^\sharp)$ is a split Lie $2$-algebroid.\vspace{2mm}

(ii) By \eqref{eq:MC-Cond 1} and Lemma \ref{lem:MCformulars}, we obtain
\begin{equation}\label{eq:Lie2morp1}
-l_1\circ H^\natural=H^\sharp\circ\frkl^H_1.
\end{equation}
By \eqref{eq:MC-Cond 2} and Lemma \ref{lem:MCformulars}, we have
\begin{eqnarray}
\label{eq:Lie2morp2}H^\sharp\frkl^H_2(\alpha^2,\beta^2)-l_2(H^\sharp(\alpha^2),H^\sharp(\beta^2))&=&-l_1K^\flat(\alpha^2,\beta^2),\\
\label{eq:Lie2morp3}(-H^\natural)\frkl^H_2(\alpha^2,\beta^1)-l_2(H^\sharp(\alpha^2),-H^\natural(\beta^1))&=&-K^\flat(\alpha^2,\frkl^H_1(\beta^1)).
\end{eqnarray}
By \eqref{eq:MC-Cond 3} and Lemma \ref{lem:MCformulars}, we have
\begin{eqnarray*}
&&K^\flat(\frkl^H_2(\alpha^2,\beta^2),\gamma^2)+K^\flat(\frkl^H_2(\gamma^2,\alpha^2),\beta^2)+K^\flat(\frkl^H_2(\beta^2,\gamma^2),\alpha^2)-l_2(H^\sharp(\alpha^2),K^\flat(\beta^2,\gamma^2))\nonumber\\
 &&-l_2(H^\sharp(\gamma^2),K^\flat(\alpha^2,\beta^2))-l_2(H^\sharp(\beta^2),K^\flat(\gamma^2,\alpha^2))+H^\natural\Big(-L^2_{K^\flat(\alpha^2,\beta^2)}\gamma^2-L^2_{K^\flat(\gamma^2,\alpha^2)}\beta^2\nonumber\\
 &&-L^2_{K^\flat(\beta^2,\gamma^2)}\alpha^2-2\dM K(\alpha^2,\beta^2,\gamma^2)+L^3_{H^\sharp(\alpha^2),H^\sharp(\beta^2)}\gamma^2+L^3_{H^\sharp(\beta^2),H^\sharp(\gamma^2)}\alpha^2+L^3_{H^\sharp(\gamma^2),H^\sharp(\alpha^2)}\beta^2\Big)\\
 &&+l_3(H^\sharp(\alpha^2),H^\sharp(\beta^2),H^\sharp(\gamma^2))=0,
\end{eqnarray*}
which implies that
\begin{eqnarray}
\nonumber&&-K^\flat(\frkl^H_2(\alpha^2,\beta^2),\gamma^2)-K^\flat(\frkl^H_2(\gamma^2,\alpha^2),\beta^2)-K^\flat(\frkl^H_2(\beta^2,\gamma^2),\alpha^2)-H^\natural(\frkl^{H,K}_3(\alpha^2,\beta^2,\gamma^2)=\\
\label{eq:Lie2morp4}&&l_3(H^\sharp(\alpha^2),H^\sharp(\beta^2),H^\sharp(\gamma^2))+l_2(H^\sharp(\alpha^2),-K^\flat(\beta^2,\gamma^2))+l_2(H^\sharp(\gamma^2),-K^\flat(\alpha^2,\beta^2))\\
\nonumber &&+l_2(H^\sharp(\beta^2),-K^\flat(\gamma^2,\alpha^2)).
\end{eqnarray}

Thus, $(H^\sharp,-H^\natural,-K^\flat)$ is a morphism from the split Lie $2$-algebroid $(A_{-1}^*[3],A_{-2}^*[3],\frkl^H_1,\frkl^H_2,\frkl^{H,K}_3,\frka_H)$ to the split Lie $2$-algebroid $(A_{-2},A_{-1},l_1,l_2,l_3,a)$.\vspace{2mm}

(iii) Note that $(\huaA,\huaA^*[3])$ is a split Lie $2$-bialgebroid if and only if
$$\Poisson{\mu+\gamma_{H,K}-\mu^{(2,1,1)},\mu+\gamma_{H,K}-\mu^{(2,1,1)}}=0.$$
Since $(\huaA,\mu)$ and $(\huaA^*[3],\gamma_{H,K})$ are split Lie $2$-algebroids, the above equality is equivalent to
\begin{eqnarray*}
\Poisson{\mu^{(1,2,1)},\Poisson{\mu^{(1,2,1)},H}}&=&0,\\
\Poisson{\mu^{(1,2,1)},\Poisson{\mu^{(1,2,1)},K}}+\half\Poisson{\mu^{(1,2,1)},\Poisson{\Poisson{\mu^{(0,3,1)}, H},H}}&=&0,\\
\Poisson{\Poisson{\mu^{(1,2,1)},H},\mu^{(0,3,1)}}&=&0.
\end{eqnarray*}
By the fact that $\half\Poisson{\mu^{(1,2,1)},\mu^{(1,2,1)}}+\Poisson{\mu^{(2,1,1)},\mu^{(0,3,1)}}=0$ and $\Poisson{H,\mu^{(2,1,1)}}=0$, we have
\begin{eqnarray*}
  \Poisson{\mu^{(1,2,1)},\Poisson{\mu^{(1,2,1)},H}} &=&-\Poisson{\mu^{(2,1,1)},\Poisson{\mu^{(0,3,1)},H}}.
\end{eqnarray*}
By the fact that $\Poisson{\mu^{(1,2,1)},\mu^{(0,3,1)}}=0$, we have
$$\Poisson{\Poisson{\mu^{(1,2,1)},H},\mu^{(0,3,1)}}=-\Poisson{\mu^{(1,2,1)},\Poisson{\mu^{(0,3,1)},H}}.$$
By the fact that $\Poisson{\mu^{(1,2,1)},\mu^{(0,3,1)}}=0$ and \eqref{eq:MC-Cond 2}, we have
\begin{eqnarray*}
&&\Poisson{\mu^{(1,2,1)},\Poisson{\mu^{(1,2,1)},K}}+\half\Poisson{\mu^{(1,2,1)},\Poisson{\Poisson{\mu^{(0,3,1)}, H},H}}\\
&=&\Poisson{\Poisson{\mu^{(1,2,1)},\Poisson{\mu^{(0,3,1)},H}},H}-2\Poisson{\mu^{(0,3,1)},\Poisson{\mu^{(2,1,1)},K}}.
\end{eqnarray*}
Thus, $(\huaA,\huaA^*[3])$ is a split Lie $2$-bialgebroid if and only if
$$\Poisson{\mu^{(2,1,1)},\Poisson{\mu^{(0,3,1)},H}}=0,\quad\Poisson{\mu^{(1,2,1)},\Poisson{\mu^{(0,3,1)},H}}=0,\quad\Poisson{\mu^{(0,3,1)},\Poisson{\mu^{(2,1,1)},K}}=0.$$
We finish the proof.\qed

\begin{rmk}
  According to Proposition \ref{thm:Lie2biC2}, given a Lie $2$-algebroid $(A_{-2},A_{-1},l_1,l_2,l_3,a)$, then $(A_{-2}\oplus A^*_{-1},A_{-1}\oplus A^*_{-2},\partial,\rho,(\cdot,\cdot)_+,\diamond,\Omega)$ is a $\LWX$ $2$-algebroid. It is natural to expect that the graph of $(H^\sharp, -H^\natural)$ is a ``Dirac structure''. However, it is straightforward to see that  the graph of $(H^\sharp, -H^\natural)$ is not closed under the operation $\diamond$ anymore. Thus, it is not a strict Dirac structure defined in Definition \ref{defi:sDirac}. But by Theorem \ref{thm:MCE1}, we can deduce that the graph of $(H^\sharp, -H^\natural)$ is a weak Dirac structure defined in Definition \ref{defi:wDirac}. We will prove this result for the more general case of  split Lie $2$-bialgebroids in  the next subsection.
\end{rmk}

\subsection{Weak Dirac structures and Maurer-Cartan elements associated to a split Lie 2-bialgebroid}
Let $\huaA=(A_{-2},A_{-1},l_1,l_2,l_3,a)$ be a split Lie $2$-algebroid with the structure function $\mu$ and $\huaA^*[3]=(A_{-1}^*[3],A_{-2}^*[3],\frkl_1,\frkl_2,\frkl_3,\frka)$ a split Lie $2$-algebroid with the structure function $\gamma$ such that $(\huaA,\huaA^*[3])$ is a split Lie $2$-bialgebroid. Let $H\in\Gamma(A_{-1}\odot A_{-2})$ and $K\in \Gamma(\wedge^3A_{-2})$. Define $\Lambda$ to be the degree $4$ function on $T^*[3](A_{-1}\oplus A_{-2})$ by\footnote{{The function $\Lambda$ is obtained by taking the projection of  $e^{\{H+K,\cdot\}}(\mu+\gamma-\mu^{(2,1,1)})$ to the subspace $C^{(2,1,1)}(\huaM)\oplus C^{(1,1,2)}(\huaM)\oplus C^{(0,1,3)}(\huaM)$ of $C^\infty(\huaM)$. See the footnote $4$ for more explanation of $\gamma_{H,K}$.}}
  $$\Lambda=\gamma+\gamma_{H,K}-\mu^{(2,1,1)},$$
where $\gamma_{H,K}$ is given by \eqref{eq:gammaHK}. Write $\Lambda=\Lambda^{(2,1,1)}+\Lambda^{(1,1,2)}+\Lambda^{(0,1,3)}$, where
  \begin{eqnarray*}
  \Lambda^{(2,1,1)}&=&\gamma^{(2,1,1)}=\mu^{(2,1,1)},\\
  \Lambda^{(1,1,2)}&=&\gamma^{(1,1,2)}+\Poisson{\mu^{(1,2,1)},H},\\
  \Lambda^{(0,1,3)}&=&\gamma^{(0,1,3)}+\Poisson{\mu^{(1,2,1)},K}+\half \Poisson{\Poisson{\mu^{(0,3,1)},H},H}.
\end{eqnarray*}
For all $\alpha^2,\beta^2,\gamma^2\in\Gamma(A^*_{-2}),~\alpha^1,\beta^1\in\Gamma(A^*_{-1})$ and $ f\in C^\infty(M)$, define
\begin{eqnarray}\label{eq:r-matrix-bi}
\left\{\begin{array}{rcl}
\tl^H_1(\alpha^1)&=&-\Poisson{\mu^{(2,1,1)},\alpha^1}=l_1^*(\alpha^1),\\
\tilde{\frka}_H(\alpha^2)(f)&=&-\Poisson{\Poisson{\gamma^{(1,1,2)}+\Poisson{\mu^{(1,2,1)},H},\alpha^2},f}=(\frka+\frka_H)(\alpha^2 )(f),\\
{\tl^H_2(\alpha^2,\beta^2)}&=&-\Poisson{\Poisson{\gamma^{(1,1,2)}+\Poisson{\mu^{(1,2,1)},H},\alpha^2},\beta^2}=(\frkl_2+\frkl_2^H)(\alpha^2,\beta^2),\\
{\tl^H_2(\alpha^2,\beta^1)}&=&-\Poisson{\Poisson{\gamma^{(1,1,2)}+\Poisson{\mu^{(1,2,1)},H},\alpha^2},\beta^1}=(\frkl_2+\frkl_2^H)(\alpha^2,\beta^1),\\
\tl^{H,K}_3(\alpha^2,\beta^2,\gamma^2)&=&-\Poisson{\Poisson{\Poisson{\gamma^{(0,1,3)}+\Poisson{\mu^{(1,2,1)},K},\alpha^2},\beta^2},\gamma^2}
-\half\Poisson{\Poisson{\Poisson{\Poisson{\Poisson{\mu^{(0,3,1)},H},H},\alpha^2},\beta^2},\gamma^2}\\
&=&(\frkl_3+\frkl_3^{H,K})(\alpha^2,\beta^2,\gamma^2),
\end{array}\right.
\end{eqnarray}
where $\frka_H, \frkl_2^H$ and $\frkl_3^{H,K}$ are given in Lemma \ref{lem:formulalH}.

The following result is a higher analogue of \cite[Section 6]{lwx}.
\begin{pro}\label{pro:bidual}
  With the above notations, if $H+K$ satisfies the following Maurer-Cartan type equation:
 \begin{eqnarray}\label{eq:GMC}
-\delta_{*}(H+K)+\half[H+K,H+K]_S+\frac{1}{6}[H+K,H+K,H+K]_S=0,
 \end{eqnarray}
where $\delta_*$ is the differential corresponding to the split Lie $2$-algebroid $\huaA^*[3]$, then we have
  $\Poisson{\Lambda,\Lambda}=0.$
Consequently, $\huaA^*[3]=(A_{-1}^*[3],A_{-2}^*[3],\tl_1^H,\tl_2^H,\tl_3^{H,K},\tilde{\frka}_H)$ is a split Lie $2$-algebroid.
\end{pro}
\pf First, note that \eqref{eq:GMC} is equivalent to the following equations
\begin{eqnarray}
  \bar{\delta}_{\ast}H&=&0,\label{eq:GMC-Cond 1}\\
 -\bar{\delta}_{\ast}K-\dM_{*}H+\half [H,H]_S&=&0,\label{eq:GMC-Cond 2}\\
 -\hat{\delta}_{\ast} H-\dM_{*}K+ {[H,K]_S}+\frac{1}{6}[H,H,H]_S&=&0.\label{eq:GMC-Cond 3}
\end{eqnarray}

On the other hand, $\Poisson{\Lambda,\Lambda}=0$ is equivalent to
$$\Poisson{\Lambda^{(1,1,2)},\Lambda^{(2,1,1)}}=0,\quad \Poisson{\Lambda^{(1,1,2)},\Lambda^{(1,1,2)}}+2\Poisson{\Lambda^{(2,1,1)},\Lambda^{(0,1,3)}}=0,\quad\Poisson{\Lambda^{(1,1,2)},\Lambda^{(0,1,3)}}=0.$$

By $\Poisson{\gamma^{(1,1,2)},\gamma^{(2,1,1)}}=0$, $\Poisson{\mu^{(1,2,1)},\gamma^{(2,1,1)}}=0$ and \eqref{eq:GMC-Cond 1}, we have
\begin{eqnarray*}
\Poisson{\Lambda^{(1,1,2)},\Lambda^{(2,1,1)}}&=&\Poisson{\gamma^{(1,1,2)}+\Poisson{\mu^{(1,2,1)},H},\gamma^{(2,1,1)}}=\Poisson{\Poisson{\mu^{(1,2,1)},H},\gamma^{(2,1,1)}}\\
&=&\Poisson{\mu^{(1,2,1)},\Poisson{H,\gamma^{(2,1,1)}}}=-\Poisson{\mu^{(1,2,1)},\bar{\delta}_{\ast}H}=0.
\end{eqnarray*}

By $\Poisson{\gamma^{(1,1,2)},\gamma^{(1,1,2)}}+2\Poisson{\gamma^{(2,1,1)},\gamma^{(0,1,3)}}=0$, $\Poisson{\mu^{(1,2,1)},\mu^{(1,2,1)}}+2\Poisson{\mu^{(2,1,1)},\mu^{(0,3,1)}}=0$,\\ $\Poisson{\gamma^{(1,1,2)},\mu^{(1,2,1)}}=0$, \eqref{eq:GMC-Cond 1} and \eqref{eq:GMC-Cond 2}, we have
\begin{eqnarray*}
&&\Poisson{\Lambda^{(1,1,2)},\Lambda^{(1,1,2)}}+2\Poisson{\Lambda^{(2,1,1)},\Lambda^{(0,1,3)}}\\
&=&\Poisson{\gamma^{(1,1,2)},\gamma^{(1,1,2)}}+2\Poisson{\gamma^{(2,1,1)},\gamma^{(0,1,3)}}+\half\Poisson{\gamma^{(1,1,2)},\Poisson{\mu^{(1,2,1)},H}}+\Poisson{\Poisson{\mu^{(1,2,1)},H},\Poisson{\mu^{(1,2,1)},H}}\\
&&+2\Poisson{\gamma^{(2,1,1)},\Poisson{\mu^{(1,2,1)},K}}+\Poisson{\gamma^{(2,1,1)},\Poisson{\Poisson{\mu^{(0,3,1)},H},H}}\\
&=&-2\Poisson{\mu^{(1,2,1)},\dM_* H}+\half\Poisson{\Poisson{\Poisson{\mu^{(1,2,1)},\mu^{(1,2,1)}},H},H}-\Poisson{\mu^{(1,2,1)},\Poisson{\Poisson{\mu^{(1,2,1)},H},H}}\\
&&-2\Poisson{\mu^{(1,2,1)},\bar{\delta}_* K}+\Poisson{\Poisson{\Poisson{\gamma_{2},\mu^{(0,3,1)}},H},H}\\
&=&2\Poisson{\mu^{(1,2,1)},\bar{\delta}_* K}+\Poisson{\mu^{(1,2,1)},\Poisson{\Poisson{\mu^{(1,2,1)},H},H}}+\Poisson{\Poisson{\half\Poisson{\mu^{(1,2,1)},\mu^{(1,2,1)}}+\Poisson{\mu^{(2,1,1)},\mu^{(0,3,1)}},H},H}\\
&&-\Poisson{\mu^{(1,2,1)},\Poisson{\Poisson{\mu^{(1,2,1)},H},H}}-2\Poisson{\mu^{(1,2,1)},\bar{\delta}_* K}\\
&=&0.
\end{eqnarray*}

By $\Poisson{\mu^{(1,2,1)},\mu^{(1,2,1)}}+2\Poisson{\mu^{(2,1,1)},\mu^{(0,3,1)}}=0$, $\Poisson{\gamma^{(1,1,2)},\mu^{(0,3,1)}}=0$ and \eqref{eq:GMC-Cond 2}, we have
\begin{eqnarray*}
  \Poisson{\Poisson{\mu^{(1,2,1)},\Poisson{\mu^{(1,2,1)},K}},H}&=&\Poisson{\Poisson{\half\Poisson{\mu^{(1,2,1)},\mu^{(1,2,1)}},K},H}\\
  &=&-\Poisson{\Poisson{\Poisson{\mu^{(2,1,1)},\mu^{(0,3,1)}},K},H}=-\Poisson{\Poisson{\mu^{(2,1,1)},K},\Poisson{\mu^{(0,3,1)},H}}\\
  &=&\half\Poisson{\Poisson{\Poisson{\mu^{(1,2,1)},H},H},\Poisson{\mu^{(0,3,1)},H}}+\Poisson{\Poisson{\gamma^{(1,1,2)},H},\Poisson{\mu^{(0,3,1)},H}}\\
  &=&\half\Poisson{\Poisson{\Poisson{\mu^{(1,2,1)},H},H},\Poisson{\mu^{(0,3,1)},H}}-\half\Poisson{\gamma^{(1,1,2)},\Poisson{\Poisson{\mu^{(0,3,1)},H},H}}.
\end{eqnarray*}
By the fact that $\Poisson{\mu^{(1,2,1)},\mu^{(0,3,1)}}=0$, we have
\begin{eqnarray*}
 &&\frac{1}{6}\Poisson{\mu^{(1,2,1)},\Poisson{\Poisson{\Poisson{\mu^{(0,3,1)},H},H},H}}\\
 &=&-\half\Poisson{\Poisson{\Poisson{\mu^{(1,2,1)},H},H},\Poisson{\mu^{(0,3,1)},H}}-\half \Poisson{\Poisson{\mu^{(1,2,1)},H},\Poisson{\Poisson{\mu^{(0,3,1)},H},H}}.
\end{eqnarray*}
By $\Poisson{\gamma^{(1,1,2)},\gamma^{(0,1,3)}}=0$, $\Poisson{\gamma^{(1,1,2)},\mu^{(1,2,1)}}=0$, $\Poisson{\mu^{(1,2,1)},\gamma^{(0,1,3)}}=0$ and \eqref{eq:GMC-Cond 3}, we have
\begin{eqnarray*}
  \Poisson{\Lambda^{(1,1,2)},\Lambda^{(0,1,3)}}&=&\Poisson{\gamma^{(1,1,2)},\Poisson{\mu^{(1,2,1)},K}}+\half\Poisson{\gamma^{(1,1,2)},\Poisson{\Poisson{\mu^{(0,3,1)},H},H}}+\Poisson{\Poisson{\mu^{(1,2,1)},H},\gamma^{(0,1,3)}}\\
  &&+\Poisson{\Poisson{\mu^{(1,2,1)},H},\Poisson{\mu^{(1,2,1)},K}}+\half \Poisson{\Poisson{\mu^{(1,2,1)},H},\Poisson{\Poisson{\mu^{(0,3,1)},H},H}}\\
  &=&-\Poisson{\mu^{(1,2,1)},\dM_*K}-\Poisson{\mu^{(1,2,1)},\hat{\delta}_*H}+\Poisson{\Poisson{\mu^{(1,2,1)},H},\Poisson{\mu^{(1,2,1)},K}}\\
  &&+\half\Poisson{\gamma^{(1,1,2)},\Poisson{\Poisson{\mu^{(0,3,1)},H},H}}+\half \Poisson{\Poisson{\mu^{(1,2,1)},H},\Poisson{\Poisson{\mu^{(0,3,1)},H},H}}\\
  &=&\Poisson{\mu^{(1,2,1)},\Poisson{\Poisson{\mu^{(1,2,1)},H},K}}+\frac{1}{6}\Poisson{\mu^{(1,2,1)},\Poisson{\Poisson{\Poisson{\mu^{(0,3,1)},H},H},H}}\\
  &&-\Poisson{\mu^{(1,2,1)},\Poisson{\Poisson{\mu^{(1,2,1)},H},K}}+\Poisson{\Poisson{\mu^{(1,2,1)},\Poisson{\mu^{(1,2,1)},K}},H}\\
  &&+\half\Poisson{\gamma^{(1,1,2)},\Poisson{\Poisson{\mu^{(0,3,1)},H},H}}+\half \Poisson{\Poisson{\mu^{(1,2,1)},H},\Poisson{\Poisson{\mu^{(0,3,1)},H},H}}\\
  &=&-\half\Poisson{\Poisson{\Poisson{\mu^{(1,2,1)},H},H},\Poisson{\mu^{(0,3,1)},H}}-\half \Poisson{\Poisson{\mu^{(1,2,1)},H},\Poisson{\Poisson{\mu^{(0,3,1)},H},H}}\\
  &&+\half\Poisson{\Poisson{\Poisson{\mu^{(1,2,1)},H},H},\Poisson{\mu^{(0,3,1)},H}}-\half\Poisson{\gamma^{(1,1,2)},\Poisson{\Poisson{\mu^{(0,3,1)},H},H}}\\
  &&+\half\Poisson{\gamma^{(1,1,2)},\Poisson{\Poisson{\mu^{(0,3,1)},H},H}}+\half \Poisson{\Poisson{\mu^{(1,2,1)},H},\Poisson{\Poisson{\mu^{(0,3,1)},H},H}}=0.
\end{eqnarray*}

Thus, we have $\Poisson{\Lambda,\Lambda}=0$. By Theorem \ref{thm:L2A-MST}, $\huaA^*[3]=(A_{-1}^*[3],A_{-2}^*[3],\tl_1^H,\tl_2^H,\tl_3^H,\tilde{\frka}_H)$ is a split Lie $2$-algebroid. \qed\vspace{3mm}

Define vector bundles $G_{-1}$ and $G_{-2}$ by
$$G_{-1}=\{\alpha^2 +H^\sharp(\alpha^2)\mid\alpha^2\in\Gamma(A^*_{-2})\},\quad G_{-2}=\{\alpha^1-H^\natural(\alpha^1)\mid\alpha^1\in\Gamma(A^*_{-1})\},$$
which are subbundles of $A_{-1}\oplus A_{-2}^*$ and $A_{-2}\oplus A_{-1}^*$ respectively. Denote by $\huaG=G_{-1}\oplus G_{-2}$. Define $l_1^\huaG:G_{-2}\longrightarrow G_{-1}$, $l^\huaG_2:\Gamma(G_i)\times \Gamma(G_j)\longrightarrow \Gamma(G_{i+j+1})$, $l^\huaG_3:\wedge^3G_{-1}\longrightarrow G_{-2}$ and $a^\huaG:G_{-1}\longrightarrow TM$ by
\begin{eqnarray}\label{eq:r-matrix-bi-g}
\left\{\begin{array}{rcl}
\lhg_1(\alpha^1-H^\natural(\alpha^1)) &=&l_1^*(\alpha^1)+H^\sharp(l_1^*(\alpha^1)),\\
\lhg_2(\alpha^2 +H^\sharp(\alpha^2),\beta^2 +H^\sharp(\beta^2))&=&(\frkl_2+\frkl_2^H)(\alpha^2,\beta^2)+H^\sharp\big((\frkl_2+\frkl_2^H)(\alpha^2,\beta^2)\big),\\
\lhg_2(\alpha^2 +H^\sharp(\alpha^2),\beta^1 -H^\natural(\beta^1))&=&(\frkl_2+\frkl_2^H)(\alpha^2,\beta^1)-H^\natural\big((\frkl_2+\frkl_2^H)(\alpha^2,\beta^1)\big),\\
\lhg_3 (\alpha^2 +H^\sharp(\alpha^2),\beta^2 +H^\sharp(\beta^2),\gamma^2+H^\sharp(\gamma^2))&=&  (\frkl_3+\frkl_3^{H,K})(\alpha^2 ,\beta^2 ,\gamma^2)-H^\natural\big((\frkl_3+\frkl_3^{H,K})(\alpha^2 ,\beta^2 ,\gamma^2)\big),\\
a^\huaG(\alpha^2 +H^\sharp(\alpha^2)) &=&(\frka+\frka_H)(\alpha^2),
\end{array}\right.
\end{eqnarray}
for all  $\alpha^2,\beta^2,\gamma^2\in\Gamma(A^*_{-2}),~\alpha^1,\beta^1\in\Gamma(A^*_{-1})$.

By Proposition \ref{pro:bidual}, we have
\begin{cor}
 Let $(\huaA,\huaA^*[3])$ be a split Lie $2$-bialgebroid, $H\in\Gamma(A_{-1}\odot A_{-2})$ and $K\in \Gamma(\wedge^3A_{-2})$. If $H+K$ satisfies the Maurer-Cartan type equation \eqref{eq:GMC}, then $\huaG_H=(G_{-2},G_{-1},\lhg_1,\lhg_2,\lhg_3,a^\huaG)$  is a split Lie $2$-algebroid.
\end{cor}

It is not hard to see that
$$\delta_{*}H\in\Gamma(\odot^2A_{-1})\oplus\Gamma(\wedge^2A_{-2}\odot A_{-1})\oplus \Gamma(\wedge^4A_{-2}),\qquad \delta_{*}K\in \Gamma(\wedge^2A_{-2}\odot A_{-1})\oplus \Gamma(\wedge^4A_{-2}).$$
We need the following preparation before we give the main result in this subsection.
\begin{lem}
  For all $\alpha^2,\beta^2,\gamma^2,\zeta^2\in\Gamma(A^*_{-2})$ and $\alpha^1,\beta^1,\gamma^1\in\Gamma(A^*_{-1})$,  we have
  \begin{eqnarray}
  \label{eq:GMC-formula 1}\bar{\delta}_{\ast}H(\alpha^1,\beta^1)&=& \langle-[H]_S(\alpha^1),\beta^1\rangle;\\
   \label{eq:GMC-formula 2}\bar{\delta}_{\ast}K(\alpha^2,\beta^2,\gamma^1)&=&\langle-[K]_S(\alpha^2,\beta^2),\gamma^1\rangle=\langle-[K]_S(\alpha^2,\gamma^1),\beta^2\rangle;\\
\label{eq:GMC-formula 3}\dM_{*}H(\alpha^2,\beta^2,\gamma^1)&=&\langle\huaL^1_{\alpha^2}H^\sharp(\beta^2)-\huaL^1_{\beta^2}H^\sharp(\alpha^2)-H^\sharp(\frkl_2(\alpha^2,\beta^2)),\gamma^1\rangle\\
    \label{eq:GMC-formula 33}&=&\langle\huaL^1_{\alpha^2}H^\natural(\gamma^1)-\iota_{\gamma^1}\delta_{*}H^\sharp(\alpha^2)-H^\natural(\frkl_2(\alpha^2,\gamma^1)),\beta^2\rangle;\\
\label{eq:GMC-formula 4}\hat{\delta}_{*}H(\alpha^2,\beta^2,\gamma^2,\zeta^2)&=& \langle-\huaL^3_{\alpha^2,\beta^2}H^\sharp(\gamma^2)-\huaL^3_{\beta^2,\gamma^2}H^\sharp(\alpha^2)
  -\huaL^3_{\gamma^2,\alpha^2}H^\sharp(\beta^2)\\
  \nonumber&&-H^\natural(l_3(\alpha^2,\beta^2,\gamma^2)),\zeta^2\rangle;\\
  \label{eq:GMC-formula 5}\dM_{*}K(\alpha^2,\beta^2,\gamma^2,\zeta^2)&=&\langle\huaL^1_{\alpha^2}K^\flat(\beta^2,\gamma^2)+\iota_{\gamma^2}\delta_{*}(K^\flat(\alpha^2,\beta^2))-\huaL^1_{\beta^2}K^\flat(\alpha^2,\gamma^2)\\
 &&-K^\flat(\frkl^H_2(\alpha^2,\beta^2),\gamma^2)
 \nonumber-K^\flat(\beta^2,\frkl^H_2(\alpha^2,\gamma^2))+K^\flat(\alpha^2,\frkl^H_2(\beta^2,\gamma^2)),\zeta^2\rangle.
  \end{eqnarray}
  \end{lem}

  \pf The proof  is similar to that of Lemma \ref{lem:MCformulars}. We omit the details. \qed\vspace{3mm}

  Now we are ready to give the main result in this paper, which says that the graph of a Maurer-Cartan element is a weak Dirac structure defined in Definition \ref{defi:wDirac}.

\begin{thm}
  Let $(\huaA,\huaA^*[3])$ be a split Lie $2$-bialgebroid, $H\in\Gamma(A_{-1}\odot A_{-2})$ and $K\in \Gamma(\wedge^3A_{-2})$. If $H+K$ satisfies the Maurer-Cartan type equation \eqref{eq:GMC}, then $(\ii_1,\ii_2,-\widetilde{K^\flat})$ is a morphism from the  Lie $2$-algebra $(\Gamma(G_{-2}),\Gamma(G_{-1}),\lhg_1,\lhg_2,\lhg_3)$  to the Leibniz $2$-algebra $(\Gamma(E_{-2}),\Gamma(E_{-1}),\partial,\diamond,\Omega)$ underlying the $\LWX$ $2$-algebroid   given in Proposition \ref{thm:Lie2biC2}, where $\ii_1$ and $\ii_2$ are inclusion maps from $G_{-1}$ and $G_{-2}$ to $E_{-1}$ and $E_{-2}$ respectively and $\widetilde{K^\flat}:\wedge^2G_{-1}\longrightarrow E_{-2}$ is defined by
  $$\widetilde{K^\flat}(\alpha^2 +H^\sharp(\alpha^2),\beta^2 +H^\sharp(\beta^2))=K^\flat(\alpha^2,\beta^2).$$
  Consequently, the split Lie $2$-algebroid $\huaG_H=(G_{-2},G_{-1},\lhg_1,\lhg_2,\lhg_3,a^\huaG)$ is a weak Dirac structure of the $\LWX$ $2$-algebroid   given in Proposition \ref{thm:Lie2biC2}.
\end{thm}
\pf First by the fact that $H$ is symmetric, i.e. $H(\alpha^2,\alpha^1)=H(\alpha^1,\alpha^2)$, it is obvious that the graded subbundle $\huaG_H$ is maximal isotropic.

Then by \eqref{eq:GMC-Cond 1}, we have
$-l_1\circ H^\natural=H^\sharp\circ l^*_1,$
which implies that
\begin{equation}\label{eq:M1}
  \ii_1\circ \lhg_1=\partial \circ \ii_2.
\end{equation}
By \eqref{eq:MC-formular1} and \eqref{eq:GMC-formula 3}, we have
\begin{eqnarray*}
  &&(-\bar{\delta}_{*}K-\dM_{*}H+\half[H,H]_S)(\alpha^2,\beta^2,-)=-\dM_{*}H(\alpha^2,\beta^2,-)+([K]_S+\half[H,H]_S(\alpha^2,\beta^2,-))\\
  &&\quad=H^\sharp(\frkl_2(\alpha^2,\beta^2))-\huaL^1_{\alpha^2}H^\sharp(\beta^2)+\huaL^1_{\beta^2}H^\sharp(\alpha^2)+H^\sharp\frkl^H_2(\alpha^2,\beta^2)-l_2(H^\sharp(\alpha^2),H^\sharp(\beta^2))+l_1K^\flat(\alpha^2,\beta^2).
\end{eqnarray*}
Thus, by \eqref{eq:GMC-Cond 2}, we deduce that
\begin{eqnarray}
\nonumber &&\lhg_2(\alpha^2 +H^\sharp(\alpha^2),\beta^2 +H^\sharp(\beta^2))-(\alpha^2 +H^\sharp(\alpha^2))\diamond (\beta^2 +H^\sharp(\beta^2))\\
 \nonumber&=&H^\sharp(\frkl_2(\alpha^2,\beta^2)+\frkl^H_2(\alpha^2,\beta^2))-\huaL^1_{\alpha^2}H^\sharp(\beta^2)+\huaL^1_{\beta^2}H^\sharp(\alpha^2)
   -l_2(H^\sharp(\alpha^2),H^\sharp(\beta^2))\\
  \label{eq:M2} &=&-l_1K^\flat(\alpha^2,\beta^2)=-l_1\widetilde{K^\flat}(\alpha^2 +H^\sharp(\alpha^2),\beta^2 +H^\sharp(\beta^2)).
\end{eqnarray}
Similarly, by \eqref{eq:MC-formular2} and \eqref{eq:GMC-formula 33}, we have
\begin{eqnarray*}
  (-\bar{\delta}_{*}K-\dM_{*}H+\half[H,H]_S)(\alpha^2,\beta^1,-)&=&-\dM_{*}H(\alpha^2,\beta^1,-)+([K]_S+\half[H,H]_S(\alpha^2,\beta^1,-))\\
  &=&-\huaL^1_{\alpha^2}H^\natural(\beta^1)+\iota_{\beta^1}\dM_{*}H^\sharp(\alpha^2)+H^\natural(\frkl_2(\alpha^2,\beta^1))\\
  &&+H^\natural(\frkl^H_2(\alpha^2,\beta^1))-l_2(H^\sharp(\alpha^2),H^\natural(\beta^1))-K^\flat(\alpha^2,l^*_1(\beta^1)).
\end{eqnarray*}
Thus, by \eqref{eq:GMC-Cond 2}, we deduce that
\begin{equation}\label{eq:M3}
  \lhg_2(\alpha^2 +H^\sharp(\alpha^2),\beta^1 -H^\natural(\beta^1))-(\alpha^2 +H^\sharp(\alpha^2))\diamond (\beta^1 -H^\natural(\beta^1))=-\widetilde{K^\flat}(\alpha^2 +H^\sharp(\alpha^2),\lhg_1(\beta^1 -H^\natural(\beta^1))).
\end{equation}
By \eqref{eq:MC-formular3}, \eqref{eq:MC-formular4}, \eqref{eq:GMC-formula 4} and \eqref{eq:GMC-formula 5}, we have
\begin{eqnarray*}
   &&(-\hat{\delta}_{\ast} H-\dM_{*}K+ {[H,K]_S}+\frac{1}{6}[H,H,H]_S)(\alpha^2,\beta^2,\gamma^2,-)\\
   &=&\huaL^3_{\alpha^2,\beta^2}H^\sharp(\gamma^2)+\huaL^3_{\beta^2,\gamma^2}H^\sharp(\alpha^2)+\huaL^3_{\gamma^2,\alpha^2}H^\sharp(\beta^2)+H^\natural(\frkl_3(\alpha^2,\beta^2,\gamma^2)) \\ && -\huaL^1_{\alpha^2}K^\flat(\beta^2,\gamma^2) -\iota_{\gamma^2}\dM_{*}(K^\flat(\alpha^2,\beta^2))+\huaL^1_{\beta^2}K^\flat(\alpha^2,\gamma^2)\\
  &&+K^\flat(\frkl_2(\alpha^2,\beta^2),\gamma^2)+K^\flat(\beta^2,\frkl_2(\alpha^2,\gamma^2))-K^\flat(\alpha^2,\frkl_2(\beta^2,\gamma^2))\\
 &&+K^\flat(\frkl^H_2(\alpha^2,\beta^2),\gamma^2)+K^\flat(\frkl^H_2(\gamma^2,\alpha^2),\beta^2)+K^\flat(\frkl^H_2(\beta^2,\gamma^2),\alpha^2)\nonumber\\
 &&-l_2(H^\sharp(\alpha^2),K^\flat(\beta^2,\gamma^2))-l_2(H^\sharp(\gamma^2),K^\flat(\alpha^2,\beta^2))-l_2(H^\sharp(\beta^2),K^\flat(\gamma^2,\alpha^2))\nonumber\\
 &&+H^\natural(\frkl^{H,K}_3(\alpha^2,\beta^2,\gamma^2))+l_3(H^\sharp(\alpha^2),H^\sharp(\beta^2),H^\sharp(\gamma^2)).
\end{eqnarray*}
Thus, by \eqref{eq:GMC-Cond 3}, we deduce that
\begin{eqnarray}
 \nonumber &&-\lhg_3(\alpha^2 +H^\sharp(\alpha^2),\beta^2 +H^\sharp(\beta^2),\gamma^2+H^\sharp(\gamma^2))-(\alpha^2 +H^\sharp(\alpha^2))\diamond \widetilde{K^\flat}(\beta^2 +H^\sharp(\beta^2),\gamma^2+H^\sharp(\gamma^2))\\
  \nonumber &&+(\beta^2 +H^\sharp(\beta^2))\diamond \widetilde{K^\flat}(\alpha^2 +H^\sharp(\alpha^2),\gamma^2+H^\sharp(\gamma^2))-\widetilde{K^\flat}(\alpha^2 +H^\sharp(\alpha^2),\beta^2+H^\sharp(\beta^2))\diamond (\gamma^2+H^\sharp(\gamma^2))\\
 \nonumber &&+\widetilde{K^\flat}(\lhg_2(\alpha^2,\beta^2)+H^\sharp\lhg_2(\alpha^2,\beta^2),\gamma^2+H^\sharp(\gamma^2))
  -\widetilde{K^\flat}(\alpha^2+H^\sharp(\alpha^2),\lhg_2(\beta^2,\gamma^2)+H^\sharp\lhg_2(\beta^2,\gamma^2))\\
\nonumber&& +\widetilde{K^\flat}(\beta^2+H^\sharp(\beta^2),\lhg_2(\alpha^2,\gamma^2)+H^\sharp\lhg_2(\alpha^2,\gamma^2))+\Omega(\alpha^2 +H^\sharp(\alpha^2),\beta^2 +H^\sharp(\beta^2),\gamma^2+H^\sharp(\gamma^2))\\
\label{eq:M4} &=&0.
\end{eqnarray}
By \eqref{eq:M1}-\eqref{eq:M4}, we deduce that $(\ii_1,\ii_2,-\widetilde{K^\flat})$ is a morphism from the Lie 2-algebra $(\Gamma(G_{-2}),\Gamma(G_{-1}),\lhg_1,\lhg_2,\lhg_3)$  to the Leibniz $2$-algebra $(\Gamma(E_{-2}),\Gamma(E_{-1}),\partial,\diamond,\Omega)$.

Finally, it is obvious that
$$
\rho(\alpha^2 +H^\sharp(\alpha^2))=a(H^\sharp(\alpha^2))+\frka(\alpha^2)=a^\huaG(\alpha^2 +H^\sharp(\alpha^2)).
$$
Therefore,  the split Lie $2$-algebroid $(G_{-2},G_{-1},\lhg_1,\lhg_2,\lhg_3,a^\huaG)$ is a weak Dirac structure of the $\LWX$ $2$-algebroid $(E_{-2},E_{-1},\partial,\rho,(\cdot,\cdot)_+,\diamond,\Omega)$ given in Proposition \ref{thm:Lie2biC2}. \qed

\subsection{Examples}

In this subsection, we give some examples of Theorem \ref{thm:MCE1} including the string Lie 2-algebra and  split Lie 2-algebroids constructed from   integrable distributions and left-symmetric algebroids.
\begin{ex}{\rm
  Let $(\g,[\cdot,\cdot]_\g)$ be a   semisimple Lie algebra and $B(\cdot,\cdot)$  its killing form. Recall that the string Lie $2$-algebra $(\Real[2], \g[1],l_1,l_2,l_3)$ is given by
\begin{eqnarray*}
 l_1=0,\quad l_2(x,y)=[x,y]_\g,\quad l_2(x,v)=0,\quad
  l_3(x,y,z)=B(x,[y,z]_\g),\quad\forall~x,y,z\in\g,~v\in\Real.
\end{eqnarray*}
We have $\g\odot\R=\g$ and $\wedge^3\R=0$. It is straightforward to see that any $h\in\g$ is a Maurer-Cartan element of the associated homotopy Poisson algebra. Furthermore, $h^\sharp:\Real^*\cong \Real\longrightarrow \g$  and $h^\natural:\g^*\longrightarrow\R$ are given by
$$
h^\sharp(s)=sh,\quad h^\natural(\alpha)=\langle h,\alpha\rangle,\quad\forall s\in\R, ~\alpha\in\g^*.
$$
The Lie 2-algebra $(\g^*[2],\R[1],\frkl_1^h,\frkl_2^h,\frkl_3^{h})$ given in Theorem \ref{thm:MCE1} (i) is given by
\begin{eqnarray*}
 \frkl_1^h=0,\quad  \frkl_2^h(s,t)=0,\quad \frkl_2^h(s,\beta)=\ad^*_{s h}\beta,\quad
  \frkl_3^{h}(s,t,w)=0,\quad\forall~\beta\in\g^*,~s,t,w\in\Real.
\end{eqnarray*}

Moreover, the Lie $2$-algebras $(\g^*,\Real,\frkl_1^h,\frkl_2^h,\frkl_3^h)$ and $(\Real,\g,l_1,l_2,l_3)$ define a Lie 2-bialgebra, whose double is the Lie 2-algebra $(\Real\oplus \g^*[2],\g\oplus \Real[1],\partial, \diamond,\Omega)$ given by
\begin{eqnarray*}
\left\{\begin{array}{rcl}
\partial&=&0,\\
(x,s)\diamond(y,t)&=&([x,y]_\g,0),\\
(x,s)\diamond(u,\alpha)&=&(u,\alpha)\diamond(x,s)=(0,\ad^*_x\alpha+\ad^*_{sh}\alpha),\\
\Omega((x,s),(y,t),(z,r))&=&(B(x,[y,z]_\g),-r B^\sharp([x,y]_\g)-s B^\sharp([y,z]_\g)-t B^\sharp([z,x]_\g)),
\end{array}\right.
\end{eqnarray*}
for all $(x,s),(y,t),(z,r)\in\g\oplus\R$ and $(u,\alpha)\in\R\oplus \g^*$ and $B^\sharp:\g\longrightarrow \g^*$ is given by
$\langle B^\sharp(x),y\rangle=B(x,y).$
}
\end{ex}

Let $(A,[\cdot,\cdot]_A,a_A)$ be a Lie algebroid and $\nabla:\Gamma(A)\times\Gamma(E)\longrightarrow \Gamma(E)$ a representation of Lie   $A$ on a vector bundle $E$. Then it is straightforward to see that $\huaA=(E[2],A[1],l_1,l_2,l_3,a)$ is a split Lie $2$-algebroid, where  $ l_1=0,~ l_3=0, ~ a=a_A$ and
  \begin{eqnarray*}
  l_2(X,Y)=[X,Y]_A, \quad l_2(X,e)=\nabla_Xe, \quad\forall~X,Y\in\Gamma(A),~e\in\Gamma(E).
  \end{eqnarray*}

  \begin{pro}\label{ex:example semi-prod 1}
    With above notations, let $H\in\Gamma(A\odot E)$ and $K\in\Gamma(\wedge^3 E)$.
    \begin{itemize}
      \item[\rm(i)] If $[H,H]_S=0$, then $(E^*,\frkl_2^H,\frka_H)$ is a Lie algebroid, where $\frka_H=a_A\circ H^\sharp$ and $ \frkl_2^H$ is given by
     \begin{equation}\label{eq:l2H}
      \frkl_2^H(e^*_1,e^*_2)=\nabla^*_{H^\sharp(e^*_1)}e^*_2-\nabla^*_{H^\sharp(e^*_2)}e^*_1,\quad \forall e^*_1,e^*_2\in\Gamma(E^*),
      \end{equation}
      Here $\nabla^*$ is the dual representation of $\nabla$ on $E^*$.
      \item[\rm(ii)]  If $[H,H]_S=0$, then $[H,K]_S=0$ if and only if $K$ is a $3$-cocycle on the Lie algebroid $(E^*,\frkl_2^H,\frka_H)$.
    \end{itemize}
  \end{pro}
  \pf (i) By \eqref{eq:r-matrix2}, it is straightforward to deduce that  $\frkl_2^H$ is given by \eqref{eq:l2H}. Since $l_1=0,$ we obtain $\frkl_1^H=l_1^*=0$. Thus, $\frkl_2^H$ satisfies the Jacobi identity and $(E^*,\frkl_2^H,\frka_H)$ is a Lie algebroid.

  (ii) For all $e^*_1,e^*_2,e^*_3,e^*_4\in\Gamma(E^*)$, we have
  \begin{eqnarray*}
 &&\langle[H,K]_S(e^*_1,e^*_2,e^*_3),e^*_4\rangle\\
 &=&\langle K^\flat(\frkl^H_2(e^*_1,e^*_2),e^*_3)+K^\flat(\frkl^H_2(e^*_3,e^*_1),e^*_2)+K^\flat(\frkl^H_2(e^*_2,e^*_3),e^*_1)-l_2(H^\sharp(e^*_1),K^\flat(e^*_2,e^*_3))\\
 &&-l_2(H^\sharp(e^*_3),K^\flat(e^*_1,e^*_2))-l_2(H^\sharp(e^*_2),K^\flat(e^*_3,e^*_1))-H^\natural(L^2_{K^\flat(e^*_1,e^*_2)}e^*_3+L^2_{K^\flat(e^*_3,e^*_1)}e^*_2\\
 &&+L^2_{K^\flat(e^*_2,e^*_3)}e^*_1+2\dM K(e^*_1,e^*_2,e^*_3)),e^*_4\rangle\\
 &=&K(\frkl_2^H(e^*_1,e^*_2),e^*_3,e^*_4)+K(\frkl_2^H(e^*_3,e^*_1),e^*_2,e^*_4)+K(\frkl_2^H(e^*_2,e^*_3),e^*_1,e^*_4)-\langle \nabla^*_{H^\sharp(e^*_1)} K^\flat(e^*_2,e^*_3),e^*_4\rangle\\
 &&-\langle \nabla^*_{H^\sharp(e^*_3)} K^\flat(e^*_1,e^*_2),e^*_4\rangle-\langle \nabla^*_{H^\sharp(e^*_2)} K^\flat(e^*_3,e^*_1),e^*_4\rangle+\langle \nabla^*_{H^\sharp(e^*_4)} K^\flat(e^*_1,e^*_2),e^*_3\rangle\\
 &&+\langle \nabla^*_{H^\sharp(e^*_4)} K^\flat(e^*_3,e^*_1),e^*_2\rangle+\langle \nabla^*_{H^\sharp(e^*_4)} K^\flat(e^*_2,e^*_3),e^*_1\rangle-2\frka_H(e^*_4)K(e^*_1,e^*_2,e^*_3)
 \end{eqnarray*}
 \begin{eqnarray*}
 &=&K(\frkl_2^H(e^*_1,e^*_2),e^*_3,e^*_4)+K(\frkl_2^H(e^*_3,e^*_1),e^*_2,e^*_4)+K(\frkl_2^H(e^*_2,e^*_3),e^*_1,e^*_4)-\frka_H(e^*_1)K(e^*_2,e^*_3,e^*_4)\\
 &&+K(e^*_2,e^*_3,\nabla^*_{H^\sharp(e^*_1)}e^*_4)-\frka_H(e^*_3)K(e^*_1,e^*_2,e^*_4)+K(e^*_1,e^*_2,\nabla^*_{H^\sharp(e^*_3)}e^*_4)-\frka_H(e^*_2)K(e^*_3,e^*_1,e^*_4)\\
 &&+K(e^*_3,e^*_1,\nabla^*_{H^\sharp(e^*_2)}e^*_4)+\frka_H(e^*_4)K(e^*_1,e^*_2,e^*_3)-K(e^*_1,e^*_2,\nabla^*_{H^\sharp(e^*_4)}e^*_3)+\frka_H(e^*_4)K(e^*_2,e^*_3,e^*_1)\\
 &&-K(e^*_3,e^*_1,\nabla^*_{H^\sharp(e^*_4)}e^*_2)+\frka_H(e^*_4)K(e^*_2,e^*_3,e^*_1)-K(e^*_2,e^*_3,\nabla^*_{H^\sharp(e^*_4)}e^*_1)-2\frka_H(e^*_4)K(e^*_1,e^*_2,e^*_3)\\
 &=&K(\frkl_2^H(e^*_1,e^*_2),e^*_3,e^*_4)+K(\frkl_2^H(e^*_3,e^*_1),e^*_2,e^*_4)+K(\frkl_2^H(e^*_2,e^*_3),e^*_1,e^*_4)+K(e^*_2,e^*_3,\frkl_2^H(e^*_1,e^*_4))\\
 &&+K(e^*_1,e^*_2,\frkl_2^H(e^*_3,e^*_4))+K(e^*_3,e^*_1,\frkl_2^H(e^*_2,e^*_4))-\frka_H(e^*_1)K(e^*_2,e^*_3,e^*_4)-\frka_H(e^*_3)K(e^*_1,e^*_2,e^*_4)\\
 &&-\frka_H(e^*_2)K(e^*_3,e^*_1,e^*_4)+\frka_H(e^*_4)K(e^*_1,e^*_2,e^*_3)\\
 &=&-(\dM^HK)(e^*_1,e^*_2,e^*_3,e^*_4),
  \end{eqnarray*}
  where $d^H$ is the coboundary operator on the Lie algebroid $(E^*,\frkl_2^H,\frka_H)$ with the coefficient in the trivial representation. Thus, $ [H,K]_S=0$ if and only if $K$ is a 3-cocycle. \qed

\begin{cor}\label{cor:ex}

Let $H\in\Gamma(A\odot E)$ and $K\in\Gamma(\wedge^3 E)$ such that
$[H,H]_S=0,\quad[H,K]_S=0.$
 \emptycomment{ By a direct calculation in local coordinates, we find that {\bf $\bar{H}$ is a Poisson structure on the Lie algebroid $A\ltimes E$ if and only if $[H,H]_S=0$ for the split Lie $2$-algebroid $\huaA[-3]$.}}
Then $\huaA^*[3]=( A^*,E^*,\frkl^H_1,\frkl_2^H,\frkl_3^{H,K},\frka_H)$ is a split Lie $2$-algebroid, where $\frkl^H_1=0,$ and
\begin{eqnarray*}
{\frkl_2^H(e^*_1,e^*_2)}&=&\nabla^*_{H^\sharp(e^*_1)}e^*_2-\nabla^*_{H^\sharp(e^*_2)}e^*_1,\\
{\frkl_2^H(e^*_1,\beta)}&=&\frkL_{H^\sharp(e^*_1)}\beta+\langle\nabla_\cdot H^\natural(\beta),e^*_1\rangle-\dM^A H(e^*_1,\beta), \\
\frkl_3^{H,K}(e^*_1,e^*_2,e^*_3)&=&\langle\nabla_\cdot K^\flat(e^*_1,e^*_2),e^*_3\rangle+\langle\nabla_\cdot K^\flat(e^*_3,e^*_1),e^*_2\rangle+\langle\nabla_\cdot K^\flat(e^*_2,e^*_3),e^*_1\rangle-2\dM^A K(e^*_1,e^*_2,e^*_3),
\end{eqnarray*}
 for all $e^*_1,e^*_2,e^*_3\in\Gamma(E^*),~\beta\in\Gamma(A^*)$. Here $\frkL_{X}:\Gamma(A^*)\longrightarrow\Gamma(A^*)$ and $\dM^A$  are the Lie derivative  and  the differential for the Lie algebroid $A$, respectively. Furthermore, $(H^\sharp,-H^\natural,-K^\flat)$ is a morphism form the split Lie $2$-algebroid $\huaA^*[3]$ to the split Lie $2$-algebroid $\huaA$ and $(\huaA,\huaA^*[3])$ is a split Lie $2$-bialgebroid.

\end{cor}

According to Proposition \ref{ex:example semi-prod 1} and Corollary \ref{cor:ex}, we can give the following example in which the Lie algebroid is given by an integral distribution and the representation is given by the Lie derivative on its normal bundle.

\emptycomment{
\begin{ex}
  Let $\huaF$ be an integral distribution on a manifold $M$ and $N:=TM/ \huaF$ be the normal bundle of $\huaF$. It is well-known that $\nabla:\Gamma(\huaF)\times \Gamma(N)\longrightarrow \Gamma(N) $ given by
  $$\nabla_X(Y+\huaF)=[X,Y]_{TM}+\huaF,\quad\forall~X\in\Gamma(\huaF),Y\in\frkX(M)$$
  is a representation of the Lie algebroid $F$ on $N$. Then $\huaA=(A_{-2}=\huaF,A_{-1}=N ,l_1,l_2,l_3,a)$ is a split Lie $2$-algebroid, where  $ l_1=0, l_3=0,  a=\ii$ and
  \begin{eqnarray*}
  l_2(X_1,X_2)=[X_1,X_2]_{TM}, \quad l_2(X,Y+\huaF)=\nabla_X(Y+\huaF), \quad\forall~X,X_1,X_2\in\Gamma(\huaF),~Y\in\frkX(M).
  \end{eqnarray*}
\end{ex}
}

\begin{ex}{\rm
  Let $\huaF\subset TM$ be an integral distribution on a manifold $M$ and $\huaF^\bot\subset T^*M$ its conormal bundle. Then $(\huaF^\bot[2],\huaF[1],l_1,l_2,l_3,a)$ is a Lie 2-algebroid, where $l_1=0,~l_3=0,$ $a:\huaF\longrightarrow TM$ is the inclusion map and $l_2$ is given by
  \begin{eqnarray*}
    l_2(X_1,X_2)=[X_1,X_2], \quad l_2(X,\xi)=\frkL_X\xi, \quad \forall X_1,X_2, X\in\Gamma(\huaF),~\xi\in\Gamma(\huaF^\bot),
  \end{eqnarray*}
 where $\frkL$ is the usual Lie derivative on $M$.
 \emptycomment{\liu{Here $\frkL$ means that $$\langle\frkL_X\xi,\overline{Y}\rangle=X\langle\xi,Y\rangle-\langle\xi,\overline{[X,Y]}\rangle.$$
 The fact that $\frkL$ is the usual Lie derivative on the manifold $M$ is not true.
 } }Then $H\in\huaF\odot\huaF^\bot$ satisfies $[H,H]_S=0$ if and only if
  $$
  H^\sharp(L^1_{H^\sharp(\alpha_1)}\alpha_2-L^1_{H^\sharp(\alpha_2)}\alpha_1)=[H^\sharp(\alpha_1),H^\sharp(\alpha_2)],\quad\forall~\alpha_1,\alpha_2\in\Gamma((\huaF^\bot)^*),
  $$
where $L^1:\huaF\times (\huaF^\bot)^*\longrightarrow (\huaF^\bot)^*$ is the Lie derivative on the Lie 2-algebroid defined by \eqref{eq:L1}. In fact, there is a natural isomorphism between $TM/\huaF$ and $(\huaF^\bot)^*$. For any $Y\in\Gamma(TM)$, we denote by $\overline{Y}$ its image in $\Gamma(TM/\huaF)$ of the natural projection $\pr:TM\longrightarrow TM/\huaF$. Then we have
$$
L^1_{X}\overline{Y}=\overline{[X,Y]}, \quad \forall X\in\Gamma(\huaF), Y\in\Gamma(TM).
$$
Such an $H\in\huaF\odot\huaF^\bot$ induces a Lie algebroid   $((\huaF^\bot)^*,[\cdot,\cdot]_H,\frka_H)$, where
  $$[\alpha_1,\alpha_2]_H=L^1_{H^\sharp(\alpha_1)}\alpha_2-L^1_{H^\sharp(\alpha_2)}\alpha_1,\quad \frka_H=H^\sharp,\quad\forall~\alpha_1,\alpha_2\in\Gamma((\huaF^\bot)^*).$$

  Let $K\in\Gamma(\wedge^3\huaF^\bot)$ be a 3-cocycle on the Lie algebroid $((\huaF^\bot)^*,[\cdot,\cdot]_H,\frka_H)$. Then there is an induced Lie 2-algebroid   $(\huaF^*[2],(\huaF^\bot)^*[1],\frkl^H_1=0,[\cdot,\cdot]_H,[\cdot,\cdot,\cdot]_{H,K},\frka_H=H^\sharp)$, where
  \begin{eqnarray*}
  [\alpha_1,\alpha_2]_H&=&L^1_{H^\sharp(\alpha_1)}\alpha_2-L^1_{H^\sharp(\alpha_2)}\alpha_1,\\
{[\alpha,\theta]_H}&=&\frkL^\huaF_{H^\sharp(\alpha)}\theta-L^2_{H^\natural(\theta)}\alpha-d^\huaF H(\alpha,\theta),\\
{[\alpha_1,\alpha_2,\alpha_3]_{H,K}}&=&-2d^\huaF( K(\alpha_1,\alpha_2,\alpha_3))-L^2_{K^\flat(\alpha_1,\alpha_2)}\alpha_3-L^2_{K^\flat(\alpha_3,\alpha_1)}\alpha_2-L^2_{K^\flat(\alpha_2,\alpha_3)}\alpha_1,
  \end{eqnarray*}
for all $\alpha,\alpha_1,\alpha_2,\alpha_3\in\Gamma((\huaF^\bot)^*),~\theta\in\Gamma(\huaF^*)$. Here  $\frkL^\huaF:\huaF\times \huaF^*\longrightarrow \huaF^*$ and $d^\huaF:\wedge^k\huaF^*\longrightarrow \wedge^{k+1}\huaF^*$ are the Lie derivative and the differential for the Lie algebroid $\huaF$ respectively. Furthermore, it is straightforward to deduce that the relation between $L^2$ and the Lie derivative $\frkL$ is given by
$$
\langle L^2_{\xi}\alpha,X\rangle=-\langle\frkL_X\xi,\alpha\rangle,\quad \forall X\in\Gamma(\huaF), \xi\in\Gamma(\huaF^\bot), \alpha\in\Gamma((\huaF^\bot)^*).
$$

}\end{ex}

The notion of a left-symmetric algebroid, also called Koszul-Vinberg algebroid, was introduced in \cite{LiuShengBaiChen,Boyom1,Boyom2} as a geometric generalization of a left-symmetric algebra (pre-Lie algebra).
  Let $(A,\cdot_A, a_A)$ be a left-symmetric algebroid. Define  a skew-symmetric bilinear bracket operation $[\cdot,\cdot]_A$ on $\Gamma(A)$ by
  $$
  [x,y]_A=x\cdot_A y-y\cdot_A x,\quad \forall ~x,y\in\Gamma(A).
  $$
Then, $(A,[\cdot,\cdot]_A,a_A)$ is a Lie algebroid, and denoted by
$A^c$, called the {\bf sub-adjacent Lie algebroid} of
 $(A,\cdot_A,a_A)$. Furthermore, $L:A\longrightarrow \frkD(A)$ defined by $L_XY=X\cdot_A Y$  gives a
  representation of the Lie algebroid  $A^c$ on $A$, where $\frkD(A)$ denotes the first order covariant differential operator bundle of the vector bundle $A$.

\begin{ex}{\rm
  Let $(A,\cdot_A,a_A)$ be a left-symmetric algebroid. Then $\huaA=(A^*[2], A[1],l_1,l_2,l_3,a_A)$ is a split Lie $2$-algebroid, where $l_1=0,~ l_3=0$ and $l_2$ is given by
  \begin{eqnarray*}
    l_2(X,Y)=[X,Y]_A,\quad  l_2(X,\eta)=L^*_X\eta,\quad\forall~X,Y\in\Gamma(A),~\eta\in\Gamma(A^*).
  \end{eqnarray*}
Here $L^*$ is the dual representation of $L$.  Let $H\in A\odot A^*$ be given by
  $$
  H(X,\xi)=\langle X,\xi\rangle,\quad \forall X\in\Gamma(A),~\xi\in\Gamma(A^*).
  $$
  That is, $H^\sharp={\id}_{A}$ and $H^\natural ={\id}_{A^*}$. Then we have $[H,H]_S=0$. In fact, it follows from
  \begin{eqnarray*}
    \half[H,H]_S(X,Y,\cdot)=H^\sharp(X\cdot_AY-Y\cdot_AX)-l_2(H^\sharp(X),H^\sharp(Y))=[X,Y]_A-[X,Y]_A=0.
  \end{eqnarray*}
Furthermore, the induced Lie algebroid structure on $A=(A^*)^*$ is exactly the sub-adjacent Lie algebroid $A^c$. For a $K\in\Gamma(\wedge^3 A^*)$,  by Proposition \ref{ex:example semi-prod 1},
\emptycomment{we have
 \begin{eqnarray*}
 &&\langle[H,K]_S(X,Y,Z),W\rangle\\
 &=&\langle K^\flat(\frkl^H_2(X,Y),Z)+K^\flat(\frkl^H_2(Z,X),Y)+K^\flat(\frkl^H_2(Y,Z),X)-l_2(H^\sharp(X),K^\flat(Y,Z))\\
 &&-l_2(H^\sharp(Z),K^\flat(X,Y))-l_2(H^\sharp(Y),K^\flat(Z,X))-H^\natural(L^2_{K^\flat(X,Y)}Z+L^2_{K^\flat(Z,X)}Y\\
 &&+L^2_{K^\flat(Y,Z)}X+2\dM K(X,Y,Z)),W\rangle\\
 &=&K([X,Y]_A,Z,W)+K([Z,X]_A,Y,W)+K([Y,Z]_A,X,W)-\langle L^*_X K^\flat(Y,Z),W\rangle\\
 &&-\langle L^*_Z K^\flat(X,Y),W\rangle-\langle L^*_Y K^\flat(Z,X),W\rangle+\langle L^*_W K^\flat(X,Y),Z\rangle+\langle L^*_W K^\flat(Z,X),Y\rangle\\
 &&+\langle L^*_W K^\flat(Y,Z),X\rangle-2a_A(W)K(X,Y,Z)\\
 &=&K([X,Y]_A,Z,W)+K([Z,X]_A,Y,W)+K([Y,Z]_A,X,W)+K(Y,Z,X\cdot_A W)\\
 &&-a_A(X)K(Y,Z,W)+K(X,Y,Z\cdot_AW)-a_A(Z)K(X,Y,W)+K(Z,X,Y\cdot_AW)\\
 &&-a_A(Y)K(Z,X,W)+a_A(W)K(X,Y,Z)-K(X,Y,W\cdot_AZ)+a_A(W)K(Z,X,Y)\\
 &&-K(Z,X,W\cdot_AY)+a_A(W)K(Y,Z,X)-K(Y,Z,W\cdot_AX)-2a_A(W)K(X,Y,Z)\\
 &=&K([X,Y]_A,Z,W)+K([Z,X]_A,Y,W)+K([Y,Z]_A,X,W)+K(Y,Z,[X,W]_A)\\
 &&-a_A(X)K(Y,Z,W)+K(X,Y,[Z,W]_A)-a_A(Z)K(X,Y,W)+K(Z,X,[Y,W]_A)\\
 &&-a_A(Y)K(Z,X,W)+a_A(W)K(X,Y,Z)\\
 &=&-(\dM^AK)(X,Y,Z,W).
  \end{eqnarray*}}
$[H,K]_S=0$ if and only if $K$ is a $3$-cocycle on the sub-adjacent Lie algebroid $A^c$.
Under these conditions, $\huaA^*[3]=(A^*[2], A[1],\frkl^H_1,\frkl^H_2,\frkl^{H,K}_3,\frka_H)$ is   a split Lie $2$-algebroid, where $\frkl^H_1=0$, and
\begin{eqnarray*}
{\frkl_2^H(X,Y)}&=&[X,Y]_A,\qquad
{\frkl_2^H(X,\xi)}=L^*_{X}\xi, \\
\frkl_3^{H,K}(X,Y,Z)&=&R^*_ZK^b(X,Y)+R^*_YK^b(Z,X)+R^*_XK^b(Y,Z)+\dM^A( K(X,Y,Z)).
\end{eqnarray*}
for all $X,Y,Z\in\Gamma(A)$ and $\xi\in\Gamma(A^*)$. Here $R_X^*:A^*\longrightarrow A^*$ is the dual map of the right multiplication $R_X$, i.e. $\langle R^*_X\xi,Y\rangle=-\langle\xi,Y\cdot_A X\rangle$.}
\end{ex}

Since left-symmetric algebras are left-symmetric algebroids naturally, the following example is a special case of the above example.

\begin{ex}{\rm
 Let $(\g,\cdot)$ be a 3-dimensional left-symmetric algebra generated by the following relations
 $$e_1\cdot e_1=2e_1,\quad e_1\cdot e_2=e_2,\quad e_1\cdot e_3=e_3,\quad e_2\cdot e_3=e_3\cdot e_2=e_1,$$
 where $\{e_1,e_2,e_3\}$ is a basis of $\g$. The corresponding sub-adjacent Lie algebra structure is given by
 $$[e_1,e_2]_\g=e_2,\quad [e_1,e_3]_\g=e_3.$$
 The dual representation $L^*$ of the sub-adjacent Lie algebra $\g^c$ on $\g^*$ is given by
 $$L^*_{e_1}e^*_1=-2e^*_1,\quad L^*_{e_1}e^*_2=-e^*_2,\quad L^*_{e_1}e^*_3=-e^*_3,\quad L^*_{e_2}e^*_1=-e^*_3,\quad L^*_{e_3}e^*_1=-e^*_2,$$
 where $\{e^*_1,e^*_2,e^*_3\}$ is the dual basis.

 The dual map of the right multiplication $R$ is given by
$$R^*_{e_1}e^*_1=-2e^*_1,\quad R^*_{e_2}e^*_1=-e^*_3,\quad R^*_{e_2}e^*_2=-e^*_1,\quad R^*_{e_3}e^*_1=-e^*_2,\quad R^*_{e_3}e^*_3=-e^*_1.$$

Let $H\in \g\odot \g^*$ be given by
  $$
  H(x,\xi)=\langle x,\xi\rangle,\quad \forall x\in\g,~\xi\in\g^*.
  $$
  That is, $H=\sum_{i=1}^3e_i\odot e_i^*.$
For any constant number $k_0$, set $K=k_0e^*_1\wedge e^*_2\wedge e^*_3$ and then $K$ is a $3$-cocycle on the sub-adjacent Lie algebra $\g^c$. Thus, $(\g^*[2],\g[1],\frkl^H_1,\frkl^H_2,\frkl^{H,K}_3)$ is a Lie $2$-algebra, where $\frkl^H_1=0$, and for all $x,y\in\g,~\xi\in\g^*$,
\begin{eqnarray*}
{\frkl_2^H(x,y)}=[x,y]_\g,\quad
{\frkl_2^H(x,\xi)}=L^*_{x}\xi, \quad
\frkl_3^{H,K}(e_1,e_2,e_3)=-4k_0e^*_1.
\end{eqnarray*}
}
\end{ex}

\section*{Appendix}

{\bf The proof of Proposition \ref{pro:C2L2}:}

\pf Let $\huaA=A_{-1}\oplus A_{-2}$ and $\huaB=B_{-1}\oplus B_{-2}$. Since the pairing $S$ is nondegenerate, $B_{-1}$ is isomorphic to $A_{-2}^*$, the dual bundle of $A_{-2}$, via $ \langle \alpha^2,X^2\rangle=S(\alpha^2,X^2)$ for all $X^2\in\Gamma(A_{-2}),~\alpha^2\in\Gamma(B_{-1})$, and $B_{-2}$ is isomorphic to $A_{-1}^*$, the dual bundle of $A_{-1}$, via $ \langle \alpha^1,X^1\rangle=S(\alpha^1,X^1)$ for all $X^1\in\Gamma(A_{-1}),~\alpha^1\in\Gamma(B_{-2})$. Under this isomorphism, { the graded symmetric bilinear form $S$ is given by \eqref{eq:naturalsymform}.}
\emptycomment{ \begin{equation*}
S(X^1+\alpha^2+X^2+\alpha^1,Y^1+\beta^2+Y^2+\beta^1)=\langle X^1,\beta^1 \rangle+\langle Y^1,\alpha^1 \rangle+\langle X^2,\beta^2 \rangle+\langle Y^2,\alpha^2 \rangle,
\end{equation*}
for all $X^1,Y^1\in\Gamma(A_{-1}),~X^2,Y^2\in\Gamma(A_{-2}),~\alpha^1,\beta^1\in\Gamma(B_{-2}),~\alpha^2,\beta^2\in\Gamma(B_{-1})$.}
By Proposition \ref{pro:Dirac structures}, both $\huaA$ and $\huaB$ are split Lie $2$-algebroids, and denoted by $(\huaA;l_1,l_2,l_3,a)$ and $(\huaB;\frkl_1,\frkl_2,\frkl_3,\frka)$.
\emptycomment{Their corresponding operations are given, respectively, by $$l_1=\partial\mid_{A_{-2}},\quad a=\rho\mid_{A_{-1}},\quad l_3=\Omega\mid_{A_{-1}},$$
 and
 $$\frkl_1=\partial\mid_{B_{-2}},\quad\frka=\rho\mid_{B_{-1}},\quad \frkl_3=\Omega\mid_{B_{-1}}.$$
}
 We  use $\delta$ and $\delta_*$ to denote their  differentials   respectively.

 {By   (iii) in Definition $\ref{defi:Courant-2 algebroid}$, we deduce that $\frkl_1=l^*_1.$ By  (ii) and (iv) in Definition $\ref{defi:Courant-2 algebroid}$, we deduce that the brackets between $\Gamma(\huaA)$ and $\Gamma(\huaB)$ are given by \eqref{eq:L2Bbracket0}.}
\emptycomment{
\begin{eqnarray*}
(X^1+\alpha^2)\diamond(Y^1+\beta^2)&=&l_2(X^1,Y^1)+L^1_{X^1}\beta^2-L^1_{Y^1}\alpha^2 +\frkl_2(\alpha^2,\beta^2) +\huaL^1_{\alpha^2}Y^1-\huaL^1_{\beta^2}X^1,\\
(X^1+\alpha^2)\diamond(X^2+\alpha^1)&=&l_2(X^1,X^2)+L^1_{X^1}\alpha^1+\iota_{X^2}\dM(\alpha^2) +\frkl_2(\alpha^2,\alpha^1) +\huaL^1_{\alpha^2}X^2+\iota_{\alpha^1}\dM_*(X^1), \\
(X^2+\alpha^1)\diamond(X^1+\alpha^2)&=&l_2(X^2,X^1)+L^2_{X^2}\alpha^2+\iota_{X^1}\dM(\alpha^1) +\frkl_2(\alpha^1,\alpha^2) +\huaL^2_{\alpha^1}X^1+\iota_{\alpha^2}\dM_*(X^2).
\end{eqnarray*}}
By  (v) in Definition $\ref{defi:Courant-2 algebroid}$, we deduce that the $(A_{-2}\oplus B_{-2})$-valued $3$-form $\Omega$ is given by \eqref{eq:L2B3-form}.
\emptycomment{\begin{eqnarray*}
\Omega(X^1+\alpha^2,Y^1+\beta^2,Z^1+\gamma^2)&=&l_3(X^1,Y^1,Z^1)+L^3_{X^1,Y^1}\gamma^2+L^3_{Y^1,Z^1}\alpha^2+L^3_{Z^1,X^1}\beta^2\\
&&+\frkl_3(\alpha^2,\beta^2,\gamma^2)+\huaL^3_{\alpha^2,\beta^2}Z^1+\huaL^3_{\beta^2,\gamma^2}X^1+\huaL^3_{\gamma^2,\alpha^2}Y^1,
\end{eqnarray*}
for all $X^1,Y^1,Z^1\in\Gamma(A_{-1}),~\alpha^2,\beta^2,\gamma^2\in\Gamma(B_{-1})$.}

Next, we will use the following two steps to show that \eqref{L2bialgebroid1} in Theorem \ref{thm:Lie2bi} holds.

{\bf Step 1:} We will show that
\begin{equation}\label{eq:Manin triple Lie2A 4}
{\delta}_*{[X^1,Y^1]}_S=-{[{\delta}_*(X^1),Y^1]}_S+{[X^1,{\delta}_*(Y^1)]}_S,\quad \forall X^1,Y^1\in\Gamma(A_{-1}).
\end{equation}

In fact, since for all $X^1\in\Gamma(A_{-1})$, $\bar{\delta}_*(X^1)=0$, we have
\begin{equation}\label{eq:Manin triple Lie2A 1}
\bar{\delta}_*{[X^1,Y^1]}_S=-{[\bar{\delta}_*(X^1),Y^1]}_S+{[X^1,\bar{\delta}_*(Y^1)]}_S.
\end{equation}
\emptycomment{By \eqref{eq:DFbracketleft}, we have
\begin{eqnarray*}
  0=Df\diamond X^1=(\delta f+\delta_\ast f)\diamond X^1=l_2(\delta_\ast f,X^1)+\huaL^2_{\delta f}X^1,
\end{eqnarray*}
which implies that
\begin{eqnarray}\label{eq:CLWX2-Dirac1}
  \huaL^2_{\delta f}X^1=-l_2(X^1,\delta_\ast f).
\end{eqnarray}
Similarly, we also have
\begin{eqnarray}\label{eq:CLWX2-Dirac2}
 L^2_{\delta_* f}\alpha^2&=&-\Schouten{\alpha^2,\delta f}.
\end{eqnarray}}
For all $X^1,Y^1\in\Gamma(A_{-1})$ and $\alpha^1\in\Gamma(B_{-2})$, by $\rm(e_1)$ in   Definition \ref{defi:2leibniz}, we have
\begin{eqnarray*}
  X^1\diamond (Y^1\diamond\alpha^1)-(X^1\diamond Y^1)\diamond\alpha^1-Y^1\diamond (X^1\diamond\alpha^1)=\Omega(X^1,Y^1,l^*_1(\alpha^1)).
\end{eqnarray*}
By \eqref{eq:jacobi-eq1}, this condition is equivalent to
\begin{eqnarray*}
  \iota_{L^1_{Y^1}\alpha^1}\dM_{\ast}X^1+l_2(X^1,\iota_{\alpha^1}\dM_{\ast}Y^1)-\iota_{\alpha^1}\dM_{\ast}l_2(X^1,Y^1)-\iota_{L^1_{Y^1}\alpha^1}\dM_{\ast}X^1-l_2(X^1,\iota_{\alpha^1}\dM_{\ast}Y^1)=0.
\end{eqnarray*}
By direct calculation, we have
\begin{eqnarray*}
  \langle\iota_{L^1_{Y^1}\alpha^1}\dM_{\ast}X^1,\alpha^2\rangle&=&\frka(\alpha^2)a(Y^1)\langle X^1,\alpha^1\rangle-\frka(\alpha^2)\langle \alpha^1,l_2(X^1,Y^1) \rangle-\langle X^1,\frkl_2(L^1_{Y^1}\alpha^1,\alpha^2)\rangle;\\
   \langle l_2(X^1,\iota_{\alpha^1}\dM_{\ast}Y^1),\alpha^2\rangle&=&a(X^1)\frka(\alpha^2)\langle Y^1,\alpha^1\rangle-a(X^1)\langle \frkl_2(\alpha^1,\alpha^2),Y^1 \rangle-\frka(L^1_{X^1}\alpha^2)\langle Y^1,\alpha^1\rangle\\
   &&+\langle Y^1,\frkl_2(L^1_{X^1}\alpha^2,\alpha^1)\rangle;\\
  \langle \iota_{\alpha^1}[\dM_{\ast}X^1,Y^1]_S, \alpha^2\rangle &=&a(Y^1)\dM_{\ast}X^1(\alpha^1,\alpha^2)-\dM_{\ast}X^1(L^1_{Y^1}\alpha^1,\alpha^2)-\dM_{\ast}X^1(\alpha^1,L^1_{Y^1}\alpha^2)\\
  &=&a(Y^1)\frka(\alpha^2)\langle X^1,\alpha^1\rangle-a(Y^1)\langle \frkl_2(\alpha^1,\alpha^2),X^1 \rangle-\frka(\alpha^2)a(Y^1)\langle X^1,\alpha^1\rangle\\
  &&+\frka(\alpha^2)\langle \alpha^1,l_2(Y^1,X^1) \rangle+\langle X^1,\frkl_2(L^1_{Y^1}\alpha^1,\alpha^2)\rangle-\frka(L^1_{Y^1}\alpha^2)\langle X^1,\alpha^1\rangle\\
  &&+\langle X^1,\frkl_2(L^1_{Y^1}\alpha^2,\alpha^1)\rangle.
\end{eqnarray*}
Then by the above formulas, we have
\begin{eqnarray*}
  0&=&\iota_{L^1_{Y^1}\alpha^1}\dM_{\ast}X^1+l_2(X^1,\iota_{\alpha^1}\dM_{\ast}Y^1)-\iota_{\alpha^1}\dM_{\ast}l_2(X^1,Y^1)-\iota_{L^1_{Y^1}\alpha^1}\dM_{\ast}X^1-l_2(X^1,\iota_{\alpha^1}\dM_{\ast}Y^1)\\
  &=&\iota_{\alpha^1}\Big(-{[\dM_*(X^1),Y^1]}_S+{[X^1,\dM_*(Y^1)]}_S-\dM_*{[X^1,Y^1]}_S\Big),
\end{eqnarray*}
which implies that
\begin{equation}\label{eq:Manin triple Lie2A 2}
\dM_*{[X^1,Y^1]}_S=-{[\dM_*(X^1),Y^1]}_S+{[X^1,\dM_*(Y^1)]}_S.
\end{equation}

For all $X^1\in\Gamma(A_{-1})$ and $\alpha^2,\beta^2,\gamma^2\in\Gamma(B_{-1})$, by $\rm(f)$ in   Definition \ref{defi:2leibniz}, we have
\begin{eqnarray*}
&&\alpha^2\diamond \Omega(\beta^2,\gamma^2,X^1)-\beta^2\diamond \Omega(\alpha^2,\gamma^2,X^1)+\gamma^2\diamond \Omega(\alpha^2,\beta^2,X^1)-\frkl_3(\alpha^2,\beta^2,\gamma^2)\diamond X^1\\
&&-\Omega(\frkl_2(\alpha^2,\beta^2),\gamma^2,X^1)-\Omega(\beta^2,\frkl_2(\alpha^2,\gamma^2),X^1)-\Omega(\beta^2,\gamma^2,\alpha^2\diamond X^1)\\
&&+\Omega(\alpha^2,\frkl_2(\beta^2,\gamma^2),X^1)+\Omega(\alpha^2,\gamma^2,\beta^2\diamond X^1)-\Omega(\alpha^2,\beta^2,\gamma^2\diamond X^1)=0,
\end{eqnarray*}
 which is equivalent to
\begin{eqnarray}
  \nonumber&&\iota_{\huaL^3_{\beta^2,\gamma^2}X^1}\dM\alpha^2-\iota_{\huaL^3_{\alpha^2,\gamma^2}X^1}\dM\beta^2+\iota_{\huaL^3_{\alpha^2,\beta^2}X^1}\dM\gamma^2-\iota_{X^1}\dM\frkl_3(\alpha^2,\beta^2,\gamma^2)\\
  \label{eq:Manin triple f1}&&+\frkl_3(\alpha^2,\beta^2,L^1_{X^1}\gamma^2)+\frkl_3(\gamma^2,\alpha^2,L^1_{X^1}\beta^2)+\frkl_3(\beta^2,\gamma^2,L^1_{X^1}\alpha^2)=0.
\end{eqnarray}
By direct calculation, we have
\begin{eqnarray}
 \label{eq:Manin triple f2} &&\langle\iota_{\huaL^3_{\beta^2,\gamma^2}X^1}\dM\alpha^2-\iota_{\huaL^3_{\alpha^2,\gamma^2}X^1}\dM\beta^2+\iota_{\huaL^3_{\beta^2,\beta^2}X^1}\dM\gamma^2-\iota_{X^1}\dM\frkl_3(\alpha^2,\beta^2,\gamma^2),Y^1\rangle\\
  \nonumber&=&\langle-\frkl_3(\gamma^2,\alpha^2,L^1_{Y^1}\beta^2)-\frkl_3(\alpha^2,\beta^2,L^1_{Y^1}\gamma^2)-\frkl_3(\beta^2,\gamma^2,L^1_{Y^1}\alpha^2),X^1\rangle\\
 \nonumber&&-a(X^1)\langle\frkl_3(\alpha^2,\beta^2,\gamma^2),Y^1 \rangle+a(Y^1)\langle\frkl_3(\alpha^2,\beta^2,\gamma^2),X^1 \rangle+\langle\frkl_3(\alpha^2,\beta^2,\gamma^2),l_2(X^1,Y^1)\rangle.
\end{eqnarray}
On the other hand, we have
\begin{eqnarray*}
  [\hat{\delta}_{\ast}X^1,Y^1]_S(\alpha^2,\beta^2,\gamma^2)&=&-a(Y^1)\langle\frkl_3(\alpha^2,\beta^2,\gamma^2),X^1 \rangle+\langle\frkl_3(\gamma^2,\alpha^2,L^1_{Y^1}\beta^2)+\frkl_3(\alpha^2,\beta^2,L^1_{Y^1}\gamma^2)\\
  &&+\frkl_3(\beta^2,\gamma^2,L^1_{Y^1}\alpha^2),X^1\rangle,\\
  \hat{\delta}_{\ast}[X^1,Y^1]_S(\alpha^2,\beta^2,\gamma^2)&=&-\langle\frkl_3(\alpha^2,\beta^2,\gamma^2),l_2(X^1,Y^1)\rangle.
\end{eqnarray*}
Then by \eqref{eq:Manin triple f1} and \eqref{eq:Manin triple f2}, we have
\begin{eqnarray*}
  &&\Big(-{[\hat{\delta}_*(X^1),Y^1]}_S+{[X^1,\hat{\delta}_*(Y^1)]}_S-\hat{\delta}_*{[X^1,Y^1]}_S\Big)(\alpha^2,\beta^2,\gamma^2)\\
  &=&a(Y^1)\langle\frkl_3(\alpha^2,\beta^2,\gamma^2),X^1 \rangle-\langle\frkl_3(\gamma^2,\alpha^2,L^1_{Y^1}\beta^2)+\frkl_3(\alpha^2,\beta^2,L^1_{Y^1}\gamma^2)+\frkl_3(\beta^2,\gamma^2,L^1_{Y^1}\alpha^2),X^1\rangle\\
  &&-a(X^1)\langle\frkl_3(\alpha^2,\beta^2,\gamma^2),Y^1 \rangle+\langle\frkl_3(\gamma^2,\alpha^2,L^1_{X^1}\beta^2)+\frkl_3(\alpha^2,\beta^2,L^1_{X^1}\gamma^2)+\frkl_3(\beta^2,\gamma^2,L^1_{X^1}\alpha^2),Y^1\rangle\\
  &&+\langle\frkl_3(\alpha^2,\beta^2,\gamma^2),l_2(X^1,Y^1)\rangle=0,
\end{eqnarray*}
which implies that
\begin{equation}\label{eq:Manin triple Lie2A 3}
\hat{\delta}_*{[X^1,Y^1]}_S=-{[\hat{\delta}_*(X^1),Y^1]}_S+{[X^1,\hat{\delta}_*(Y^1)]}_S.
\end{equation}

Thus, by \eqref{eq:Manin triple Lie2A 1}, \eqref{eq:Manin triple Lie2A 2} and \eqref{eq:Manin triple Lie2A 3},  \eqref{eq:Manin triple Lie2A 4} follows immediately. \vspace{3mm}

{\bf Step 2:} We will show that
\begin{equation}\label{eq:Manin triple Lie2A 42}
{\delta}_*{[X^1,Y^2]}_S=-{[{\delta}_*(X^1),Y^2]}_S+{[X^1,{\delta}_*(Y^2)]}_S,\quad \forall X^1\in\Gamma(A_{-1}), Y^2\in\Gamma(A_{-2}).
\end{equation}

For $X^1\in\Gamma(A_{-1}),~Y^2\in\Gamma(A_{-2})$,  by $\rm(a)$ in Definition \ref{defi:2leibniz}, we have
\begin{eqnarray*}
  \langle\bar{\delta}_*{[X^1,Y^2]}_S+{[\bar{\delta}_*(X^1),Y^2]}_S-{[X^1,\bar{\delta}_*(Y^2)]}_S,\alpha^1\rangle
  =\langle-l_1(l_2(X^1,Y^2)+l_2(X^1,l_1(Y^2)),\alpha^1\rangle
  =0,
\end{eqnarray*}
which implies that
\begin{equation}\label{eq:Manin triple Lie2A 12}
\bar{\delta}_*{[X^1,Y^2]}_S=-{[\bar{\delta}_*(X^1),Y^2]}_S+{[X^1,\bar{\delta}_*(Y^2)]}_S.
\end{equation}

For all $X^1\in\Gamma(A_{-1}),~Y^2\in\Gamma(A_{-2})$ and $\alpha^2\in\Gamma(B_{-1})$, by $\rm(e_1)$ in Definition \ref{defi:2leibniz}, we have
\begin{eqnarray*}
  X^1\diamond (\alpha^2\diamond Y^2)-(X^1\diamond \alpha^2)\diamond Y^2-\alpha^2\diamond (X^1\diamond Y^2)=\Omega(X^1,\alpha^2,l_1(Y^2)).
\end{eqnarray*}
By \eqref{eq:interi-eq3}, this condition is equivalent to
\begin{eqnarray*}
l_2(X^1,\huaL^1_{\alpha^2}Y^2)+\iota_{\iota_{Y^2}\dM\alpha^2}\dM_{\ast}X^1-\huaL^1_{L^1_{X^1}\alpha^2}Y^2  +l_2(\huaL^1_{\alpha^2}X^1,Y^2)-\huaL^1_{\alpha^2}l_2(X^1,Y^2)=0.
\end{eqnarray*}
By direct calculation, on the one hand, we have
\begin{eqnarray*}
  \langle l_2(X^1,\huaL^1_{\alpha^2}Y^2),\beta^2\rangle&=&a(X^1)\frka(\alpha^2)\langle Y^2,\beta^2\rangle-a(X^1)\langle Y^2,\frkl_2(\alpha^2,\beta^2)\rangle-\frka(\alpha^2)a(X^1)\langle Y^2,\beta^2\rangle\\
  &&+\frka(\alpha^2)\langle \beta^2,l_2(X^1,Y^2)\rangle+\langle Y^2,\frkl_2(\alpha^2,L^1_{X^1}\beta^2)\rangle;\\
  \langle\iota_{\iota_{Y^2}\dM\alpha^2}\dM_{\ast}X^1,\beta^2\rangle&=&a(\huaL^1_{\beta^2}X^1)\langle Y^2,\alpha^2\rangle-\frka(\beta^2)\langle \alpha^2,l_2(X^1,Y^2)\rangle-\langle X^1,\frkl_2(\beta^2,L^2_{Y^2}\alpha^2)\rangle;\\
  \langle\huaL^1_{L^1_{X^1}\alpha^2}Y^2,\beta^2\rangle&=&\frka(L^1_{X^1}\alpha^2)\langle Y^2,\beta^2\rangle-\langle Y^2,\frkl_2(L^1_{X^1}\alpha^2,\beta^2)\rangle;\\
 \langle l_2(\huaL^1_{\alpha^2}X^1,Y^2),\beta^2\rangle&=&\frka(\alpha^2)\langle \beta^2,l_2(X^1,Y^2)\rangle+\langle X^1,\frkl_2(\alpha^2,L^2_{Y^2}\beta^2)\rangle;\\
\langle \huaL^1_{\alpha^2}l_2(X^1,Y^2),\beta^2\rangle&=&\frka(\alpha^2)\langle \beta^2,l_2(X^1,Y^2)\rangle-\langle l_2(X^1,Y^2),\frkl_2(\alpha^2,\beta^2)\rangle.
\end{eqnarray*}
On the other hand, we have
\begin{eqnarray*}
 \langle\iota_{\alpha^2}\dM_*{[X^1,Y^2]}_S,\beta^2\rangle&=&\frka(\alpha^2)\langle \beta^2,l_2(X^1,Y^2)\rangle-\frka(\beta^2)\langle \alpha^2,l_2(X^1,Y^2)\rangle-\langle l_2(X^1,Y^2),\frkl_2(\alpha^2,\beta^2)\rangle;\\
\langle\iota_{\alpha^2} {[\dM_*(X^1),Y^2]}_S,\beta^2\rangle&=&-\dM_{\ast}X^1(L^2_{Y^2}\alpha^2,\beta^2)+\dM_{\ast}X^1(L^2_{Y^2}\beta^2,\alpha^2)\\
&=&\frka(\beta^2)\langle \alpha^2,l_2(X^1,Y^2)\rangle+\langle X^1,\frkl_2(\beta^2,L^2_{Y^2}\alpha^2)\rangle-\frka(\alpha^2)\langle \beta^2,l_2(X^1,Y^2)\rangle\\
&&-\langle X^1,\frkl_2(\alpha^2,L^2_{Y^2}\beta^2)\rangle;\\
\langle\iota_{\alpha^2}{[X^1,\dM_*(Y^2)]}_S,\beta^2\rangle&=&a(X^1)\frka(\alpha^2)\langle Y^2,\beta^2\rangle-a(X^1)\frka(\beta^2)\langle Y^2,\alpha^2\rangle-a(X^1)\langle Y^2,\frkl_2(\alpha^2,\beta^2)\rangle\\
&&-\frka(L^1_{X^1}\alpha^2)\langle Y^2,\beta^2\rangle+\frka(\beta^2)a(X^1)\langle Y^2,\alpha^2\rangle-\frka(\beta^2)\langle \alpha^2,l_2(X^1,Y^2)\rangle\\
&&+\langle Y^2,\frkl_2(L^1_{X^1}\alpha^2,\beta^2)\rangle-\frka(\alpha^2)a(X^1)\langle Y^2,\beta^2\rangle+\frka(\alpha^2)\langle \beta^2,l_2(X^1,Y^2)\rangle\\
&&+\frka(L^1_{X^1}\beta^2)\langle Y^2,\alpha^2\rangle+\langle Y^2,\frkl_2(\alpha^2,L^1_{X^1}\beta^2)\rangle.
\end{eqnarray*}
Then by the above formulas and $\rho(X^1\diamond \alpha^2)=[\rho(X^1),\rho(\alpha^2)]$, we have
\begin{eqnarray*}
  0&=&l_2(X^1,\huaL^1_{\alpha^2}Y^2)+\iota_{\iota_{Y^2}\dM\alpha^2}\dM_{\ast}X^1-\huaL^1_{L^1_{X^1}\alpha^2}Y^2  +l_2(\huaL^1_{\alpha^2}X^1,Y^2)-\huaL^1_{\alpha^2}l_2(X^1,Y^2)\\
  &=&\iota_{\alpha^2}\Big(-{[\dM_*(X^1),Y^2]}_S+{[X^1,\dM_*(Y^2)]}_S-\dM_*{[X^1,Y^2]}_S\Big),
\end{eqnarray*}
which implies that
\begin{equation}\label{eq:Manin triple Lie2A 22}
\dM_*{[X^1,Y^2]}_S=-{[\dM_*(X^1),Y^2]}_S+{[X^1,\dM_*(Y^2)]}_S.
\end{equation}

By a direct calculation, we have
 \begin{eqnarray*}
{[\hat{\delta}_*(X^1),Y^2]}_S(\alpha^2,\beta^2,\gamma^2)={[X^1,\hat{\delta}_*(Y^2)]}_S(\alpha^2,\beta^2,\gamma^2)=\hat{\delta}_*{[X^1,Y^2]}_S(\alpha^2,\beta^2,\gamma^2)=0,
\end{eqnarray*}
which implies that
\begin{equation}\label{eq:Manin triple Lie2A 32}
\hat{\delta}_*{[X^1,Y^2]}_S=-{[\hat{\delta}_*(X^1),Y^2]}_S+{[X^1,\hat{\delta}_*(Y^2)]}_S.
\end{equation}
Thus, by \eqref{eq:Manin triple Lie2A 12}, \eqref{eq:Manin triple Lie2A 22} and \eqref{eq:Manin triple Lie2A 32},  \eqref{eq:Manin triple Lie2A 42} follows immediately.

By \eqref{eq:Manin triple Lie2A 4} and \eqref{eq:Manin triple Lie2A 42}, for all $X,Y\in\Gamma(\huaA[-3])$, we have
$$\delta_*{[X,Y]}_S=-{[\delta_*(X),Y]}_S+(-1)^{|X|}{[X,\delta_*(Y)]}_S,$$
which implies that \eqref{L2bialgebroid1} in Theorem \ref{thm:Lie2bi} holds.

Similarly, we can show that \eqref{L2bialgebroid2} also holds. Therefore, $(\huaA,\huaB)$ is a split Lie 2-bialgebroid.\qed

Jiefeng Liu

 School of Mathematics and Statistics, Northeast Normal University, Changchun 130024, Jilin, China

 Email:liujf12@126.com\\

 Yunhe Sheng

 Department of Mathematics, Jilin University, Changchun 130012, Jilin, China

Email:shengyh@jlu.edu.cn


{\footnotesize
\begin{thebibliography}{999}

\bibitem[AKSZ]{aksz}
M. Alexandrov, M. Kontsevich, A. Schwarz and O. Zaboronsky, The geometry of the master equation and topological quantum field theory, {\em Internat. J. Modern Phys. A} 12 (1997), no. 7, 1405-1429.

\bibitem[AP]{ammardefiLeibnizalgebra}
 M. Ammar and N. Poncin, Coalgebraic Approach to the Loday Infinity Category, Stem
Differential for $2n$-ary Graded and Homotopy Algebras, \emph{Ann.
Inst. Fourier (Grenoble).} 60 (1) (2010),  355-387.

\bibitem[AM]{AM}
A. Alekseev and E. Meinrenken, Dirac structures and Dixmier-Douady bundles, \emph{Int. Math. Res. Not.} IMRN  2012,  no. 4, 904–956.

\bibitem[BC]{baez:2algebras}
 J. C. Baez and A. S. Crans, Higher-Dimensional Algebra VI: Lie
 2-Algebras, \emph{Theory  Appl.  Categ.} 12 (2004),
 492-528.

 \bibitem[BSZ]{BSZ}
C.  Bai, Y.  Sheng and C.  Zhu, Lie 2-bialgebras, \emph{Comm. Math. Phys.} 320 (2013), no. 1, 149-172.

\bibitem[BV]{BV}
D. Bashkirov and A. A. Voronov, On homotopy Lie bialgebroids,  arXiv:1612.02026.

\bibitem[BP]{BP}
G. Bonavolont\`a and N. Poncin, On the category of Lie $n$-algebroids, \emph{J. Geom. Phys.} 73 (2013), 70-90.

 \bibitem[BCSX]{derivedPoisson}
R. Bandiera, Z. Chen, M. Stiénon and P  Xu,
Shifted derived Poisson manifolds associated with Lie pairs, \emph{Comm. Math. Phys.} (2019). https://doi.org/10.1007/s00220-019-03457-w.

\bibitem[Bru]{Bruce}
A. Bruce,  From $L_\infty$-algebroids to higher Schouten/Poisson structures, \emph{Rep. Math. Phys.} 67 (2011), no. 2, 157-177.

\bibitem[Bur]{Bur}
H. Bursztyn,  A brief introduction to Dirac manifolds,  Geometric and topological methods for quantum field theory, 4–38, Cambridge Univ. Press, Cambridge, 2013.

\bibitem[BCG]{BCG}
H. Bursztyn, G. Cavalcanti and M. Gualtieri,  Reduction of Courant algebroids and generalized complex structures, \emph{Adv. Math.}  211  (2007),  no. 2, 726–765.

\bibitem[BC]{BC}
H. Bursztyn and M. Crainic,  Dirac geometry, quasi-Poisson actions and D/G -valued moment maps, \emph{J. Diff. Geom.}  82  (2009),  no. 3, 501–566.

\bibitem[CGM]{CGM}
A. Cabrera, M. Gualtieri and E. Meinrenken, Dirac geometry of the holonomy fibration, \emph{Comm. Math. Phys.}  355  (2017),  no. 3, 865-904.

\bibitem[CGdS]{CGdS}
 J. F. Cariñena, J. Grabowski, J. de Lucas and C. Sardón, Dirac-Lie systems and Schwarzian equations, \emph{J. Diff. Equa.}  257  (2014),  no. 7, 2303–2340.

\bibitem[CF]{CF} A. S. Cattaneo and G. Felder,  Relative formality theorem and quantisation of coisotropic submanifolds, \emph{Adv. Math.} {\bf 208} (2007), no. 2, 521-548.

\bibitem[CSX1]{CSX1}Z. Chen, M. Sti\'enon and P. Xu, Weak Lie 2-bialgebras, \emph{J. Geom. Phys.} 68 (2013), 59-68.


\bibitem[Cou]{CourantDirac}%
T. Courant, Dirac manifolds,\emph{ Trans. Amer. Math. Soc.} 319 (1990), 631-661.

\bibitem[Get]{getzler:higher-derived}
E. Getzler, \newblock {Higher derived brackets}, \newblock arXiv:1010.5859v1.

\bibitem[Gua]{gualtieri}
M. Gualtieri,  Generalized complex geometry, \emph{Ann. of Math.} (2) 174 (2011), no. 1, 75-123.

\bibitem[Hit]{hitchin}
N. J. Hitchin, Generalized Calabi-Yau manifolds, \emph{Q. J. Math.} 54 (2003), no. 3,
281-308.


\bibitem[IU]{Ikeda}
N. Ikeda and K. Uchino, QP-structures of degree 3 and 4D topological field theory, \emph{Comm. Math. Phys.} 303 (2011), no. 2, 317-330.



\bibitem[Jotz13]{Jotz13} M. Jotz and T. Ratiu,  Dirac optimal reduction, \emph{Int. Math. Res. Not.} IMRN  2013,  no. 1, 84–155.

\bibitem[Jot15]{Jotz2}
M. Jotz Lean, N-manifolds of degree $2$ and metric double vector bundles, arXiv:1504.00880.

\bibitem[Jot19]{Jotz3}
M. Jotz Lean, Lie $2$-algebroids and matched pairs of $2$-representations - a geometric approach,  \emph{Pacific J. Math.}  301 (1) (2019),  143-188.

\bibitem[KV]{VoronovHigherP}
H. M. Khudaverdian and Th. Voronov, Higher Poisson brackets and differential forms, \emph{Geometric methods in physics,} 203-215, AIP Conf. Proc., 1079, \emph{Amer. Inst. Phys., Melville, NY,} 2008.

\bibitem[KS96]{Kosmann-Schwarzbach}
 Y. Kosmann-Schwarzbach, From Poisson algebras to Gerstenhaber algebras, {\em Ann. Inst. Fourier (Grenoble) } 46 (1996), no. 5, 1243-1274.

\bibitem[KS05]{Kosmann-Schwarzbach05}
Y. Kosmann-Schwarzbach, Quasi, twisted, and all that. . . in Poisson geometry and Lie algebroid theory. {\em The breadth
of symplectic and Poisson geometry}, 363-389, Progr. Math., 232, {\em Birkhäuser Boston, Boston, MA,} 2005.

  \bibitem[KS13]{Schwarzbach4}
Y. Kosmann-Schwarzbach, Courant algebroids. A short history, \emph{SIGMA Symmetry Integrability Geom. Methods Appl.} 9 (2013), Paper 014, 8 pp.


\bibitem[Kra]{olga}
O.~Kravchenko,
\newblock Strongly homotopy {L}ie bialgebras and {L}ie
quasi-bialgebras,
\newblock {\em Lett. Math. Phys.} 81 (1) (2007), 19-40.

\bibitem[LM95]{LadaMartin}
T. Lada and M. Markl,
\newblock Strongly homotopy {L}ie algebras,
\newblock {\em Comm. Algebra} 23 (6) (1995), 2147-2161.

\bibitem[LS]{stasheff:introductionSHLA} T. Lada and J. Stasheff, Introduction to sh Lie algebras for
physicists, \emph{Int. J. Theor. Phys.} 32 (7) (1993), 1087-1103.

\bibitem[LSX]{LSX}
H. Lang, Y. Sheng and X. Xu, Strong homotopy Lie algebras, homotopy Poisson manifolds and Courant algebroids, \emph{Lett. Math. Phys.} 107 (5) (2017.05), 861-885.

\bibitem[LiM09]{DLB}
D. Li-Bland and E. Meinrenken, Courant algebroid and Poisson
geometry, \emph{Int. Math. Res. Not.} 11 (2009),
2106-2145.

\bibitem[LiM14]{LiM}
D. Li-Bland and E. Meinrenken, Dirac Lie groups, \emph{Asian J. Math.}  18  (2014),  no. 5, 779–815.

\bibitem[Liu]{LiuDirac}
Z.-J. Liu, Some remarks on Dirac structures and Poisson reduction,
{\em Poisson Geometry}, Banach Center Publications, vol 51, 2000,
165-173.

\bibitem[LSBC]{LiuShengBaiChen}
J. Liu, Y. Sheng, C.  Bai and Z.   Chen,    Left-symmetric algebroids,   \emph{ Math. Nachr.  }  289 (14-15) (2016), 1893-1908.

\bibitem[LS]{LiuSheng}
J. Liu and Y. Sheng, QP-structures of degree 3 and $\LWX$ 2-algebroids,  \emph{J. Symplectic Geom.} 17 (6) (2019), 1853-1891.

\bibitem[LWX97]{lwx}
Z. Liu, A. Weinstein and P. Xu,
\newblock Manin triples for {L}ie bialgebroids,
\newblock {\em J. Diff. Geom.} 45 (3) (1997), 547-574.

 \bibitem[Liv]{livernet}
M. Livernet, Homologie des alg$\rm\grave{e}$bres stables de matrices sur une $A_\infty$-alg$\rm\grave{e}$bre, \emph{C. R. Acad. Sci. Paris S$\rm\acute{e}$r. I
Math.} 329 (2) (1999), 113-116.

\bibitem[Mac]{General theory of Lie groupoid and Lie algebroid}
K. C. H. Mackenzie, General theory of Lie groupoids and Lie
algebroids.  \emph{Lecture Note Series, $213$. London Mathematical
Society.} Cambridge University Press, Cambridge, 2005.


\bibitem[MX]{MackenzieX:1994}
 K. Mackenzie and P. Xu, Lie bialgebroids and Poisson groupoids,
\emph{Duke Math. J.} 73 (2) (1994), 415-452.

\bibitem[Meh]{Mehta}
R. A. Mehta, On homotopy Poisson actions and reduction of symplectic $Q$-manifolds, \emph{Diff. Geom. Appl.}  {\bf 29} (2011), no. 3, 319-328.


\bibitem[MZ]{MZ}
R. A. Mehta and M. Zambon,   $L_\infty$-algebra actions, \emph{Diff. Geom. Appl.} 30 (2012), no. 6, 576-587.

\bibitem[Mei18]{Mei18}
E. Meinrenken,  Poisson geometry from a Dirac perspective, \emph{Lett. Math. Phys.}  108  (2018),  no. 3, 447–498.

\bibitem[Mei17]{Mei17}
E. Meinrenken, Dirac actions and Lu's Lie algebroid, \emph{Transform. Groups}  22  (2017),  no. 4, 1081–1124.

\bibitem[Ngu01]{Boyom1}
M. Nguiffo Boyom,  Cohomology of Koszul-Vinberg algebroids and Poisson manifolds. I, \emph{Banach Center Publ.} 54 (2001), 99-110.

\bibitem[Ngu05]{Boyom2}
M. Nguiffo Boyom, KV-cohomology of Koszul-Vinberg algebroids and Poisson manifolds, \emph{Internat. J. Math.} 16 (2005), no. 9, 1033-1061.


\bibitem[Roy]{Roytenbergphdthesis}
D. Roytenberg,
\newblock {Courant algebroids, derived brackets and even symplectic
  supermanifolds},
\newblock PhD thesis, UC Berkeley, 1999, arXiv:math.DG/9910078.

\bibitem[Roy02A]{roytwist}
D. Roytenberg, Quasi-Lie bialgebroids and twisted Poisson manifolds, \emph{Lett. Math. Phys.} 61 (2002), 123-137.

\bibitem[Roy02]{royt}
D. Roytenberg, On the structure of graded symplectic supermanifolds and {C}ourant
  algebroids, In {\em Quantization, Poisson brackets and beyond (Manchester,
  2001)}, volume 315 of {\em Contemp. Math.} pages 169-185. Amer. Math. Soc.,  Providence, RI, 2002.


\bibitem[Roy07B]{RoyCF}
 D.  Roytenberg, AKSZ-BV formalism and Courant algebroid-induced topological field theories, \emph{Lett. Math. Phys.} 79 (2007), no. 2, 143-159.

\bibitem[RW98]{rw}
D. Roytenberg and A. Weinstein,
\newblock Courant algebroids and strongly homotopy {L}ie algebras,
\newblock {\em Lett. Math. Phys.} 46 (1) (1998), 81-93.

\bibitem[Sh]{sheng}
Y. Sheng, The first Pontryagin class of a quadratic Lie 2-algebroid, \emph{Comm. Math. Phys.}  362 (2) (2018),  689-716.

\bibitem[SL]{Leibniz2al}Y. Sheng and Z. Liu, Leibniz $2$-algebras and twisted Courant algebroids, \emph{ Comm. Algebra } 41 (2013), no. 5, 1929-1953.

\bibitem[SZ12]{sz:intsemi}
Y. Sheng and C. Zhu, Integration of semidirect product Lie 2-algebras, \emph{Int. J. Geom. Methods Mod. Phys.} Vol. 9, No. 5 (2012), 1250043.


\bibitem[SZ17]{sz} Y. Sheng and C. Zhu, Higher extensions of Lie algebroids, \emph{Comm. Contemp. Math.}  19 (3) (2017), 1650034, 41 pages.


\bibitem[Sta]{Stasheff1} M. Schlessinger and J. Stasheff, The Lie algebra
structure of tangent cohomology and deformation theory, \emph{J.
Pure Appl. Algebra} 38 (1985), 313-322.



  \bibitem[Vor05]{Voronov1}
T. Voronov, Higher derived brackets and homotopy algebras, \emph{J. Pure Appl. Algebra} 202 (2005), no. 1-3, 133-153.

\bibitem[Vor10]{Voronov:2010halgd}
T. Voronov,
\newblock {Q-manifolds and Higher Analogs of Lie Algebroids},
\newblock XXIX Workshop on Geometric Methods in Physics. AIP CP 1307, pp.
  191-202, Amer. Inst. Phys., Melville, NY, 2010.

\bibitem[Zam]{zambon:l-infty}
M. Zambon, \newblock {$L_\infty$-algebras and higher analogues of Dirac structures and Courant algebroids}, \newblock \emph{J. Symplectic Geom.} 10 (2012), no. 4, 563-599.

\end{thebibliography}
}
\end{document}